\documentclass[11pt]{article}
\usepackage[letterpaper,margin=1in]{geometry}
\usepackage[hidelinks]{hyperref}
\usepackage{amsmath,amssymb,amsthm,mathtools}
\usepackage{mathrsfs}
\usepackage[T1]{fontenc}
\usepackage{lmodern}
\usepackage{enumitem}

\numberwithin{equation}{section}
\usepackage{microtype}

\newtheorem{theorem}{Theorem}[section]
\newtheorem{proposition}[theorem]{Proposition}
\newtheorem{lemma}[theorem]{Lemma}

\newtheorem*{remark*}{Remark}
\newtheorem{corollary}[theorem]{Corollary}

\newtheorem{assumption}{Assumption}[section]
\newcommand{\R}{\mathbb R}
\newcommand{\Mplus}{\mathcal M_+}
\newcommand{\dd}{\,\mathrm d}
\title{Large deviation principle of error of Euler--Maruyama  method  for stochastic differential equations}
\author{Diancong Jin\\[4pt]
  \parbox[t]{0.94\textwidth}{\centering\small
    School of Mathematics and Statistics,\\
    Huazhong University of Science and Technology, Wuhan 430074, China;\\[3pt]
    Hubei Key Laboratory of Engineering Modeling and Scientific Computing,\\
    Huazhong University of Science and Technology, Wuhan 430074, China.\\[4pt]
    \href{mailto:jindc@hust.edu.cn}{\texttt{jindc@hust.edu.cn}}}}
\date{}
\begin{document}
\maketitle

\begin{abstract}	
We study large deviations for the error of the Euler–Maruyama method applied to a one-dimensional stochastic differential equation with $C_b^3$ coefficients and a uniformly elliptic diffusion coefficient. We prove that the error at the terminal time, without rescaling, satisfies a large deviation principle (LDP) at speed $n$, where $n$ is the number of steps, with a good rate function that vanishes only at the origin. In particular, the probability that the error exceeds a fixed level decays exponentially in $n$. Two effects determine the rate function. A single Brownian increment of fixed nonzero size has
finite cost at speed $n$, and the cumulative effect of smaller Brownian increments depends on the
states at which their squares accumulate. To capture both effects, we introduce, for a cutoff $\varepsilon\in(0,1]$, an empirical enhanced control that records the state-dependent quadratic variation of the small increments together with the time, state and size of each increment exceeding $\varepsilon$. We prove an LDP for the empirical enhanced control for every fixed $\varepsilon$, and deduce the LDP for the error by the contraction principle and an exponentially good approximation as $\varepsilon\to0$.

\end{abstract}
\medskip
\noindent\textbf{Keywords:} Euler--Maruyama scheme, terminal error,
large deviations, empirical enhanced control, exponentially good approximation.
\section{Introduction}
\label{sec:introduction}
\subsection{Setting and background}
Consider the following scalar stochastic differential equation (SDE) 
\begin{align}
	\dd X_t=\sigma(X_t)\dd W_t+b(X_t)\dd t,~t\in(0,1],\qquad X_0=x_0, \label{eq:tag-1}
\end{align}
where $W=\{W_t,~t\in[0,1]\}$ is a standard Brownian motion defined on a complete probability space \((\Omega,\mathcal F,\mathbb P)\). For \(t_k=k/n\) and
\(\Delta_kW=W_{t_{k+1}}-W_{t_k}\), consider the Euler--Maruyama method for \eqref{eq:tag-1}:
\begin{align}
X^n_{k+1}&=X^n_k+\sigma(X^n_k)\Delta_kW+\frac1n b(X^n_k),
\qquad X^n_0=x_0. \label{eq:tag-2}
\end{align}
Define the terminal error $U^n_1$ by
\begin{align}
U^n_1&=X_1-X^n_n. \label{eq:tag-3}
\end{align}
In this work, we study how fast \(\mathbb P\{|U_1^n|\ge\delta\}\) decays as \(n\to\infty\) for fixed \(\delta>0\). More precisely, we investigate whether the LDP holds for the terminal error $U^n_1$ at speed $n$.

We impose the following assumption.
\begin{assumption}[Uniformly elliptic coefficients]\label{ass:coeff}
\[
\sigma,b\in C_b^3(\R),\qquad
0<\sigma_{\mathrm{min}}\leq\sigma(x)\leq\sigma_{\mathrm{max}}<\infty,
\qquad x\in\R.
\]
\end{assumption}
\noindent Here \(C_b^3(\R)\) means that the functions and their derivatives up to
order three are bounded and continuous. The two bounds on \(\sigma\) make the change of coordinate in
\eqref{eq:tag-7} invertible on the whole real line.

Euler~\cite{Euler1768} approximated solutions of ordinary differential
equations by taking each step with the slope at its left endpoint.
Maruyama's stochastic discretization~\cite{Maruyama} likewise holds the
coefficients fixed over each step. In~\eqref{eq:tag-2}, the drift
\(b(X_k^n)\) is multiplied by the step size \(1/n\) and the diffusion
coefficient \(\sigma(X_k^n)\) by the Brownian increment \(\Delta_k W\).
Pardoux and Talay~\cite{PardouxTalay} studied pathwise and mean-square
convergence of stochastic discretization schemes and compared them in
simulations. Under smoothness and growth conditions on the coefficients
and the test function \(f\), Talay and Tubaro~\cite{TalayTubaro} obtained
an expansion of the weak Euler error
\(\mathbb E[f(X_n^n)]-\mathbb E[f(X_1)]\) in powers of the step size.

There has been some work concerning the probabilistic asymptotic behavior of errors of the Euler method.   For equations driven by Brownian motion, Kurtz and
Protter~\cite{KurtzProtterEuler} studied the limiting distribution of
the normalized Euler error. Jacod and Protter~\cite{JacodProtter}
extended this analysis to equations driven by more general
semimartingales. Cl\'ement~\cite{Clement} studied pseudo-moderate
deviations of the Euler error process for one-dimensional diffusions.
For Brownian diffusions, the limiting distributions
in~\cite{KurtzProtterEuler,JacodProtter} describe fluctuations of the
Euler error on the \(n^{-1/2}\) scale. The pseudo-moderate deviation
estimates in~\cite{Clement} give exponential probability bounds for
deviations on the larger scale \(a_n/\sqrt{n}\), where \(a_n\to\infty\)
and  \(a_n=o(n^{1/6})\).
Here we ask how fast probabilities such as
\(\mathbb P\{|U_1^n|\geq\delta\}\) decay as \(n\to\infty\), with both
the terminal time and the threshold \(\delta>0\) fixed.
For linear SDEs with additive noise, Hong, Liang and
Sheng~\cite[Theorem~2.3(iii)]{HLS2023} considered a class of
schemes that includes Euler--Maruyama. They proved that each component
of the unscaled terminal error satisfies a full LDP at speed \(n^2\),
where \(n\) is the number of time steps. However, to the best of our knowledge, no results on the LDP of numerical errors for nonlinear SDEs are available.

\subsection{Main result and contributions}
Our main result is as follows.
\begin{theorem}[Main result]\label{thm:main}
Under Assumption~\ref{ass:coeff}, the terminal errors \(U^n_1\) satisfy a
full LDP on \(\R\) at speed \(n\) with good rate function \(J_{\rm un}\)
given by the variational formula \eqref{eq:tag-6}.
Moreover, \(J_{\rm un}(e)=0\) if and only if \(e=0\).
\end{theorem}

For every fixed \(\delta>0\), Theorem \ref{thm:main} gives an exponential upper bound for \(\mathbb P\{|U_1^n|\ge\delta\}\) at speed \(n\). More precisely, it holds that
\[
\limsup_{n\to\infty}\frac1n\log\mathbb P\{|U_1^n|\ge\delta\}
\le -c_\delta,\qquad c_\delta=\inf_{|e|\ge\delta}J_{\rm un}(e)\in(0,\infty].
\]
Proposition~\ref{prop:rate} also shows
that \(J_{\rm un}=\lim_{\varepsilon\downarrow0}J_\varepsilon\),
where \(J_\varepsilon\) is the rate function for the terminal error
approximation with cutoff \(\varepsilon\).

As a consistency check, suppose that \(b(x)\equiv b_0\) and
\(\sigma(x)\equiv\sigma_0>0\). Then it holds \(U_1^n=0\). In this case the Lamperti transform is simply
\(\Theta(x)=x/\sigma_0\), and for every state label \(\xi\) we have
\[
\alpha\equiv0,\qquad
\beta\equiv\gamma\equiv b_0/\sigma_0,\qquad
\bar{\mathfrak j}(\xi,\delta)=\delta,\qquad
\bar{\mathfrak g}(\xi,q)=\sigma_0q.
\]
Thus it follows from \eqref{eq:controlled-error} that \(d_{A,B,\mathfrak m}(1)=0\) for every control, so its
terminal error is \(\bar{\mathfrak g}(\xi_f,0)=0\). Taking \(B=A\) and
an empty mark list gives a control of zero cost. Hence \eqref{eq:tag-6}
yields \(J_{\rm un}(0)=0\) and \(J_{\rm un}(e)=+\infty\) for every \(e\neq0\).

For every open set \(O\subset\R\) and every closed set \(F\subset\R\),
the full LDP asserts that
\[
\begin{aligned}
\liminf_{n\to\infty}\frac1n\log\mathbb P\{U_1^n\in O\}
   &\geq-\inf_{e\in O}J_{\rm un}(e),\\
\limsup_{n\to\infty}\frac1n\log\mathbb P\{U_1^n\in F\}
   &\leq-\inf_{e\in F}J_{\rm un}(e).
\end{aligned}
\]
For instance, if an interval \(D\) satisfies
\(\inf_{D^\circ}J_{\rm un}=\inf_{\overline D}J_{\rm un}=c<\infty\), then
\(\mathbb P\{U_1^n\in D\}=\exp(-nc+o(n))\).

Let us state our main contributions as follows.
\begin{itemize}
	\item We give a full LDP for terminal error of the Euler--Maruyama method applied to a scalar SDE, which is the first result, to the best of our knowledge, on the LDP of error of numerical methods for nonlinear SDEs. The LDP gives an exponential upper bound for \(\mathbb P\{|U_1^n|\ge\delta\}\) at speed \(n\).
	
	
	\item  We approximate the terminal Euler error by tracking the states visited by the scheme, combining the effects of small Brownian increments, and recording the size and sign of each large increment. We use this approximation to derive the large deviation rate function for the terminal error.
\end{itemize}

\subsection{Proof strategy}
\label{sub:steps}
\noindent\textbf{Step 1. Lamperti transform and an error decomposition.}

We first apply the Lamperti transform to make the diffusion coefficient equal to one. Write
\begin{align}
\Theta(x)&=\int_0^x\frac{\dd u}{\sigma(u)},\qquad
Y_t=\Theta(X_t),\qquad Z_k^n=\Theta(X_k^n),\label{eq:tag-7}\\
\alpha(y)&=\sigma'(x),\qquad
\gamma(y)=\frac{b(x)}{\sigma(x)},\qquad
\beta(y)=\gamma(y)-\tfrac12\alpha(y),\quad x=\Theta^{-1}(y).
\label{eq:tag-8}
\end{align}
By It\^o's formula, the noise coefficient becomes
\(\Theta'(x)\sigma(x)=1\). The transformed equation and its Euler update
are therefore
\begin{align}
\dd Y_t&=\dd W_t+\beta(Y_t)\dd t,\qquad
Y_0=y_0:=\Theta(x_0),\label{eq:tag-9}\\
Y_{k+1}^n&=Y_k^n+\Delta_kW+\tfrac1n\beta(Y_k^n),\qquad
Y_0^n=y_0.\label{eq:aux-lamperti-recursion}
\end{align}
Then denote
\(d_k^n:=Y_k^n-Z_k^n\). We split the terminal error into two terms:
\begin{equation}\label{eq:U1n}
U_1^n=
\bigl[\Theta^{-1}(Y_1)-\Theta^{-1}(Y_n^n)\bigr]
+\bigl[\Theta^{-1}(Z_n^n+d_n^n)-\Theta^{-1}(Z_n^n)\bigr].
\end{equation}
The first term on the right-hand side  of \eqref{eq:U1n}  is negligible at speed \(n\) due to  
the Lipschitz bound on \(\Theta^{-1}\)  and Lemma~\ref{lem:bridge} which states that for any $\eta>0$,
$$\limsup_{n\to\infty}\frac{1}{n}\log\mathbb P\{|Y_1-Y_n^n|>\eta\}=-\infty.$$
To estimate the second term in \eqref{eq:U1n}, we next
construct an equation to approximate \(d_n^n\).

\medskip\noindent\textbf{Step 2. Use empirical controls to approximate \(d_n^n\).}

\noindent\emph{Separate small increments from large ones.}
Fix \(0<\varepsilon\leq1\) and set
\(L_n^\varepsilon=\{0\leq k<n:|\Delta_kW|\geq\varepsilon\}\).
Subtracting the two updates gives the recursion
\begin{equation}\label{eq:terminal-difference-recursion}
d_{k+1}^n=d_k^n+\frac{\beta(Z_k^n+d_k^n)-\gamma(Z_k^n)}{n}
                     +f_k+\eta_k,\qquad d_0^n=0,
\end{equation}
where, for \(x=X_k^n\),
\begin{equation}\label{eq:fk}
f_k=
\begin{cases}
\tfrac12\sigma'(x)(\Delta_kW)^2,& k\notin L_n^\varepsilon,\\[1mm]
\Delta_kW-\bigl[\Theta(x+\sigma(x)\Delta_kW)-\Theta(x)\bigr],
     & k\in L_n^\varepsilon .
\end{cases}
\end{equation}
Here, $\eta_k$ is the remainder defined by \eqref{eq:step-remainder}. 
For a small Brownian increment, the linear Brownian terms cancel and Taylor's
formula gives the first line of \eqref{eq:fk}.
For a large
increment the second line of \eqref{eq:fk} keeps the full nonlinear dependence on its
signed size. The Gaussian tail shows why such increments must be retained:
for each fixed \(a>0\), the law \(\Delta_kW\sim N(0,1/n)\) gives
\[
\mathbb P\{|\Delta_kW|\geq a\}
=\exp\!\left(-\tfrac12na^2+o(n)\right).
\]
Thus even one increment of fixed size has finite cost at speed \(n\).

\medskip\noindent\textit{Record where each increment acts.}
The coefficient \(\tfrac12\sigma'(X_k^n)\) in \eqref{eq:fk} shows why
we need to record the starting state along with the square of a small
increment. Lemma~\ref{lem:compactification} collects the quantities
needed from a transformed state \(y\) into a label \(\iota(y)\) in a
compact space \(K_\sigma\). We record times and labels together on
\(K_\sigma^0=[0,1]\times K_\sigma\).
Write \(\xi_k^n=\iota(Z_k^n)\) for the label of the state at the start
of step \(k\). Let \(\delta_z\) be unit mass at \(z\) and \(\lambda\) Lebesgue measure
on \([0,1]\). The three
terms to be recorded are the time spent at each label, the squared
small increments, and the signed large increments. The corresponding
empirical quantities are
\begin{align}
A_n&=\sum_{k=0}^{n-1}
   \bigl(\lambda\big|_{[t_k,t_{k+1})}\bigr)\otimes\delta_{\xi_k^n},
   \label{eq:empirical-occupation}\\
B_n^\varepsilon&=\sum_{k\notin L_n^\varepsilon}
   (\Delta_kW)^2\,\delta_{(t_k,\xi_k^n)},\label{eq:empirical-variance}\\
\mathfrak m_n^\varepsilon&=
   \bigl((t_k,\xi_k^n,\Delta_kW):k\in L_n^\varepsilon\bigr)
   \quad\text{in increasing order of }k.\label{eq:empirical-marks}
\end{align}
For example, a small increment contributes \((\Delta_kW)^2\) to
\(B_n^\varepsilon\) at precisely the label that determines its coefficient
in \eqref{eq:fk}. Every step, including one with a large increment,
contributes its full interval to \(A_n\). Adding the final label \(\xi_n^n=\iota(Z_n^n)\) gives
the random control
\(\zeta_n^\varepsilon=(A_n,B_n^\varepsilon,\mathfrak m_n^\varepsilon,\xi_n^n)\).
Even if the underlying Euler states escape to infinity, compactness of
\(K_\sigma\) allows their labels to have convergent subsequences.
This compactness is a key ingredient in proving the exponential tightness of \(\zeta_n^\varepsilon\), combined  with the tail bound for \(\sum_{k=0}^{n-1}(\Delta_k W)^2\), in the proof of Theorem \ref{thm:cutoffldp}.

For this cutoff, the space \(\mathcal C_\varepsilon\) consists of all
controls of the form
\[
\zeta=(A,B,\mathfrak m,\xi_f),\qquad
\mathfrak m=((t_j,\xi_j,\delta_j))_{j=1}^{N},\qquad
0\leq t_1\leq\cdots\leq t_N\leq1,\quad |\delta_j|\geq\varepsilon.
\]
The probability measure \(A\) on \(K_\sigma^0\) specifies how time is
spent at the state labels, with \(A([0,t]\times K_\sigma)=t\).
The finite nonnegative measure \(B\) on the same space specifies how
the sum of squared small increments is distributed over times and labels.
Each mark records the time \(t_j\), starting label \(\xi_j\in K_\sigma\),
and signed size \(\delta_j\) of a large increment. The list may be
empty (\(N=0\)), and \(\xi_f\in K_\sigma\) is the terminal label.

\medskip\noindent\textit{Equation with empirical controls.}
To define the drift and increment
terms at every label in \(K_\sigma\), we use the continuous extensions provided by
Lemma~\ref{lem:compactification}. For \(\xi=\iota(y)\), \(x=\Theta^{-1}(y)\),
and \(q\in\R\), these satisfy
\[
\bar\alpha(\xi)=\alpha(y),\qquad
\bar\gamma(\xi)=\gamma(y),\qquad
\bar\beta(\xi,q)=\beta(y+q),\qquad
\bar{\mathfrak j}(\xi,q)=\Theta(x+\sigma(x)q)-y.
\]
The first three entries of \(\zeta\) determine the scalar
error path \(d=d_{A,B,\mathfrak m}\). In the resulting equation, the drift is integrated against \(A\) and
the quadratic coefficient against \(B\). Each mark
contributes a jump at its listed time. With \(d(0-)=0\), the path solves
\begin{equation}\label{eq:controlled-error}
\begin{aligned}
d(t)&=\int_{[0,t]\times K_\sigma}
       [\bar\beta(\xi,d(r))-\bar\gamma(\xi)]\,A(\dd r,\dd\xi)\\
    &\quad+\frac12\int_{[0,t]\times K_\sigma}\bar\alpha(\xi)\,B(\dd r,\dd\xi)
       +\sum_{j:t_j\leq t}[\delta_j-\bar{\mathfrak j}(\xi_j,\delta_j)] .
\end{aligned}
\end{equation}
The controlled equation gives a terminal displacement \(d(1)\) in Lamperti coordinates. Its value in the original coordinate is \(G(\zeta):=\bar{\mathfrak g}(\xi_f,d(1))\). 
For the empirical control \(\zeta_n^\varepsilon\), evaluating \(G\) gives the following approximation to the terminal Euler error \(U_1^n=X_1-X_n^n\).
Since \(\bar{\mathfrak g}(\iota(y),q)=\Theta^{-1}(y+q)-\Theta^{-1}(y)\),
the empirical endpoint in the original coordinate is
\begin{equation}\label{eq:terminal-cutoff-endpoint}
\widehat U_1^{n,\varepsilon}:=G(\zeta_n^\varepsilon)
=\Theta^{-1}\!\bigl(Z_n^n+d_{A_n,B_n^\varepsilon,\mathfrak m_n^\varepsilon}(1)\bigr)
 -\Theta^{-1}(Z_n^n).
\end{equation}
Comparing the empirical form of \eqref{eq:controlled-error} with the
recursion \eqref{eq:terminal-difference-recursion} suggests the approximation of
\(d_n^n\) by \(d_{A_n,B_n^\varepsilon,\mathfrak m_n^\varepsilon}(1)\).
The proof of Lemma~\ref{lem:quotient} bounds the resulting error.

\medskip\noindent\textbf{Step 3. LDPs of $\zeta_n^\varepsilon$ and $\widehat U_1^{n,\varepsilon}$.}\par\nopagebreak[4]

Let \(\mathcal C=\bigcup_{0<\eta\leq1}\mathcal C_\eta\). The cost of a
control \(\zeta=(A,B,\mathfrak m,\xi_f)\in\mathcal C\) is
\begin{equation}\label{eq:IAB}
I\colon\mathcal C\to[0,\infty],\qquad
I(\zeta)\mathrel{=}I(A,B,\mathfrak m,\xi_f)
\mathrel{\coloneqq}I_{\rm sv}(A,B)+\frac12\sum_{j=1}^{N}\delta_j^2,
\end{equation}
where the first term is given by~\eqref{eq:Isv-definition}.
Writing \(B=vA+B^\perp\) separates the part with density \(v\) relative
to \(A\) from the mass \(B^\perp\) carried by a set of \(A\)-measure
zero. Their contributions give
\[
I_{\rm sv}(A,B)=\int_{K_\sigma^0}\ell(v)\,\dd A
                 +\tfrac12 B^\perp(K_\sigma^0),\qquad
\ell(u)=
\begin{cases}
\tfrac12(u-1-\log u),&u>0,\\
+\infty,&u\leq0.
\end{cases}
\]
For example, if the list is empty and \(B=uA\) for a constant \(u>0\),
the entire cost is \(\ell(u)\). For each mark \(j\), observing a Brownian increment near the prescribed value \(\delta_j\) contributes \(\delta_j^2/2\) to the rate function.
At a fixed cutoff \(\varepsilon\), the rate function for
\(\zeta_n^\varepsilon\) is the restriction
\(I_\varepsilon:=\left.I\right|_{\mathcal C_\varepsilon}\).

Theorem~\ref{thm:cutoffldp} proves
the LDP for \(\zeta_n^\varepsilon\), using changes of measure
as in \cite{BoueDupuis,BudhirajaDupuisMaroulas,DupuisEllis}. For the lower bound, fix a neighbourhood \(O\) of a finite-cost control
\(\zeta\). We first approximate \((A,B)\) by regular pairs as in
Lemma~\ref{lem:regularize}. Lemma~\ref{lem:marked-approximation} then
gives a control in \(O\) with piecewise constant
state and variance profiles. The argument of changes of measure and a relative entropy estimate in Lemma \ref{lem:recovery-tilted-law} then give
\(\liminf_{n\to\infty}n^{-1}\log\mathbb P\{\zeta_n^\varepsilon\in O\}
\geq-I_\varepsilon(\zeta)\) under the original law.

For the local upper bound, we fix the indices of the \(N\) large
increments. Conditional exponential moments of the squared small
increments in Lemma \ref{lem:gaussianexp} yield the cost \(I_{\rm sv}(A,B)\) through its dual formula given by Lemma \ref{lem:Isv-dual}. 
The Gaussian density gives a cost of \(\delta_j^2/2\) for each large
increment. Exponential tightness allows us to choose a compact set
outside which the probability decays at any prescribed exponential rate.
Its intersection with a closed set \(F\) is covered by finitely many
of the neighbourhoods used for the local upper bound. This gives the
upper bound for \(F\).

The continuity
of \(G\) and the contraction principle then give the LDP of  $\widehat U_1^{n,\varepsilon}$ for any fixed  $\varepsilon\in(0,1]$. 
For \(e\in\R\), taking the infimum of the costs of all controls in
\(\mathcal C_\varepsilon\) with terminal error \(G(\zeta)=e\) gives
\begin{equation}
J_\varepsilon(e)
=\inf\left\{
I_\varepsilon(A,B,\mathfrak m,\xi_f):
(A,B,\mathfrak m,\xi_f)\in\mathcal C_\varepsilon,\
\bar{\mathfrak g}\bigl(\xi_f,d_{A,B,\mathfrak m}(1)\bigr)=e
\right\},
\label{eq:tag-5}
\end{equation}
with the convention \(\inf\varnothing=+\infty\).

\medskip\noindent\textbf{Step 4. Remove the cutoff.}

Lemma~\ref{lem:quotient} bounds the error between $\Theta^{-1}(Z_n^n+d_n^n)-\Theta^{-1}(Z_n^n)$ and $\widehat{U}_1^{n,\varepsilon}$ on the events $\{\mu_n([0,1])\le M\}$, $M>0$. 
Here \(\mu_n([0,1])=\sum_{k=0}^{n-1}(\Delta_kW)^2\). For \(M>1\),
\eqref{eq:tag-13} bounds \(\mathbb P\{\mu_n([0,1])\geq M\}\) by
\(e^{-n\ell(M)}\), with \(\ell(M)\to\infty\) as \(M\to\infty\).
Together with the estimate for the first term of \eqref{eq:U1n} in
Lemma~\ref{lem:bridge}, these two bounds give the exponentially good
approximation of \(U_1^n\) by
\(\widehat U_1^{n,\varepsilon}\) in Proposition~\ref{prop:expapprox}.

Proposition~\ref{prop:approximation-principle} now transfers the
fixed-cutoff LDPs to \(U_1^n\). Proposition~\ref{prop:rate} identifies
the resulting rate at each \(e\in\R\) as
\(\lim_{\varepsilon\downarrow0}J_\varepsilon(e)\). As the cutoff decreases,
the spaces \(\mathcal C_\varepsilon\) increase to \(\mathcal C\), so the
infimum in \eqref{eq:tag-5} becomes
\begin{equation}
J_{\rm un}(e)
=\inf\left\{
I(A,B,\mathfrak m,\xi_f):
(A,B,\mathfrak m,\xi_f)\in\mathcal C,\
\bar{\mathfrak g}\bigl(\xi_f,d_{A,B,\mathfrak m}(1)\bigr)=e
\right\}.
\label{eq:tag-6}
\end{equation}

\emph{Organization.} The paper is organized as follows. Section~\ref{sec:marked-controls}
derives the one-step error estimates in Lamperti coordinates, constructs
the compact space of state labels, and proves the comparison between
the controlled endpoint and the discrete error.
Section~\ref{sec:terminal-ldp} proves Theorem~\ref{thm:main} using
Theorem~\ref{thm:cutoffldp}, which is proved in
Section~\ref{sec:fixed-cutoff-ldp}. Appendix~\ref{app:wellposedness}
proves well-posedness and continuity of the controlled equation.
\section{Discrete error decomposition and compactified controls}
\label{sec:marked-controls}

We begin with some notations. For a Polish space \(E\), write \(\mathcal P(E)\) for the Borel probability
measures and \(\Mplus(E)\) for the finite nonnegative Borel measures on
\(E\). Both spaces carry the weak topology. We write \(\mu_n\Rightarrow\mu\)
for convergence against every bounded continuous real-valued function
\(f\in C_b(E)\), as in \eqref{eq:weaktoM}.
We denote Lebesgue measure on \([0,1]\) by \(\lambda\) and the unit point
mass at \(z\) by \(\delta_z\). For measurable sets \(B,D\), the restriction
\(\mu\big|_D\) satisfies \((\mu\big|_D)(B)=\mu(B\cap D)\).
The pushforward of \(\mu\) by a measurable map \(T\) is defined by
\((T_\#\mu)(B)=\mu(T^{-1}(B))\). In particular, \(\mathcal L_Q(X)=X_\#Q\)
is the law of a random element \(X\) under \(Q\).
For probability measures \(Q,P\) on the same measurable space, the relative
entropy \(H(Q\mid P)\) is \(\int\log(\mathrm dQ/\mathrm dP)\dd Q\) if
\(Q\ll P\) and \(+\infty\) otherwise. The normal law \(N(m,v)\) has mean
\(m\) and variance \(v\). A path \(x\) on \([0,1]\) is c\`adl\`ag if it is
right-continuous on \([0,1)\) and the left limit
\(x(t-)=\lim_{s\uparrow t}x(s)\) exists for every \(t\in(0,1]\).

We first justify the two Euler updates used in the introduction
and derive the remainder in \eqref{eq:terminal-difference-recursion}.
Differentiating \eqref{eq:tag-7} gives
\(\Theta'=1/\sigma\) and \(\Theta''=-\sigma'/\sigma^2\). Hence the drift
in It\^o's formula for \(Y_t\) is
\[
\Theta'(X_t)b(X_t)+\tfrac12\Theta''(X_t)\sigma(X_t)^2
=\frac{b(X_t)}{\sigma(X_t)}-\tfrac12\sigma'(X_t)=\beta(Y_t).
\]
This gives \eqref{eq:tag-9}, whose Euler update is
\eqref{eq:aux-lamperti-recursion}. The bounds
\(1/\sigma_{\mathrm{max}}\leq\Theta'\leq1/\sigma_{\mathrm{min}}\) show that
\(\Theta\) maps \(\R\) increasingly onto \(\R\). Its inverse satisfies
\((\Theta^{-1})'=\sigma\circ\Theta^{-1}\) and is globally Lipschitz
with constant \(\sigma_{\mathrm{max}}\). The chain rule and
Assumption~\ref{ass:coeff} likewise show that \(\alpha,\gamma,\beta\)
in \eqref{eq:tag-8} are bounded and globally Lipschitz.

For the original Euler scheme in the same coordinate, write
\(Z_{k+1}^n=\mathcal S_{n,\Delta_kW}(Z_k^n)\), where
\begin{equation}
\mathcal S_{n,\delta}(y)
=
\Theta\!\left(
\Theta^{-1}(y)
+\sigma(\Theta^{-1}(y))\delta
+\frac1n b(\Theta^{-1}(y))
\right),
\quad
y,\delta\in\R,
\label{eq:transformed-euler-update}
\end{equation}
and define the one-step remainder
\begin{equation}
r_n(y,\delta)
:=
\mathcal S_{n,\delta}(y)-y-\delta-\frac1n\gamma(y)
+\frac12\alpha(y)\delta^2.
\label{eq:lamperti-one-step-remainder}
\end{equation}
Subtracting
\(Z_{k+1}^n=\mathcal S_{n,\Delta_kW}(Z_k^n)\) from
\eqref{eq:aux-lamperti-recursion}, using \(Y_k^n=Z_k^n+d_k^n\), and
expressing \(\mathcal S_{n,\Delta_kW}(Z_k^n)-Z_k^n\) through
\eqref{eq:lamperti-one-step-remainder}, we obtain, for every
\(0\leq k<n\),
\begin{equation}
d_{k+1}^n
=
d_k^n
+\frac1n
\bigl[\beta(Z_k^n+d_k^n)-\gamma(Z_k^n)\bigr]
+\frac12\alpha(Z_k^n)(\Delta_kW)^2
-r_n(Z_k^n,\Delta_kW).
\label{eq:discrete-error-small}
\end{equation}
The first line of \eqref{eq:fk} uses this expansion when
\(k\notin L_n^\varepsilon\). On those steps,
\(|\Delta_kW|^3\leq\varepsilon(\Delta_kW)^2\), so the cubic terms
can be summed using the total squared increments. Taylor's formula
gives the uniform estimate
\begin{equation}
|r_n(y,\delta)|
\leq
C\left(
|\delta|^3+\frac{|\delta|}{n}+\frac1{n^2}
\right),
\qquad y,\delta\in\R,\quad n\geq1,
\label{eq:lamperti-remainder-global}
\end{equation}
where \(C\) depends only on the bounds in Assumption~\ref{ass:coeff}.


For \(k\in L_n^\varepsilon\), the cubic term in the bound
\eqref{eq:lamperti-remainder-global} need not tend to zero as the mesh
is refined. Returning to the original updates and subtracting
\(Z_{k+1}^n=\mathcal S_{n,\Delta_kW}(Z_k^n)\) from
\eqref{eq:aux-lamperti-recursion} gives
\begin{equation}
d_{k+1}^n
=
d_k^n+\frac1n\beta(Z_k^n+d_k^n)+\Delta_kW
-\bigl[\mathcal S_{n,\Delta_kW}(Z_k^n)-Z_k^n\bigr].
\label{eq:discrete-error-marked}
\end{equation}

To separate the nonlinear increment effect from the drift contribution,
define
\begin{equation}\label{eq:Jdeltay}
\mathcal S_\delta(y)
:=
\Theta\bigl(
\Theta^{-1}(y)+\sigma(\Theta^{-1}(y))\delta
\bigr),
\qquad
y,\delta\in\R.
\end{equation}
Unlike \(\mathcal S_{n,\delta}\), the map \(\mathcal S_\delta\) omits the drift increment
\(\frac1n b(\Theta^{-1}(y))\). Since \(\Theta'\) and \(b\) are bounded,
the mean value theorem gives
\begin{equation}
\bigl|\mathcal S_{n,\delta}(y)-\mathcal S_\delta(y)\bigr|
\leq
\frac{\|\Theta'\|_\infty\|b\|_\infty}{n},
\qquad
y,\delta\in\R.
\label{eq:lamperti-drift-perturbation}
\end{equation}

The second line of \eqref{eq:fk} can now be written as
\(f_k=\Delta_kW-[\mathcal S_{\Delta_kW}(Z_k^n)-Z_k^n]\).
Using \eqref{eq:discrete-error-small} on small steps and
\eqref{eq:discrete-error-marked} on large ones gives the following remainder in
\eqref{eq:terminal-difference-recursion},
\begin{equation}\label{eq:step-remainder}
\eta_k=
\begin{cases}
-r_n(Z_k^n,\Delta_kW),&k\notin L_n^\varepsilon,\\[1mm]
\dfrac1n\gamma(Z_k^n)
-\bigl[\mathcal S_{n,\Delta_kW}(Z_k^n)-\mathcal S_{\Delta_kW}(Z_k^n)\bigr],
&k\in L_n^\varepsilon.
\end{cases}
\end{equation}
Indeed, on a large step the sum of \(f_k\) and \(\eta_k\) is
\(\Delta_kW+\gamma(Z_k^n)/n-[\mathcal S_{n,\Delta_kW}(Z_k^n)-Z_k^n]\).
The term \(\gamma(Z_k^n)/n\) cancels its counterpart in the drift of
\eqref{eq:terminal-difference-recursion}. On small steps,
\eqref{eq:lamperti-remainder-global} controls \(\eta_k\), while on large
steps \eqref{eq:lamperti-drift-perturbation} gives \(|\eta_k|\leq C/n\).

After omitting \(\eta_k\), the recursion
\eqref{eq:terminal-difference-recursion} and the coordinate conversion in
the second bracket of \eqref{eq:U1n} depend on \(y=Z_k^n\) through
the values
\[
\alpha(y),\quad
\gamma(y),\quad
\beta(y+q),\quad
\Theta^{-1}(y+q)-\Theta^{-1}(y),\quad
\mathcal S_q(y)-y,
\qquad q\in\R.
\]
Since the state space \(\R\) is not compact, we next encode these values
by a state label
\(\xi=\iota(y)\) in a compact space \(K_\sigma\).
\begin{lemma}[Coefficient compactification]
\label{lem:compactification}
There exist a compact metric space \(K_\sigma\) and a continuous map
\(\iota:\R\to K_\sigma\) with dense image such that each of the quantities
\begin{equation}
\alpha(y),\quad
\gamma(y),\quad
\beta(y+q),\quad
\Theta^{-1}(y+q)-\Theta^{-1}(y),\quad
{\mathcal S}_q(y)-y,
\qquad q\in\mathbb Q,
\label{eq:coefficient-observables}
\end{equation}
is determined by the state label \(\iota(y)\), with \(q\) held fixed in
the last three families. Each induced real-valued function on
\(\iota(\R)\) has a unique continuous extension to \(K_\sigma\).

Denote the extensions of the first two functions by
\(\bar\alpha,\bar\gamma:K_\sigma\to\R\). For every \(y\in\R\), they satisfy
\[
\begin{aligned}
\bar\alpha(\iota(y))
&=\alpha(y)=\sigma'\!\left(\Theta^{-1}(y)\right),\\
\bar\gamma(\iota(y))
&=\gamma(y)=\frac{b\!\left(\Theta^{-1}(y)\right)}
                  {\sigma\!\left(\Theta^{-1}(y)\right)}.
\end{aligned}
\]

For each \(q\in\mathbb Q\), denote the extensions of the last three
functions in \eqref{eq:coefficient-observables} by
\(\bar\beta(\cdot,q)\), \(\bar{\mathfrak g}(\cdot,q)\), and
\(\bar{\mathfrak j}(\cdot,q)\), respectively. The resulting three maps on
\(K_\sigma\times\mathbb Q\) extend uniquely in the parameter \(q\) to
jointly continuous maps from \(K_\sigma\times\R\) to \(\R\).
Denoting these maps by the same symbols, we have, for every \(y,q\in\R\),
\[
\begin{aligned}
\bar\beta(\iota(y),q)
&=\beta(y+q),\\
\bar{\mathfrak g}(\iota(y),q)
&=\Theta^{-1}(y+q)-\Theta^{-1}(y),\\
\bar{\mathfrak j}(\iota(y),q)
&={\mathcal S}_q(y)-y.
\end{aligned}
\]
\end{lemma}

\begin{proof}[Proof of Lemma \ref{lem:compactification}]
We first record Lipschitz bounds in the parameter. Let
\(L_\beta\) be the global Lipschitz constant of \(\beta\). Since
\((\Theta^{-1})'=\sigma\circ\Theta^{-1}\leq\sigma_{\mathrm{max}}\), and
since, with \(x=\Theta^{-1}(y)\),
\[
\frac{\partial}{\partial q}\mathcal S_q(y)
=\frac{\sigma(x)}{\sigma(x+\sigma(x)q)}
\in\Bigl(0,\frac{\sigma_{\mathrm{max}}}{\sigma_{\mathrm{min}}}\Bigr],
\]
we have, for \(q,q'\in\R\) and uniformly in \(y\in\R\),
\begin{align}\label{eq:beta1}
|\beta(y+q)-\beta(y+q')|
&\leq L_\beta|q-q'|,\\\label{eq:beta2}
|\Theta^{-1}(y+q)-\Theta^{-1}(y+q')|
&\leq \sigma_{\mathrm{max}}\,|q-q'|,\\\label{eq:beta3}
|{\mathcal S}_q(y)
-{\mathcal S}_{q'}(y)|
&\leq
\frac{\sigma_{\mathrm{max}}}{\sigma_{\mathrm{min}}}|q-q'|.
\end{align}

Let \(\mathscr F\) denote the countable family of functions in
\eqref{eq:coefficient-observables}, with \(q\in\mathbb Q\).
Each \(f\in\mathscr F\) is continuous. Each \(f\) is also
bounded, since \(\alpha,\gamma,\beta\) are bounded, and \eqref{eq:beta2}
and \eqref{eq:beta3} with \(q'=0\) and \(\mathcal S_0(y)=y\) give
\(|\Theta^{-1}(y+q)-\Theta^{-1}(y)|\leq\sigma_{\mathrm{max}}|q|\) and
\(|\mathcal S_q(y)-y|\leq(\sigma_{\mathrm{max}}/\sigma_{\mathrm{min}})|q|\).

Enumerate \(\mathscr F=\{\chi_i:i\in\mathbb N\}\) and define
\[
\iota(y)
:=
\bigl(\chi_i(y)\bigr)_{i\geq1}
\in
\mathcal P:=\prod_{i\geq1}
[-\|\chi_i\|_\infty,\|\chi_i\|_\infty],
\qquad
 K_\sigma:=\overline{\iota(\R)},
\]
where the closure is taken in the product topology.
By Tychonoff's theorem, \(\mathcal P\) is compact. Since
\(\mathscr F\) is countable, the product topology on \(\mathcal P\) is
induced by the metric
\[
\mathrm{d}_{K_\sigma}(\xi,\eta)
:=
\sum_{i=1}^{\infty}
2^{-i}
\bigl(|\xi_{i}-\eta_{i}|\wedge 1\bigr).
\]
Hence the closed subset \(K_\sigma\) of \(\mathcal P\) is a
compact metric space. The map \(\iota\) is continuous because each of its
coordinates is continuous, and \(\iota(\R)\) is dense in \(K_\sigma\) by
definition. For each \(i\), the coordinate projection
\(\xi\mapsto\xi_i\) is continuous on \(K_\sigma\) and equals \(\chi_i\) on
\(\iota(\R)\). This gives continuous extensions of all \(f\in\mathscr F\),
which are unique by density. This defines \(\bar\alpha,\bar\gamma\) and, for
\(q\in\mathbb Q\),
\(\bar\beta(\cdot,q)\), \(\bar{\mathfrak g}(\cdot,q)\), and
\(\bar{\mathfrak j}(\cdot,q)\).

It remains to extend the parameter \(q\) from \(\mathbb Q\) to \(\R\).
Let \(q,q'\in\mathbb Q\). The continuous function
\(\xi\mapsto|\bar\beta(\xi,q)-\bar\beta(\xi,q')|\) is bounded by
\(L_\beta|q-q'|\) on the dense set \(\iota(\R)\) by \eqref{eq:beta1}, and
hence on \(K_\sigma\). By \eqref{eq:beta2} and \eqref{eq:beta3}, the same
holds for \(\bar{\mathfrak g}\) and \(\bar{\mathfrak j}\), with constants
\(\sigma_{\mathrm{max}}\) and
\(\sigma_{\mathrm{max}}/\sigma_{\mathrm{min}}\). Consequently, for every
\(q\in\R\) and every rational sequence \(q_m\to q\),
\[
\bar\beta(\cdot,q_m),\qquad
\bar{\mathfrak g}(\cdot,q_m),\qquad
\bar{\mathfrak j}(\cdot,q_m)
\]
are uniformly Cauchy on \(K_\sigma\). Their limits are independent of the
approximating sequence and define the three families for all \(q\in\R\).
As uniform limits, these maps are continuous on
\(K_\sigma\) for each fixed \(q\), and the three Lipschitz bounds in the
parameter extend to all \(q,q'\in\R\). Moreover, since \(\beta\) is
continuous, \(\bar\beta(\iota(y),q)=\lim_m\beta(y+q_m)=\beta(y+q)\) for all
\(y,q\in\R\), and the same argument gives the identities for
\(\bar{\mathfrak g}\) and \(\bar{\mathfrak j}\).

Finally, if
$
(\xi_m,q_m)\to(\xi,q)$
in $K_\sigma\times\R$,
then
\[
\begin{aligned}
|\bar\beta(\xi_m,q_m)-\bar\beta(\xi,q)|
&\leq
|\bar\beta(\xi_m,q_m)-\bar\beta(\xi_m,q)|
+
|\bar\beta(\xi_m,q)-\bar\beta(\xi,q)|\\
&\leq
L_\beta|q_m-q|
+
|\bar\beta(\xi_m,q)-\bar\beta(\xi,q)|
{\rightarrow}0.
\end{aligned}
\]
The same argument applies to
\(\bar{\mathfrak g}\) and \(\bar{\mathfrak j}\), with parameter-Lipschitz
constants \(\sigma_{\mathrm{max}}\) and
\(\sigma_{\mathrm{max}}/\sigma_{\mathrm{min}}\), respectively. This proves the asserted joint
continuity.
\end{proof}
For a Polish space \(E\) and \(\mu_n,\mu\in\Mplus(E)\),
\begin{equation}\label{eq:weaktoM}
\mu_n\Rightarrow\mu
\quad\Longleftrightarrow\quad
\int_E \varphi\dd\mu_n
{\rightarrow}
\int_E \varphi\dd\mu
\qquad
\text{for every }\varphi\in C_b(E).
\end{equation}
Set
$
K_\sigma^0:=[0,1]\times K_\sigma.
$
The size of the increments is measured by the scalar
quadratic-variation measure
\begin{equation}
\mu_n:=\sum_{k=0}^{n-1}(\Delta_kW)^2\delta_{t_{k+1}}.
\label{eq:tag-10}
\end{equation}
The measure \(\mu_n\) assigns mass \((\Delta_kW)^2\) to \(t_{k+1}\), following the convention in \cite{DGS2017}, which we use in Proposition~\ref{prop:qv}.

We now verify how the empirical measures and marks in
\eqref{eq:empirical-occupation}--\eqref{eq:empirical-marks} encode
\eqref{eq:terminal-difference-recursion}. Recall that
\(\xi_k^n=\iota(Z_k^n)\). On a single grid interval the drift equals
\[
\frac{\beta(Z_k^n+d_k^n)-\gamma(Z_k^n)}{n}
=\int_{[t_k,t_{k+1})\times K_\sigma}
   [\bar\beta(\xi,d_k^n)-\bar\gamma(\xi)]\,A_n(\dd r,\dd\xi).
\]
For a small increment, \(B_n^\varepsilon\) has an atom of mass
\((\Delta_kW)^2\) at \((t_k,\xi_k^n)\). Its contribution to the
second integral in \eqref{eq:controlled-error} is therefore
\(f_k=\tfrac12\bar\alpha(\xi_k^n)(\Delta_kW)^2\). For a large increment,
the corresponding mark contributes
\(f_k=\Delta_kW-\bar{\mathfrak j}(\xi_k^n,\Delta_kW)\) to the sum in that
equation.

We also need the measure that includes the squares of all increments,
\[
B_n:=\sum_{k=0}^{n-1}(\Delta_kW)^2\delta_{(t_k,\xi_k^n)}
    =B_n^\varepsilon+\sum_{k\in L_n^\varepsilon}
                    (\Delta_kW)^2\delta_{(t_k,\xi_k^n)} .
\]
The scalar measure \(\mu_n\) puts each 
$(\Delta_k W)^2$ at \(t_{k+1}\).
\(B_n\) puts $(\Delta_k W)^2$ at \(t_k\) with the label that \(A_n\) uses throughout
\([t_k,t_{k+1})\). In particular \(B_n(K_\sigma^0)=\mu_n([0,1])\).

Any weak limit \((A,B)\) of \((A_n,B_n^\varepsilon)\) consists of
a probability measure \(A\) with Lebesgue time marginal and a finite
nonnegative measure \(B\) on \(K_\sigma^0\). We therefore define
\begin{equation}
\mathcal Q_\sigma
:=
\left\{
(A,B)\in\mathcal P(K_\sigma^0)\times\Mplus(K_\sigma^0):
(\pi_1)_\#A=\lambda
\right\},
\label{eq:measure-pair-space}
\end{equation}
where \(\pi_1:K_\sigma^0\to[0,1]\) denotes the projection
\(\pi_1(t,\xi)=t\). In particular,
\((A_n,B_n^\varepsilon)\in\mathcal Q_\sigma\).
Since \(K_\sigma^0\) is compact and metrizable,
\(\mathcal P(K_\sigma^0)\) and
\(\Mplus(K_\sigma^0)\), equipped with their weak topologies, are
Polish spaces. Since \(\pi_1\) is continuous, the map
\(A\mapsto(\pi_1)_\#A\) is continuous for the weak topology. Hence
\(\mathcal Q_\sigma\) is closed in
\(\mathcal P(K_\sigma^0)\times
\Mplus(K_\sigma^0)\), and is therefore Polish.

To fix a metric for this topology, choose a sequence
\((\varphi_j)_{j\geq1}\) dense in the closed unit ball of
\(C(K_\sigma^0)\) for the uniform norm, with \(\varphi_1\equiv1\).
For \(\mu,\nu\in\Mplus(K_\sigma^0)\), set
\begin{equation}\label{eq:dw}
d_{\mathrm w}(\mu,\nu)
:=
\sum_{j=1}^{\infty}2^{-j}
\left(
\left|\int_{K_\sigma^0}\varphi_j\,\dd\mu
-\int_{K_\sigma^0}\varphi_j\,\dd\nu\right|\wedge1
\right).
\end{equation}
Convergence in \(d_{\mathrm w}\) implies convergence of the
total masses, by the choice \(\varphi_1\equiv1\), and hence convergence of
\(\int\varphi\,\dd\mu\) for every \(\varphi\in C(K_\sigma^0)\), by
density. Thus \(d_{\mathrm w}\) induces the weak topology on
\(\Mplus(K_\sigma^0)\). We equip \(\mathcal Q_\sigma\) with the metric
\[
d_{\mathcal Q_\sigma}\bigl((A,B),(A',B')\bigr)
:=
d_{\mathrm w}(A,A')+d_{\mathrm w}(B,B'),
\]
which induces the restriction of the product weak topology.
In particular,
\[
d_{\mathcal Q_\sigma}\bigl((A_m,B_m),(A,B)\bigr){\rightarrow}0
\quad\Longleftrightarrow\quad
A_m\Rightarrow A\ \text{and}\ B_m\Rightarrow B.
\]

Fix \(0<\varepsilon\leq1\) and let \(n_\ell\to\infty\),
allowing the sample point \(\omega_\ell\) to vary with \(\ell\).
The realized lists \(\mathfrak m_{n_\ell}^{\varepsilon}(\omega_\ell)\)
converge to a finite ordered list
\begin{equation}
\mathfrak{m}
=
((t_j,\xi_j,\delta_j))_{j=1}^{N},
\qquad
0\leq t_1\leq\cdots\leq t_N\leq1,
\qquad
|\delta_j|\geq\varepsilon,
\label{eq:ordered-mark-list}
\end{equation}
precisely when they eventually have length \(N\) and their
\(j\)th marks converge to \((t_j,\xi_j,\delta_j)\) for each
\(j=1,\ldots,N\). Here \(\xi_j\in K_\sigma\) and
\(N\in\mathbb N_0:=\{0,1,2,\ldots\}\). The case \(N=0\) corresponds to the
empty list.
For \((A,B)\in\mathcal Q_\sigma\) and a list satisfying
\eqref{eq:ordered-mark-list}, the drift and the two contributions to the
error give \eqref{eq:controlled-error} with initial value \(d(0-)=0\)
before any jump at time zero. This use of measures is related to the
measure-control viewpoint of \cite{DalMasoRampazzo} and to relaxed
controls \cite{Warga}. Appendix~\ref{app:wellposedness} proves that the
equation has a unique c\`adl\`ag solution \(d_{A,B,\mathfrak m}\).
For the empirical controls
\eqref{eq:empirical-occupation}--\eqref{eq:empirical-marks},
\eqref{eq:controlled-error} becomes
\begin{equation}
\begin{split}
d_{A_n,B_n^\varepsilon,\mathfrak m_n^\varepsilon}(t)
&=
\sum_{k=0}^{n-1}
\int_{[0,t]\cap[t_k,t_{k+1})}
\bigl[
\bar\beta\bigl(
\iota(Z_k^n),
d_{A_n,B_n^\varepsilon,\mathfrak m_n^\varepsilon}(r)
\bigr)
-\bar\gamma\bigl(\iota(Z_k^n)\bigr)
\bigr]\dd r\\
&\quad+
\frac12
\sum_{\substack{k\notin L_n^\varepsilon\\ t_k\leq t}}
(\Delta_kW)^2
\bar\alpha\bigl(\iota(Z_k^n)\bigr)+
\sum_{\substack{k\in L_n^\varepsilon\\ t_k\leq t}}
\bigl[
\Delta_kW
-\bar{\mathfrak j}\bigl(
\iota(Z_k^n),\Delta_kW
\bigr)
\bigr],
\end{split}
\label{eq:empirical-controlled-error}
\end{equation}
with \(d_{A_n,B_n^\varepsilon,\mathfrak m_n^\varepsilon}(0-)=0\).
The integral is continuous in \(t\). At each grid time \(t_k\), the
sums add \(f_k\). Between grid times the derivative is
\(\beta(Z_k^n+d)-\gamma(Z_k^n)\). This verifies the evolution used in
the introduction to define \(\widehat U_1^{n,\varepsilon}\) in
\eqref{eq:terminal-cutoff-endpoint}. The next lemma bounds the combined
effect of omitting \(\eta_k\) and integrating the drift exactly on each
grid interval.

\begin{lemma}
\label{lem:quotient}
There exists a constant \(C<\infty\), depending only on the
bounds in Assumption~\ref{ass:coeff}, such that, for every \(M>0\),
\(n\geq1\) and \(\varepsilon\in(0,1]\), on
\(\{\mu_n([0,1])\leq M\}\),
\begin{equation}
\left|
\Theta^{-1}(Y_n^n)-\Theta^{-1}(Z_n^n)
-\widehat U_1^{n,\varepsilon}
\right|
\leq
C\left(
\varepsilon M+\sqrt{\frac Mn}
+\frac{M+1}n+\min\Bigl\{1,\frac{M}{n\varepsilon^2}\Bigr\}
\right).
\label{eq:terminal-quotient-comparison}
\end{equation}
\end{lemma}

\begin{proof}[Proof of Lemma \ref{lem:quotient}]
Throughout the proof, \(C\) denotes a constant depending
only on the bounds in Assumption~\ref{ass:coeff}, and we work on the event
\(\{\mu_n([0,1])\leq M\}\).

We first bound \(f_k\) in \eqref{eq:fk} and \(\eta_k\) in
\eqref{eq:step-remainder}. Since \(\mathcal S_0(y)=y\),
\(\partial_\delta \mathcal S_\delta(y)|_{\delta=0}=1\), and
\(\partial_\delta^2\mathcal S_\delta(y)=\Theta''(x+\sigma(x)\delta)\sigma(x)^2\) with
\(x=\Theta^{-1}(y)\) is bounded, Taylor's formula gives
\begin{equation}
\left|
{\mathcal S}_\delta(y)-y-\delta
\right|
\leq C\delta^2,
\qquad
y,\delta\in\R.
\label{eq:lamperti-global-quadratic}
\end{equation}
Together with the boundedness of \(\alpha\), this gives
\begin{equation}\label{eq:fkbound}
\sum_{k=0}^{n-1}|f_k|
\leq C\sum_{k=0}^{n-1}(\Delta_kW)^2
=C\mu_n([0,1])\leq CM.
\end{equation}
For the remainders, each \(k\in L_n^\varepsilon\)
contributes at least \(\varepsilon^2\) to \(\mu_n([0,1])\), so
\(\#L_n^\varepsilon\leq\min\{n,M/\varepsilon^2\}\). For
\(k\in L_n^\varepsilon\), \eqref{eq:lamperti-drift-perturbation} and the
boundedness of \(\gamma\) give \(|\eta_k|\leq C/n\). For
\(k\notin L_n^\varepsilon\), we have \(|\Delta_kW|<\varepsilon\), so
\(|\Delta_kW|^3\leq\varepsilon(\Delta_kW)^2\). Combining
\eqref{eq:lamperti-remainder-global} with the Cauchy--Schwarz bound
\(\sum_k|\Delta_kW|\leq\sqrt{nM}\), we obtain
\begin{equation}
\sum_{k=0}^{n-1}|\eta_k|
\leq
C\left(
\varepsilon M
+\sqrt{\frac Mn}
+\frac{M+1}{n}
+\min\left\{1,\frac{M}{n\varepsilon^2}\right\}
\right).
\label{eq:terminal-remainder-sum}
\end{equation}

Next, we describe the solution
\(\mathsf d:=d_{A_n,B_n^\varepsilon,\mathfrak m_n^\varepsilon}\) of
\eqref{eq:empirical-controlled-error}, with \(\mathsf d(0-)=0\). Let
\(\Phi_k^s\) denote the flow generated by the vector field
\[
F_k(x):=\beta(Z_k^n+x)-\gamma(Z_k^n),\qquad
x\in\R.
\]
Since \(\beta\) and \(\gamma\) are bounded and \(\beta\) is
globally \(L_\beta\)-Lipschitz, the vector fields \(F_k\) are bounded by
\(C\) and \(L_\beta\)-Lipschitz, uniformly in \(k\) and \(n\).
By \eqref{eq:empirical-controlled-error}, \(\mathsf d\) jumps
by \(f_k\) at \(t_k\) and solves
\(\dot{\mathsf d}=F_k(\mathsf d)\) on \((t_k,t_{k+1})\). Hence
\begin{equation}
\mathsf d(t_{k+1}-)
=
\Phi_k^{1/n}\!\left(\mathsf d(t_k-)+f_k\right),
\qquad
0\leq k<n.
\label{eq:terminal-quotient-flow}
\end{equation}
Since \eqref{eq:empirical-controlled-error} has no jump at
\(t_n=1\),
\(
\mathsf d(1)=\mathsf d(t_n-).
\)

To compare the discrete recursion \eqref{eq:terminal-difference-recursion}
with the flow recursion \eqref{eq:terminal-quotient-flow},
introduce
\begin{equation}\label{eq:wrec}
w_0=0,
\qquad
w_{k+1}
=
\Phi_k^{1/n}(w_k)+f_k,
\qquad
0\leq k<n.
\end{equation}
We first compare \(w\) with \(d^n\). Since
\(|\Phi_k^r(x)-x|\leq Cr\), the Lipschitz continuity of \(F_k\) gives
\[
\left|
\Phi_k^{1/n}(x)
-x
-\frac1nF_k(x)
\right|
\leq
L_\beta\int_0^{1/n}\bigl|\Phi_k^r(x)-x\bigr|\,\dd r
\leq
\frac{C}{n^2},
\qquad x\in\R.
\]
Subtracting \eqref{eq:wrec} from
\eqref{eq:terminal-difference-recursion} and using this estimate, we
obtain
\[
|d_{k+1}^n-w_{k+1}|
\leq
\left(1+\frac{L_\beta}{n}\right)|d_k^n-w_k|
+|\eta_k|
+\frac{C}{n^2}.
\]
Since \(d_0^n=w_0=0\), the discrete Gronwall inequality yields
\begin{equation}
|d_n^n-w_n|
\leq
C\left(
\sum_{k=0}^{n-1}|\eta_k|+\frac1n
\right).
\label{eq:terminal-euler-flow-error}
\end{equation}

We next compare \(w\) with \(\mathsf d\) at the grid times.
Gronwall's inequality gives
\(
\left|\Phi_k^s(x+z)-\Phi_k^s(x)\right|
\le e^{L_\beta s}|z|
\)
for \(0\le s\le1\), uniformly in \(k\) and \(n\). Hence, using the
integral form of the flow and the \(L_\beta\)-Lipschitz continuity of
\(F_k\),
\[
\left|
\Phi_k^{s}(x+z)-\Phi_k^{s}(x)-z
\right|
\le
L_\beta\int_0^s
\left|
\Phi_k^r(x+z)-\Phi_k^r(x)
\right|\,\dd r\le
\left(e^{L_\beta s}-1\right)|z|
\leq\frac Cn|z|
\quad\text{for }s=\frac1n.
\]
Writing
\(w_{k+1}-\mathsf d(t_{k+1}-)
=[\Phi_k^{1/n}(w_k)+f_k-\Phi_k^{1/n}(w_k+f_k)]
+[\Phi_k^{1/n}(w_k+f_k)-\Phi_k^{1/n}(\mathsf d(t_k-)+f_k)]\)
by \eqref{eq:terminal-quotient-flow} and \eqref{eq:wrec}, and applying
the last two estimates to the two brackets, we obtain
\[
|w_{k+1}-\mathsf d(t_{k+1}-)|
\leq
\left(1+\frac{C}{n}\right)
|w_k-\mathsf d(t_k-)|
+\frac{C}{n}|f_k|.
\]
Since \(w_0=\mathsf d(0-)=0\), the discrete Gronwall inequality, together
with \eqref{eq:fkbound}, gives
\begin{equation}
|w_n-\mathsf d(1)|
\leq
\frac{C}{n}\sum_{k=0}^{n-1}|f_k|
\leq
\frac{CM}{n}.
\label{eq:terminal-commutator-error}
\end{equation}

Finally, since
\(\Theta^{-1}(Y_n^n)=\Theta^{-1}(Z_n^n+d_n^n)\) and \(\Theta^{-1}\) is
\(\sigma_{\mathrm{max}}\)-Lipschitz, \eqref{eq:terminal-cutoff-endpoint},
\eqref{eq:terminal-euler-flow-error} and
\eqref{eq:terminal-commutator-error} give
\begin{equation}
\begin{aligned}
\left|
\Theta^{-1}(Y_n^n)-\Theta^{-1}(Z_n^n)
-\widehat U_1^{n,\varepsilon}
\right|
&=
\left|
\Theta^{-1}(Z_n^n+d_n^n)
-\Theta^{-1}(Z_n^n+\mathsf d(1))
\right|\\
&\leq
\sigma_{\mathrm{max}}
\left(
|d_n^n-w_n|+|w_n-\mathsf d(1)|
\right)
\leq
C\left(
\sum_{k=0}^{n-1}|\eta_k|+\frac{M+1}n
\right).
\end{aligned}
\label{eq:terminal-cutoff-endpoint-proof}
\end{equation}
Combining this estimate with
\eqref{eq:terminal-remainder-sum} proves
\eqref{eq:terminal-quotient-comparison}.
\end{proof}

\section{\texorpdfstring{Proof of Theorem~\ref{thm:main}}{Proof of Theorem 1.1}}
\label{sec:terminal-ldp}

We begin with the LDP for the quadratic-variation measures
\(\mu_n\) in \eqref{eq:tag-10}. It provides the function \(\ell\) used in
the functional \(I_{\rm sv}\) below, and the tail bound \eqref{eq:tag-13} for
the total quadratic variation, which is used repeatedly.

\begin{proposition}
\label{prop:qv}
The sequence $\{\mu_n\}$ satisfies a full LDP on
$\Mplus([0,1])$, equipped with the weak topology
\eqref{eq:weaktoM}, at speed $n$ and with good rate function $I_{\rm qv}$
given by
\begin{equation}
I_{\rm qv}(\mu)=\int_0^1\ell(v(t))\dd t+\frac12\mu^s([0,1]),\quad
\mu=v(t)\,\dd t+\mu^s\in\Mplus([0,1]),
\label{eq:tag-11}
\end{equation}
where $v(t)\dd t$ is the absolutely continuous part of $\mu$ with respect to $\dd t$, $\mu^s$ is its singular part, and
\begin{equation}\label{eq:ellv}
\ell(v)=
\begin{cases}
\frac12(v-1-\log v),&v>0,\\
+\infty,&v\leq0.
\end{cases}
\end{equation}
Moreover,
\begin{equation}
I_{\rm qv}(\mu)
\geq
\ell\bigl(\mu([0,1])\bigr),
\qquad\mu\in\Mplus([0,1]),
\label{eq:tag-12}
\end{equation}
and, for every $M>1$ and $n\geq1$,
\begin{equation}
\mathbb P\{\mu_n([0,1])\geq M\}
\leq e^{-n\ell(M)}.
\label{eq:tag-13}
\end{equation}
\end{proposition}
\begin{proof}[Proof of Proposition~\ref{prop:qv}]
We specialize \cite[Remark~3.5]{DGS2017}
to unit volatility, zero drift, and zero correlation. Then the first
realized quadratic-variation coordinate there is
\[
Q_{1,t}^n
=
\sum_{j=1}^{\lfloor nt\rfloor}
\bigl(W_{j/n}-W_{(j-1)/n}\bigr)^2
=
\mu_n([0,t]),
\]
and, under the identification of bounded-variation paths
with their Stieltjes measures, the weak-* topology used in
\cite{DGS2017} is the weak topology in \eqref{eq:weaktoM}.
Hence the cited result gives the LDP for \(\{\mu_n\}\) with
rate function \eqref{eq:tag-11}.

To prove \eqref{eq:tag-12}, write
\(\mu=v(t)\,\dd t+\mu^s\), and set \(r=\int_0^1v\,\dd t\) and
\(s=\mu^s([0,1])\). Since \(\ell\) is convex, Jensen's inequality gives
\(\int_0^1\ell(v(t))\,\dd t\geq\ell(r)\). Since
\(\ell(r+s)-\ell(r)=\frac12\bigl(s-\log(1+s/r)\bigr)\leq\frac s2\) for
\(r>0\), and \(\ell(0)=+\infty\), we obtain
\(I_{\rm qv}(\mu)\geq\ell(r)+\frac s2\geq\ell(r+s)=\ell(\mu([0,1]))\).
Since \(\ell(m)\to\infty\) as \(m\to\infty\), \eqref{eq:tag-12} shows that
the sublevel sets of \(I_{\rm qv}\) have bounded total mass. They are
therefore relatively compact by Prokhorov's theorem, and closed by the
lower semicontinuity of \(I_{\rm qv}\). Hence \(I_{\rm qv}\) is good.

Finally, \(n\mu_n([0,1])\) is a sum of \(n\) independent squared standard
normal variables. For \(M>1\), let \(\theta:=\frac12(1-M^{-1})\in(0,\frac12)\).
Since \(\mathbb E\exp\{\theta n(\Delta_kW)^2\}=(1-2\theta)^{-1/2}\),
Markov's inequality gives
\[
\mathbb P\{\mu_n([0,1])\geq M\}
\leq e^{-n\theta M}(1-2\theta)^{-n/2}
=\exp\Bigl\{-n\Bigl[\theta M+\frac12\log(1-2\theta)\Bigr]\Bigr\}
=e^{-n\ell(M)},
\]
which is \eqref{eq:tag-13}.
\end{proof}

\begin{remark*}
The LDP in Proposition~\ref{prop:qv} also follows from
\cite[Theorem~3.1]{Najim} by viewing \(\mu_n\) as a weighted empirical
measure. Indeed,
\(\mu_n=\frac1n\sum_{k=0}^{n-1}(\sqrt n\,\Delta_kW)^2\delta_{t_{k+1}}\),
where for each \(n\) the weights are independent \(\chi^2_1\) random
variables. The empirical distribution of the grid points converges to
Lebesgue measure, which charges every nonempty open subset of \([0,1]\).
The log moment-generating function of each weight is
\(\Lambda(\theta)=-\frac12\log(1-2\theta)\) for \(\theta<\frac12\).
Its convex conjugate is \(\ell\) in \eqref{eq:ellv}, and the recession
function of \(\ell\) equals \(v/2\) for \(v\geq0\), giving the two terms
in \eqref{eq:tag-11}. Since every \(\mu_n\) is nonnegative and the rate
function is infinite outside the closed set \(\Mplus([0,1])\), the LDP
restricts to this set \cite[Lemma~4.1.5]{DemboZeitouni}.
\end{remark*}

Fix a cutoff \(0<\varepsilon\leq1\). For each
\(N\in\mathbb N_0\), let \(\mathcal C_{\varepsilon,N}\) be the set of all
quadruples \((A,B,\mathfrak m,\xi_f)\) such that \((A,B)\in\mathcal Q_\sigma\),
\(\xi_f\in K_\sigma\), and \(\mathfrak m\) is an ordered list of exactly
\(N\) marks as in \eqref{eq:ordered-mark-list}. For \(N=0\), the list is empty.
The fixed-cutoff enhanced control space is the disjoint union
\begin{equation}
\mathcal C_\varepsilon
:=\bigsqcup_{N=0}^{\infty}\mathcal C_{\varepsilon,N},
\label{eq:enhanced-control-space}
\end{equation}
and the empirical enhanced control is
\begin{equation}
\zeta_n^\varepsilon
=
\bigl(
A_n,B_n^\varepsilon,\mathfrak m_n^\varepsilon,\iota(Z_n^n)
\bigr)\in\mathcal C_\varepsilon.
\label{eq:empirical-control}
\end{equation}
The terminal state tag \(\xi_f\in K_\sigma\) is retained to convert the
terminal discrepancy from Lamperti coordinates back to the original
coordinates. On \(\mathcal C_{\varepsilon,N}\), we use the
topology
induced by weak convergence of \(A\) and \(B\) and coordinatewise
convergence of the marks and terminal tag. For each \(N\), \(\mathcal C_{\varepsilon,N}\) is
a closed subset of the Polish space
\(\mathcal Q_\sigma\times([0,1]\times K_\sigma\times\R)^N\times K_\sigma\),
since the constraints \(t_1\leq\cdots\leq t_N\) and
\(|\delta_j|\geq\varepsilon\) are closed. Thus \(\mathcal C_{\varepsilon,N}\) is Polish for every \(N\),
and \(\mathcal C_\varepsilon\) is Polish in the disjoint union topology.

This topology is induced by the following metric. For \(\zeta,\zeta'\in\mathcal C_{\varepsilon,N}\) with
\[
\zeta=(A,B,\mathfrak m,\xi_f),\qquad
\zeta'=(A',B',\mathfrak m',\xi_f'),
\qquad
\mathfrak m=((t_j,\xi_j,\delta_j))_{j=1}^N,\qquad
\mathfrak m'=((t_j',\xi_j',\delta_j'))_{j=1}^N,
\]
set
\begin{equation}\label{eq:joint-control-metric}
\begin{aligned}
d_{\mathcal C_\varepsilon}(\zeta,\zeta')
:=\min\Bigg\{1,\;&
d_{\mathcal Q_\sigma}\bigl((A,B),(A',B')\bigr)
+d_{K_\sigma}(\xi_f,\xi_f')\\
&+\sum_{j=1}^N\Bigl(
|t_j-t_j'|+d_{K_\sigma}(\xi_j,\xi_j')
+|\delta_j-\delta_j'|\Bigr)\Bigg\},
\end{aligned}
\end{equation}
and for \(\zeta\in\mathcal C_{\varepsilon,N}\) and
\(\zeta'\in\mathcal C_{\varepsilon,N'}\) with \(N\ne N'\), set
\(d_{\mathcal C_\varepsilon}(\zeta,\zeta'):=1\).
The cutoff-free control space is
\(
\mathcal C=\bigcup_{0<\eta\leq1}\mathcal C_\eta.
\)

The action of an enhanced control is
\eqref{eq:IAB}. Its first term extends the rate function \(I_{\rm qv}\) in
\eqref{eq:tag-11} to pairs \((A,B)\in\mathcal Q_\sigma\), with the Lebesgue
measure \(\dd t\) replaced by \(A\). For \((A,B)\in\mathcal Q_\sigma\), let
\[
B=vA+B^\perp,
\qquad B^\perp\perp A,
\]
be the Lebesgue decomposition of \(B\) with respect to \(A\), where
\(v\geq0\) is the density of the absolutely continuous part.
We define the functional \(I_{\rm sv}\) by
\begin{equation}
I_{\rm sv}(A,B)
:=
\int_{K_\sigma^0}\ell(v)\,\dd A
+\frac12 B^\perp(K_\sigma^0),\quad (A,B)\in\mathcal{Q}_\sigma.
\label{eq:Isv-definition}
\end{equation}
This definition is independent of the \(A\)-a.e.\ version of \(v\).
If \(K_\sigma\) consists of a single point, as for constant coefficients, then \(A=\dd t\)
and \(I_{\rm sv}(A,B)=I_{\rm qv}(B)\). The second term
\(\frac12\sum_j\delta_j^2\) in \eqref{eq:IAB} is the Gaussian cost of the
large increments.
\begin{theorem}[Fixed-cutoff enhanced-control LDP]
\label{thm:cutoffldp}
For every fixed \(0<\varepsilon\leq1\), the random controls
\(\zeta_n^\varepsilon\) satisfy a full LDP on
\(\mathcal C_\varepsilon\) at speed \(n\), with good rate function
\(I_\varepsilon:=\left.I\right|_{\mathcal C_\varepsilon}\),
where \(I\) is defined~in~\eqref{eq:IAB}.
\end{theorem}

Theorem~\ref{thm:cutoffldp} is proved in
Section~\ref{sec:fixed-cutoff-ldp}. To pass to the endpoint, define
\[
G:\
\mathcal C_\varepsilon\ni(A,B,\mathfrak m,\xi_f)
\longmapsto
\bar{\mathfrak g}\bigl(\xi_f,d_{A,B,\mathfrak m}(1)\bigr)\in\R.
\]
By \eqref{eq:empirical-control} and
\eqref{eq:terminal-cutoff-endpoint},
\(\widehat U_1^{n,\varepsilon}=G(\zeta_n^\varepsilon)\).

\begin{corollary}
\label{cor:fixedcutoff}
For every fixed \(0<\varepsilon\leq1\),
\(\widehat U_1^{n,\varepsilon}\) satisfies a full LDP on \(\R\) at speed
\(n\), with good rate function \(J_\varepsilon\) defined in
\eqref{eq:tag-5}.
\end{corollary}
\begin{proof}
By Proposition~\ref{prop:wellposed}, \(G\) is continuous on
\(\mathcal C_{\varepsilon,N}\) for every \(N\). These sets
form an open cover of \(\mathcal C_\varepsilon\), so \(G\) is
continuous on \(\mathcal C_\varepsilon\). Since
\(\widehat U_1^{n,\varepsilon}=G(\zeta_n^\varepsilon)\), the contraction
principle \cite[Theorem~4.2.1]{DemboZeitouni}, applied to the LDP in
Theorem~\ref{thm:cutoffldp}, yields the LDP for
\(\widehat U_1^{n,\varepsilon}\) with the good rate function
\(e\mapsto\inf\{I_\varepsilon(\zeta):\zeta\in\mathcal C_\varepsilon,\ G(\zeta)=e\}\),
which is \(J_\varepsilon\).
\end{proof}

\begin{proposition}[{\cite[Exercise~4.2.29(a), p.~136]{DemboZeitouni}}]
\label{prop:approximation-principle}
Let \((E,\rho)\) be a Polish space and let \(V_n\) and
\(V_n^{(m)}\), \(n,m\ge1\), be \(E\)-valued random variables.
Suppose that for each fixed \(m\), \((V_n^{(m)})_{n\ge1}\) satisfies a
full LDP at speed \(n\) with good rate function \(\mathcal R_m\).
Suppose also that there are couplings of \(V_n\) and \(V_n^{(m)}\) for which
\begin{equation}
\lim_{m\to\infty}\limsup_{n\to\infty}\frac1n
\log\mathbb P\!\left\{\rho(V_n,V_n^{(m)})>\delta\right\}
=-\infty,\qquad \delta>0.
\label{eq:terminal-scalar-approximation}
\end{equation}
Then \((V_n)_{n\ge1}\) is exponentially tight at speed \(n\) and satisfies
a full LDP at that speed with good rate function
\begin{equation}
\mathcal R(y):=\sup_{r>0}\liminf_{m\to\infty}
\inf_{\rho(z,y)<r}\mathcal R_m(z),\qquad y\in E.
\label{eq:terminal-scalar-relaxation}
\end{equation}
\end{proposition}

For any sequence \(\varepsilon_m\downarrow0\) in \((0,1]\),
we apply Proposition~\ref{prop:approximation-principle} on \(\R\) with
\(V_n=U_1^n\) and \(V_n^{(m)}=\widehat U_1^{n,\varepsilon_m}\).
Corollary~\ref{cor:fixedcutoff} gives the fixed-cutoff LDPs.
To verify \eqref{eq:terminal-scalar-approximation},
Proposition~\ref{prop:expapprox} uses the following estimate for
\(|Y_1-Y_n^n|\).

%
\begin{lemma}
\label{lem:bridge}
Let \(Y\) solve \eqref{eq:tag-9}, and let \(Y^n\) be defined
by \eqref{eq:aux-lamperti-recursion}. Then, for every \(\eta>0\),
\begin{equation}
\limsup_{n\to\infty}
\frac1n\log
\mathbb P
\left\{
|Y_1-Y_n^n|>\eta
\right\}
=-\infty.
\label{eq:bridge-superexponential}
\end{equation}
\end{lemma}

\begin{proof}
Set \(e_k:=Y_{t_k}-Y_k^n\) and
\(S_k:=\sup_{t_k\leq t\leq t_{k+1}}|W_t-W_{t_k}|\).
By \eqref{eq:tag-9}, \eqref{eq:aux-lamperti-recursion}
and the Lipschitz continuity of \(\beta\), for \(0\leq k<n\),
\begin{equation}\label{eq:Ytk+1}
|e_{k+1}|
\leq
\left(1+\frac{L_\beta}{n}\right)|e_k|
+
L_\beta\int_{t_k}^{t_{k+1}}
|Y_s-Y_{t_k}|\,\dd s.
\end{equation}
Since \(Y_t-Y_{t_k}=W_t-W_{t_k}+\int_{t_k}^t\beta(Y_s)\,\dd s\)
by \eqref{eq:tag-9}, and \(\beta\) is bounded,
\(\sup_{t_k\leq t\leq t_{k+1}}|Y_t-Y_{t_k}|\leq S_k+\|\beta\|_\infty/n\).
Since \(e_0=0\), inserting this bound into \eqref{eq:Ytk+1} and applying
the discrete Gronwall inequality give
\begin{equation}
|Y_1-Y_n^n|
\leq
C_1\left(\frac1n\sum_{k=0}^{n-1}S_k+\frac1n\right)
\label{eq:bridge-grid-error}
\end{equation}
for a constant \(C_1\) depending only on \(L_\beta\) and
\(\|\beta\|_\infty\).

By the independence and scaling of Brownian increments, the
variables \(\sqrt n\,S_0,\ldots,\sqrt n\,S_{n-1}\) are i.i.d.\ with the law
of \(S:=\sup_{0\leq s\leq1}|W_s|\). By the reflection principle,
\(\mathbb P\{S>x\}\leq4\,\mathbb P\{W_1>x\}\leq2e^{-x^2/2}\), and hence
\(\kappa:=\mathbb E\exp(S^2/4)<\infty\).

Fix \(\eta>0\), and let \(n>2C_1/\eta\). By
\eqref{eq:bridge-grid-error}, the event \(\{|Y_1-Y_n^n|>\eta\}\) implies
\(\sum_k\sqrt n\,S_k>c_\eta n^{3/2}\) with \(c_\eta:=\eta/(2C_1)\). By the
Cauchy--Schwarz inequality, it therefore implies
\(\sum_k nS_k^2>c_\eta^2n^2\). Independence and the exponential Markov
inequality yield
\[
\mathbb P\left\{|Y_1-Y_n^n|>\eta\right\}
\leq
\mathbb P\Bigl\{\sum_{k=0}^{n-1}nS_k^2>c_\eta^2n^2\Bigr\}
\leq
e^{-c_\eta^2n^2/4}\,\kappa^n.
\]
Hence
\(\frac1n\log\mathbb P\{|Y_1-Y_n^n|>\eta\}\leq-c_\eta^2n/4+\log\kappa\to-\infty\),
which proves \eqref{eq:bridge-superexponential}.
\end{proof}
\begin{proposition}[Exponentially good approximation]
\label{prop:expapprox}
For every \(\delta>0\),
\begin{equation}
\lim_{\varepsilon\downarrow0}
\limsup_{n\to\infty}
\frac1n
\log
\mathbb P
\left\{
|U_1^n-\widehat U_1^{n,\varepsilon}|>\delta
\right\}
=
-\infty.
\label{eq:terminal-exponential-approximation}
\end{equation}
\end{proposition}

\begin{proof}[Proof of Proposition~\ref{prop:expapprox}]
Since \(\Theta^{-1}(Y_n^n)=\Theta^{-1}(Z_n^n+d_n^n)\),
\eqref{eq:U1n} gives
\[
|U_1^n-\widehat U_1^{n,\varepsilon}|
\leq
|\Theta^{-1}(Y_1)-\Theta^{-1}(Y_n^n)|+
\left|
\Theta^{-1}(Y_n^n)-\Theta^{-1}(Z_n^n)
-\widehat U_1^{n,\varepsilon}
\right|.
\]
Fix \(M>1\), and let \(C\) be the constant in
Lemma~\ref{lem:quotient}. Choose \(\varepsilon_0\in(0,1]\) such that
\(C\varepsilon_0M<\delta/4\). For every fixed
\(0<\varepsilon\leq\varepsilon_0\), the remaining terms on the right-hand
side of \eqref{eq:terminal-quotient-comparison} tend to zero as
\(n\to\infty\). Hence, for all sufficiently large \(n\), the second term
above is at most \(\delta/2\) on \(\{\mu_n([0,1])\leq M\}\). Since
\(\Theta^{-1}\) is \(\sigma_{\mathrm{max}}\)-Lipschitz, it follows that
\[
\mathbb P
\left\{
|U_1^n-\widehat U_1^{n,\varepsilon}|>\delta
\right\}
\leq
\mathbb P\{\mu_n([0,1])>M\}
+
\mathbb P
\left\{
|Y_1-Y_n^n|>\frac{\delta}{2\sigma_{\mathrm{max}}}
\right\}.
\]
By \eqref{eq:tag-13} and Lemma~\ref{lem:bridge}, we
obtain, for every \(0<\varepsilon\leq\varepsilon_0\),
\[
\limsup_{n\to\infty}
\frac1n
\log
\mathbb P
\left\{
|U_1^n-\widehat U_1^{n,\varepsilon}|>\delta
\right\}
\leq
-\ell(M).
\]
Hence the upper limit as \(\varepsilon\downarrow0\) of the
left-hand side is at most \(-\ell(M)\). Since \(M>1\) is arbitrary and
\(\ell(M)\to\infty\) as \(M\to\infty\), this proves
\eqref{eq:terminal-exponential-approximation}.
\end{proof}

%

%
\begin{lemma}[Absorption of small marks]
\label{lem:absorb}
Let \(0<\varepsilon\leq1\) satisfy
\(\|\sigma'\|_\infty\varepsilon<1\).
For every
\(\zeta=(A,B,\mathfrak m,\xi_f)\in\mathcal C\) with
\(I(\zeta)<\infty\), there exist
\(\widetilde B\in\Mplus(K_\sigma^0)\) and a list
\(\widetilde{\mathfrak m}\) such that
$
(A,\widetilde B,\widetilde{\mathfrak m},\xi_f)
\in\mathcal C_\varepsilon,
$
$
\bar{\mathfrak g}\bigl(
\xi_f,d_{A,\widetilde B,\widetilde{\mathfrak m}}(1)
\bigr)
=
\bar{\mathfrak g}\bigl(
\xi_f,d_{A,B,\mathfrak m}(1)
\bigr)$,
and
\begin{equation}\label{eq:Abs}
I_\varepsilon(A,\widetilde B,\widetilde{\mathfrak m},\xi_f)
\leq
(1-\|\sigma'\|_\infty\varepsilon)^{-2}
I(\zeta).
\end{equation}
\end{lemma}
\begin{proof}
We replace each mark with \(|\delta_j|<\varepsilon\) by an
atom of \(B\) that produces the same jump in \eqref{eq:controlled-error}.
The jump produced by a mark \((t,\xi,\delta)\) is
\begin{equation}
h(\xi,\delta):=\delta-\bar{\mathfrak j}(\xi,\delta).
\label{eq:jump-remainder}
\end{equation}
We first bound \(h\). Set
$
\alpha_+:=\max_{\xi\in K_\sigma}\max\{\bar\alpha(\xi),0\}$ and
$\alpha_-:=\max_{\xi\in K_\sigma}\max\{-\bar\alpha(\xi),0\}$.
Let \(y\in\R\) and \(x=\Theta^{-1}(y)\). Since
\(\bar{\mathfrak j}(\iota(y),\delta)=\Theta(x+\sigma(x)\delta)-\Theta(x)\)
and \(\Theta''=-\sigma'/\sigma^2\), Taylor's formula
with integral remainder gives
\[
h(\iota(y),\delta)
=\delta^2\int_0^1(1-s)\sigma'(x+s\,\sigma(x)\delta)
\left(\frac{\sigma(x)}
{\sigma(x+s\,\sigma(x)\delta)}\right)^2\dd s.
\]
For \(|\delta|<\varepsilon\), we have
\(
\sigma(x+s\,\sigma(x)\delta)
\geq
\sigma(x)(1-\|\sigma'\|_\infty\varepsilon)>0
\)
and \(-\alpha_-\leq\sigma'\leq\alpha_+\). Since
\(\int_0^1(1-s)\,\dd s=\frac12\), it follows that
\begin{equation}\label{eq:alphaest}
-\frac12\alpha_-(1-\|\sigma'\|_\infty\varepsilon)^{-2}\delta^2
\leq h(\xi,\delta)
\leq
\frac12\alpha_+(1-\|\sigma'\|_\infty\varepsilon)^{-2}\delta^2,
\qquad |\delta|<\varepsilon,
\end{equation}
for \(\xi=\iota(y)\). Since \(\iota(\R)\) is dense in
\(K_\sigma\) and \(h(\cdot,\delta)\) is continuous,
\eqref{eq:alphaest} holds for every \(\xi\in K_\sigma\).

By compactness of \(K_\sigma\) and continuity of
\(\bar\alpha\), there exist \(\xi_+\) with \(\bar\alpha(\xi_+)=\alpha_+\)
if \(\alpha_+>0\), and \(\xi_-\) with \(\bar\alpha(\xi_-)=-\alpha_-\) if
\(\alpha_->0\). For each mark \((t_j,\xi_j,\delta_j)\) with
\(|\delta_j|<\varepsilon\), define
\[
\nu_j:=
\begin{cases}
\displaystyle
\frac{2h(\xi_j,\delta_j)}{\alpha_+}\,
\delta_{(t_j,\xi_+)},
& h(\xi_j,\delta_j)>0,\\[2mm]
\displaystyle
\frac{2|h(\xi_j,\delta_j)|}{\alpha_-}\,
\delta_{(t_j,\xi_-)},
& h(\xi_j,\delta_j)<0,\\[2mm]
0,& h(\xi_j,\delta_j)=0.
\end{cases}
\]
By \eqref{eq:alphaest}, \(h(\xi_j,\delta_j)>0\) forces
\(\alpha_+>0\), and \(h(\xi_j,\delta_j)<0\) forces \(\alpha_->0\). Hence
\(\nu_j\) is well defined, and \eqref{eq:alphaest} gives
\begin{equation}
\nu_j(K_\sigma^0)
\leq
(1-\|\sigma'\|_\infty\varepsilon)^{-2}\delta_j^2,
\qquad
\frac12\int_{K_\sigma^0}
\bar\alpha(\xi)\,\nu_j(\dd t,\dd\xi)
=
h(\xi_j,\delta_j).
\label{eq:small-mark-mass}
\end{equation}

Set \(\widetilde B:=B+\nu\) with
\(\nu:=\sum_{j:\,|\delta_j|<\varepsilon}\nu_j\), and
\[
\widetilde{\mathfrak m}
:=
\bigl((t_j,\xi_j,\delta_j)\bigr)_{j:\,|\delta_j|\geq\varepsilon}.
\]
Then \((A,\widetilde B,\widetilde{\mathfrak m},\xi_f)\in
\mathcal C_\varepsilon\). By the second identity in
\eqref{eq:small-mark-mass}, the jump \(h(\xi_j,\delta_j)\) of each deleted
mark at \(t_j\) equals the contribution of \(\nu_j\) to the \(B\)-term in
\eqref{eq:controlled-error}. Hence \((A,B,\mathfrak m)\) and
\((A,\widetilde B,\widetilde{\mathfrak m})\) give the same equation
\eqref{eq:controlled-error}, and uniqueness yields
\(d_{A,\widetilde B,\widetilde{\mathfrak m}}=d_{A,B,\mathfrak m}\). In
particular, the two terminal errors coincide.

It remains to compare the actions. Since \(\nu\) is
concentrated on the finitely many times of the deleted marks, while the
time marginal of \(A\) is Lebesgue measure, \(\nu\perp A\). Thus
\(I_{\rm sv}(A,B+\nu)=I_{\rm sv}(A,B)+\frac12\nu(K_\sigma^0)\).
By \eqref{eq:IAB} and the first bound in
\eqref{eq:small-mark-mass},
\[
\begin{aligned}
I_\varepsilon(A,\widetilde B,\widetilde{\mathfrak m},\xi_f)
&=
I(\zeta)
-\frac12\sum_{j:\,|\delta_j|<\varepsilon}\delta_j^2
+\frac12\nu(K_\sigma^0)\\
&\leq
I(\zeta)
+
\bigl((1-\|\sigma'\|_\infty\varepsilon)^{-2}-1\bigr)
\frac12\sum_{j:\,|\delta_j|<\varepsilon}\delta_j^2.
\end{aligned}
\]
Since
\(
\frac12\sum_{j:\,|\delta_j|<\varepsilon}\delta_j^2
\leq
I(\zeta),
\)
this proves \eqref{eq:Abs}.
\end{proof}
\begin{proposition}
\label{prop:rate}
For every \(e\in\R\) and every sequence
\(\varepsilon_m\downarrow0\),
\begin{equation}\label{eq:June}
J_{\rm un}(e)
=
\lim_{\varepsilon\downarrow0}J_\varepsilon(e)
=
\sup_{\rho>0}\liminf_{m\to\infty}
\inf_{|y-e|<\rho}J_{\varepsilon_m}(y).
\end{equation}
\end{proposition}

\begin{proof}
For \(0<\eta<\varepsilon\leq1\), we have
\(\mathcal C_\varepsilon\subset\mathcal C_\eta\) and
\(I_\eta(\zeta)=I_\varepsilon(\zeta)\) for every \(\zeta\in\mathcal C_\varepsilon\). Hence
\(J_\eta(e)\leq J_\varepsilon(e)\), so \(J_\varepsilon(e)\) is
nonincreasing as \(\varepsilon\downarrow0\).
Since
\(\mathcal C=\bigcup_{0<\varepsilon\leq1}\mathcal C_\varepsilon\),
the definition \eqref{eq:tag-6} gives
\(
J_{\rm un}(e)
=
\inf_{0<\varepsilon\leq1}J_\varepsilon(e)
=\lim_{\varepsilon\downarrow0}J_\varepsilon(e).
\)
This proves the first equality in \eqref{eq:June}.

To prove the second equality, we first compare
\(J_\varepsilon\) with \(J_{\rm un}\). Let \(\|\sigma'\|_\infty\varepsilon<1\).
Lemma~\ref{lem:absorb} gives, for
every \((A,B,\mathfrak m,\xi_f)\in\mathcal C\) with terminal error \(e\)
and finite action, a control
\((A,\widetilde B,\widetilde{\mathfrak m},\xi_f)\in\mathcal C_\varepsilon\)
with the same terminal error and
\(
I_\varepsilon(A,\widetilde B,\widetilde{\mathfrak m},\xi_f)
\leq
(1-\|\sigma'\|_\infty\varepsilon)^{-2}
I(A,B,\mathfrak m,\xi_f).
\)
Taking the infimum over such controls yields
\begin{equation}
J_{\rm un}(e)
\leq J_\varepsilon(e)
\leq
(1-\|\sigma'\|_\infty\varepsilon)^{-2}J_{\rm un}(e),
\qquad e\in\R.
\label{eq:cutoff-comparison}
\end{equation}
We next show that \(J_{\rm un}\) is lower semicontinuous.
Let \(e_j\to e\), and fix \(\varepsilon\) with
\(\|\sigma'\|_\infty\varepsilon<1\).
Since \(J_\varepsilon\) is lower semicontinuous by
Corollary~\ref{cor:fixedcutoff}, \eqref{eq:cutoff-comparison} gives
\[
\liminf_{j\to\infty}J_{\rm un}(e_j)
\geq
(1-\|\sigma'\|_\infty\varepsilon)^2
\liminf_{j\to\infty}J_\varepsilon(e_j)\geq
(1-\|\sigma'\|_\infty\varepsilon)^2J_\varepsilon(e)\geq
(1-\|\sigma'\|_\infty\varepsilon)^2J_{\rm un}(e).
\]
Letting \(\varepsilon\downarrow0\) proves the lower semicontinuity
of \(J_{\rm un}\).
Finally, since \(J_{\varepsilon_m}\) is nonincreasing in \(m\)
and converges pointwise to \(J_{\rm un}\), for every
\(\rho>0\),
\[
\lim_{m\to\infty}\inf_{|y-e|<\rho}J_{\varepsilon_m}(y)
=\inf_m\inf_{|y-e|<\rho}J_{\varepsilon_m}(y)=\inf_{|y-e|<\rho}\inf_mJ_{\varepsilon_m}(y)
=\inf_{|y-e|<\rho}J_{\rm un}(y).
\]
By the lower semicontinuity of \(J_{\rm un}\), the right-hand
side increases to \(J_{\rm un}(e)\) as \(\rho\downarrow0\). Taking
\(\sup_{\rho>0}\) proves the second equality in \eqref{eq:June}.
\end{proof}

\begin{proof}[Proof of Theorem~\ref{thm:main}]
\medskip\noindent\textit{Step 1.} Fix a sequence
\(\varepsilon_m\downarrow0\) in \((0,1]\). For the approximations
\(V_n^{(m)}=\widehat U_1^{n,\varepsilon_m}\), Corollary~\ref{cor:fixedcutoff}
gives the fixed-cutoff LDP with good rate function \(J_{\varepsilon_m}\),
and Proposition~\ref{prop:expapprox} verifies
\eqref{eq:terminal-scalar-approximation} with \(V_n=U_1^n\).
Proposition~\ref{prop:approximation-principle} therefore yields a full LDP
for \(U_1^n\) at speed \(n\). Its good rate function is given by
\eqref{eq:terminal-scalar-relaxation} with \(\mathcal R_m=J_{\varepsilon_m}\)
and equals \(J_{\rm un}\) by Proposition~\ref{prop:rate}.

\medskip\noindent\textit{Step 2.} We show that \(J_{\rm un}(e)=0\) if
and only if \(e=0\).

Let \(\zeta\mathrel{=}(A,B,\mathfrak m,\xi_f)\in\mathcal C\). Since
\(\bar\gamma(\xi)=\bar\beta(\xi,0)+\frac12\bar\alpha(\xi)\)
for \(\xi\in K_\sigma\), \eqref{eq:controlled-error} can be rewritten as
\begin{equation}
\begin{aligned}
d(t)
&=
\int_{[0,t]\times K_\sigma}
\bigl[
\bar\beta(\xi,d(r))-\bar\beta(\xi,0)
\bigr]
A(\dd r,\dd\xi)\\
&\quad+\frac12
\int_{[0,t]\times K_\sigma}
\bar\alpha(\xi)(B-A)(\dd r,\dd\xi)
+\sum_{j:t_j\leq t}
\bigl[
\delta_j-\bar{\mathfrak j}(\xi_j,\delta_j)
\bigr].
\end{aligned}
\label{eq:controlled-error-difference}
\end{equation}

We first show that \(J_{\rm un}(0)=0\). Fix
\(\xi_f\in K_\sigma\), and let \(A\in\mathcal P(K_\sigma^0)\) have time
marginal \(\dd t\), for instance \(A=\dd t\otimes\delta_{\xi_f}\). For
\(B=A\) and \(\mathfrak m=\varnothing\), the last two terms in
\eqref{eq:controlled-error-difference} vanish, so \(d\equiv0\) solves
\eqref{eq:controlled-error-difference}. By uniqueness in
Proposition~\ref{prop:wellposed}, \(d_{A,A,\varnothing}\equiv0\), and hence
\(\bar{\mathfrak g}(\xi_f,d_{A,A,\varnothing}(1))=\bar{\mathfrak g}(\xi_f,0)=0\).
Moreover, \(v=1\) and \(B^\perp=0\), so
\(I(A,A,\varnothing,\xi_f)=\int_{K_\sigma^0}\ell(1)\,\dd A=0\). Thus
\eqref{eq:tag-6} yields
\(
0\leq J_{\rm un}(0)\leq I(A,A,\varnothing,\xi_f)=0.
\)

Conversely, we bound the terminal error by the action.
Let \(I(\zeta)<\infty\), and let \(B=vA+B^\perp\) be the Lebesgue
decomposition. We bound the three terms in
\eqref{eq:controlled-error-difference}. Since \(\bar\beta\) is globally
\(L_\beta\)-Lipschitz in its second variable, the first term is
at most \(L_\beta\int_0^t|d(r)|\,\dd r\), because \(A\) has Lebesgue time
marginal. Since \(B-A=(v-1)A+B^\perp\) and \(\bar\alpha\) is bounded,
the second term is at most
\(
\frac12\|\bar\alpha\|_\infty
\bigl[
\int_{K_\sigma^0}|v-1|\,\dd A
+
B^\perp(K_\sigma^0)
\bigr].
\)
By \eqref{eq:lamperti-global-quadratic}, the density of
\(\iota(\R)\) in \(K_\sigma\), and the continuity of \(\bar{\mathfrak j}\) in
Lemma~\ref{lem:compactification}, we have
\(|\delta-\bar{\mathfrak j}(\xi,\delta)|\leq C\delta^2\) for all \(\xi\) and
\(\delta\), so the third term is at most \(C\sum_j\delta_j^2\). Hence
\begin{equation}
|d(t)|
\leq
L_\beta\int_0^t|d(r)|\,\dd r
+
C\left[
\int_{K_\sigma^0}|v-1|\,\dd A
+
B^\perp(K_\sigma^0)
+
\sum_j\delta_j^2
\right].
\label{eq:controlled-terminal-bound}
\end{equation}
Gronwall's inequality bounds \(|d_{A,B,\mathfrak m}(1)|\) by
\(C\) times the bracket in \eqref{eq:controlled-terminal-bound}. Since
\(\bar{\mathfrak g}(\xi,0)=0\) and \(\bar{\mathfrak g}(\xi,\cdot)\) is
\(\sigma_{\mathrm{max}}\)-Lipschitz by \eqref{eq:beta2}, we obtain
\begin{equation}
\left|
\bar{\mathfrak g}
\bigl(\xi_f,d_{A,B,\mathfrak m}(1)\bigr)
\right|
\leq
\sigma_{\mathrm{max}}|d_{A,B,\mathfrak m}(1)|
\leq
C\left[
\int_{K_\sigma^0}|v-1|\,\dd A
+
B^\perp(K_\sigma^0)
+
\sum_j\delta_j^2
\right].
\label{eq:terminal-action-bound}
\end{equation}

We next express the right-hand side in terms of the action. The elementary
estimate
\(
|v-1|
\leq
C\bigl(\sqrt{\ell(v)}+\ell(v)\bigr)
\)
for \(v>0\), together with \(A(K_\sigma^0)=1\) and the
Cauchy--Schwarz inequality, gives
\[
\int_{K_\sigma^0}|v-1|\,\dd A
\leq
C\left[
\left(
\int_{K_\sigma^0}\ell(v)\,\dd A
\right)^{1/2}
+
\int_{K_\sigma^0}\ell(v)\,\dd A
\right].
\]
Since each of \(\int\ell(v)\,\dd A\), \(\frac12B^\perp(K_\sigma^0)\)
and \(\frac12\sum_j\delta_j^2\) is at most \(I(\zeta)\) by
\eqref{eq:IAB} and \eqref{eq:Isv-definition},
\eqref{eq:terminal-action-bound} gives
\begin{equation}
\left|
\bar{\mathfrak g}
\bigl(\xi_f,d_{A,B,\mathfrak m}(1)\bigr)
\right|
\leq
C\left(
\sqrt{I(\zeta)}+I(\zeta)
\right).
\label{eq:terminal-small-action-bound}
\end{equation}

Now suppose \(J_{\rm un}(e)=0\). By \eqref{eq:tag-6}, there exist
\(\zeta_j=(A_j,B_j,\mathfrak m_j,\xi_{f,j})\in\mathcal C\) with
terminal error \(e\) and
\(\lim_{j\to\infty}I(\zeta_j)=0\).
By \eqref{eq:terminal-small-action-bound},
\(
|e|
\leq
C\big(
\sqrt{I(\zeta_j)}
+
I(\zeta_j)
\big)
\rightarrow0
\).
Hence \(e=0\).
\end{proof}

Since \(J_{\rm un}\) is good and vanishes only at zero, its
infimum over the closed set \(\{|e|\geq\delta\}\) is attained if it is
finite, and is therefore positive. Consequently, for every \(\delta>0\),
\[
c_\delta
:=
\inf_{|e|\geq\delta}J_{\rm un}(e)
>0,
\]
and the LDP upper bound gives
\[
\limsup_{n\to\infty}
\frac1n\log
\mathbb P\{|U_1^n|\geq\delta\}
\leq
-c_\delta<0.
\]
Thus the terminal Euler--Maruyama error concentrates exponentially
around zero at speed \(n\).

\section[The fixed-cutoff large deviation principle]{The fixed-cutoff large deviation principle}
\label{sec:fixed-cutoff-ldp}

Throughout this section, \(0<\varepsilon\leq1\) is fixed, and we prove Theorem~\ref{thm:cutoffldp}, namely the LDP for the
empirical controls \(\zeta_n^\varepsilon\) on \(\mathcal C_\varepsilon\)
with good rate function \(I_\varepsilon\). We follow the standard route
through a lower bound obtained by changes of measure, a local upper
bound, and exponential tightness.
Section~\ref{subsec:Isv-lsc} proves that the functional \(I_{\rm sv}\) is
lower semicontinuous and that its sublevel sets have bounded mass, which
leads to the goodness of \(I_\varepsilon\).
Section~\ref{subsec:recovery} constructs, for every control of
finite action, changes of measure under which \(\zeta_n^\varepsilon\)
concentrates near that control at an entropy cost close to its action.
Section~\ref{subsec:proof-cutoffldp} combines these results with an
exponential identity for squared Gaussian increments with predictable
weights, Lemma~\ref{lem:gaussianexp}, and proves the
theorem.

\subsection{Lower semicontinuity of \texorpdfstring{\(I_{\rm sv}\)}{I\_sv}}
\label{subsec:Isv-lsc}

The lower semicontinuity of \(I_{\rm sv}\) follows from a dual
representation. In Rockafellar's duality formula for convex integral
functionals of measures \cite[Theorem~5]{Rockafellar}, the recession
function determines the contribution of the singular part. Here
\(B^\perp\) is nonnegative, and the recession function of \(\ell\)
equals \(r/2\) for \(r\geq0\). The cited theorem assumes that the
reference measure has full support. Since \(A\) may vanish on a nonempty
open set, we give a direct proof.

\begin{lemma}[Dual representation of \(I_{\rm sv}\)]
\label{lem:Isv-dual}
For every \((A,B)\in\mathcal Q_\sigma\), the functional
\(I_{\rm sv}\) defined in \eqref{eq:Isv-definition} satisfies
\begin{equation}
I_{\rm sv}(A,B)
=
\sup_{\substack{\vartheta\in C(K_\sigma^0),
                 \sup\vartheta<1/2}}
\left\{
\int_{K_\sigma^0}\vartheta\,\dd B
+
\frac12\int_{K_\sigma^0}\log(1-2\vartheta)\,\dd A
\right\}.
\label{eq:Isv-dual}
\end{equation}
\end{lemma}

\begin{proof}
For \(r>0\), the concave function
\(q\mapsto qr+\frac12\log(1-2q)\) on \((-\infty,1/2)\) attains its maximum
at \(q=\frac12-\frac1{2r}\), with value \(\ell(r)\). For \(r=0\), its
supremum is \(+\infty=\ell(0)\). Hence the function $\ell$ in
\eqref{eq:ellv} satisfies
\begin{equation}
\ell(r)
=
\sup_{q<1/2}
\left\{
qr+\frac12\log(1-2q)
\right\},
\qquad r\geq0.
\label{eq:ell-pointwise-dual}
\end{equation}
Fix \((A,B)\in\mathcal Q_\sigma\), and let \(B=vA+B^\perp\) be the
Lebesgue decomposition used in \eqref{eq:Isv-definition}. We first
show that the right-hand side of \eqref{eq:Isv-dual} is at most
\(I_{\rm sv}(A,B)\). Let \(\vartheta\in C(K_\sigma^0)\) satisfy
\(\sup\vartheta<1/2\). By \eqref{eq:ell-pointwise-dual},
\(\vartheta v+\frac12\log(1-2\vartheta)\leq\ell(v)\) pointwise, and
\(\int\vartheta\,\dd B^\perp\leq\frac12B^\perp(K_\sigma^0)\) since \(\vartheta<\frac12\). Hence

\[
\begin{aligned}
\int_{K_\sigma^0}\vartheta\,\dd B
+\frac12
\int_{K_\sigma^0}
\log(1-2\vartheta)\,\dd A
&=
\int_{K_\sigma^0}
\left[
\vartheta v+\frac12\log(1-2\vartheta)
\right]\dd A
+
\int_{K_\sigma^0}\vartheta\,\dd B^\perp\\
&\le \int_{K_\sigma^0}\ell(v)\,\dd A
+\frac12B^\perp(K_\sigma^0)
=
I_{\rm sv}(A,B).
\end{aligned}
\]
Taking the supremum over all admissible \(\vartheta\) proves the claim.

For the reverse inequality, we restrict \(q\) in
\eqref{eq:ell-pointwise-dual} to a compact interval, realize the
restricted supremum by a measurable test function, and approximate this
function by continuous ones. For \(R\geq2\), let
\[
\ell_R(r)
:=
\sup_{-R\leq q\leq 1/2-1/R}
\left\{
qr+\frac12\log(1-2q)
\right\},
\qquad r\geq0.
\]
By concavity, this supremum is attained at the projection
\(q_R(r)\) of \(\frac12-\frac1{2r}\) onto \([-R,\frac12-\frac1R]\), with
\(q_R(0):=-R\), so that
\begin{equation}
\ell_R(r)=q_R(r)\,r+\frac12\log\bigl(1-2q_R(r)\bigr),
\qquad r\geq0.
\label{eq:ellR-attained}
\end{equation}
The map \(q_R\) is continuous on \([0,\infty)\). Moreover,
\(\ell_R\geq0\), since \(q=0\) is admissible, and \(\ell_R\) increases to
\(\ell\) as \(R\to\infty\) by \eqref{eq:ell-pointwise-dual}.
Since \(B^\perp\perp A\), there exists a Borel set
\(S\subset K_\sigma^0\) such that
$A(S)=0$ and $
B^\perp(K_\sigma^0\setminus S)=0$. Define the Borel function
\[
u_R(z)
:=
\begin{cases}
q_R(v(z)),&z\notin S,\\[2mm]
\displaystyle\frac12-\frac1R,&z\in S.
\end{cases}
\]
Since \(B=vA\) on \(K_\sigma^0\setminus S\), \(B=B^\perp\) on
\(S\) and \(A(S)=0\), the identity \eqref{eq:ellR-attained} yields
\begin{equation}\label{eq:Ksigu}
\int_{K_\sigma^0}u_R\,\dd B
+\frac12\int_{K_\sigma^0}\log(1-2u_R)\,\dd A
=
\int_{K_\sigma^0}\ell_R(v)\,\dd A
+
\left(\frac12-\frac1R\right)B^\perp(K_\sigma^0).
\end{equation}

Since \(K_\sigma^0\) is a compact metric space, continuous
functions are dense in \(L^1(K_\sigma^0,A+B)\).
Clipping continuous approximants of \(u_R\) to
\([-R,1/2-1/R]\) preserves continuity and does not increase the
\(L^1\)-distance to \(u_R\), because \(u_R\) takes values in this
interval. We thus obtain
\(\vartheta_{R,j}\in C(K_\sigma^0)\) such that
\[
-R\leq\vartheta_{R,j}\leq\frac12-\frac1R,
\qquad
\lim_{j\to\infty}\int_{K_\sigma^0}
|\vartheta_{R,j}-u_R|\,\dd(A+B)
=0.
\]
Since \(q\mapsto\frac12\log(1-2q)\) is Lipschitz on
\([-R,1/2-1/R]\), \eqref{eq:Ksigu} implies
\[
\lim_{j\to\infty}
\left[
\int_{K_\sigma^0}\vartheta_{R,j}\,\dd B
+\frac12\int_{K_\sigma^0}
\log(1-2\vartheta_{R,j})\,\dd A
\right]=
\int_{K_\sigma^0}\ell_R(v)\,\dd A
+
\left(\frac12-\frac1R\right)B^\perp(K_\sigma^0).
\]
Since \(0\leq\ell_R\uparrow\ell\), the right-hand side converges to
\(I_{\rm sv}(A,B)\) as \(R\to\infty\) by monotone convergence. Every
\(\vartheta_{R,j}\) is admissible in \eqref{eq:Isv-dual}, so
the right-hand side of \eqref{eq:Isv-dual} is at least
\(I_{\rm sv}(A,B)\).
\end{proof}
\begin{proposition}
\label{prop:sv}
The functional \(I_{\rm sv}\) is lower semicontinuous on
\(\mathcal Q_\sigma\), every \((A,B)\in\mathcal Q_\sigma\) satisfies
\begin{equation}\label{eq:IsvAB}
B(K_\sigma^0)\leq4I_{\rm sv}(A,B)+2\log2,
\end{equation}
and the sublevel sets \(\{I_{\rm sv}\leq L\}\), \(L<\infty\), are
compact.
\end{proposition}

\begin{proof}
Let \(\vartheta\in C(K_\sigma^0)\) satisfy \(\sup\vartheta<1/2\).
Then \(\vartheta\) and \(\log(1-2\vartheta)\) are bounded and continuous on
\(K_\sigma^0\), so the map
\(
(A,B)
\mapsto
\int\vartheta\,\dd B
+
\frac12\int\log(1-2\vartheta)\,\dd A
\)
is continuous under weak convergence of \(A\) and \(B\). By
Lemma~\ref{lem:Isv-dual}, \(I_{\rm sv}\) is the
supremum of these continuous functionals and is therefore lower
semicontinuous.

Since \(A\) is a probability measure, the constant test function
\(\vartheta\equiv1/4\) in \eqref{eq:Isv-dual} yields
\(I_{\rm sv}(A,B)\geq\frac14B(K_\sigma^0)-\frac12\log2\), which is
\eqref{eq:IsvAB}.
Now fix \(L<\infty\). By \eqref{eq:IsvAB}, every
\((A,B)\in\{I_{\rm sv}\leq L\}\) satisfies
\(B(K_\sigma^0)\leq4L+2\log2\). Since \(K_\sigma^0\) is compact,
Prokhorov's theorem implies that the
\(B\)-coordinates of \(\{I_{\rm sv}\leq L\}\) are relatively compact in
\(\Mplus(K_\sigma^0)\). The \(A\)-coordinates
lie in the set
\(
\{A\in\mathcal P(K_\sigma^0):(\pi_1)_\#A=\lambda\},
\)
which is compact as a closed subset of the compact space
\(\mathcal P(K_\sigma^0)\). Thus \(\{I_{\rm sv}\leq L\}\) is relatively
compact. Since \(I_{\rm sv}\) is lower semicontinuous, \(\{I_{\rm sv}\leq L\}\)
is closed in \(\mathcal Q_\sigma\), which is closed in
\(\mathcal P(K_\sigma^0)\times\Mplus(K_\sigma^0)\). Hence
\(\{I_{\rm sv}\leq L\}\) is compact.
\end{proof}
\subsection{Recovery of finite-action controls}
\label{subsec:recovery}

The lower bound is proved by changing the law of the driving noise on
the canonical Wiener space \(\mathcal W:=C_0([0,1];\R)\), with coordinate
process \(W_t(\omega):=\omega(t)\), its natural filtration
\((\mathcal F_t)\) and Wiener measure \(\mathbb P\). For \(0\leq k\leq n\),
we write \(\Delta_{<k}W:=(\Delta_0W,\ldots,\Delta_{k-1}W)\), with
\(\Delta_{<0}W:=\varnothing\).

A pair \((A,B)\in\mathcal Q_\sigma\) is called regular if \(B=vA\) with
\(r\leq v\leq R\) \(A\)-a.e.\ for some constants \(0<r<R<\infty\).
Lemma~\ref{lem:regularize} approximates a pair \((A,B)\) with
\(I_{\rm sv}(A,B)<\infty\) by regular pairs, and Lemma~\ref{lem:marked-approximation} approximates the
resulting controls by piecewise-constant ones with the prescribed marks
and terminal label.
Lemma~\ref{lem:restoring-tracking} shows that a diffusion with a strong
restoring drift reproduces the resulting state labels, and
Lemma~\ref{lem:recovery-tilted-law} realizes this diffusion by the Euler
scheme under conditional Gaussian changes of measure with explicit
entropy. The comparison of the tilted scheme with the diffusion is proved
in Appendix~\ref{app:euler-comparison}, and Step~1 of the proof of
Theorem~\ref{thm:cutoffldp} combines these steps.

\begin{lemma}[Approximation by regular pairs]
\label{lem:regularize}
For every \((A,B)\in\mathcal Q_\sigma\) with \(I_{\rm sv}(A,B)<\infty\), there
exist regular pairs \((A_m,B_m)\in\mathcal Q_\sigma\), \(m\geq2\), such that
\begin{equation}
(A_m,B_m)\Rightarrow(A,B),\qquad
I_{\rm sv}(A_m,B_m)\to I_{\rm sv}(A,B).
\label{eq:tag-5-26}
\end{equation}
\end{lemma}
\begin{proof}
Write \(B=vA+B^\perp\) with \(B^\perp\perp A\). We truncate
the density \(v\), spread the singular part \(B^\perp\) over short time
windows, and add a small multiple of the spread measure to \(A\), so that
the new density is bounded above and below.

We first truncate \(v\). Since \(B\) is finite, \(v\in L^1(A)\).
Set \(v_m^0:=(v\vee m^{-1})\wedge m\) for \(m\geq2\). Since
\(v_m^0\to v\) pointwise and \(|v_m^0-v|\leq v+1\), the dominated
convergence theorem gives \(\int_{K_\sigma^0}|v_m^0-v|\,\dd A\to0\), and in
particular \(v_m^0A\Rightarrow vA\).
Moreover, \(v_m^0=v\) on \(\{m^{-1}\leq v\leq m\}\), and
outside this set \(v_m^0\) is the endpoint of \([m^{-1},m]\) nearest to
\(v\). Since \(\ell\) is decreasing on \((0,1]\), increasing on
\([1,\infty)\), and \(m^{-1}\leq1\leq m\), this gives
\(0\leq\ell(v_m^0)\leq\ell(v)\). Since \(I_{\rm sv}(A,B)<\infty\), we have
\(\ell(v)\in L^1(A)\), and the dominated convergence theorem yields
\begin{equation}
	\int_{K_\sigma^0}\ell(v_m^0)\,\dd A
	\rightarrow
	\int_{K_\sigma^0}\ell(v)\,\dd A .
	\label{eq:tag-5-27}
\end{equation}

We next spread the singular part in time. Choose a
symmetric nonnegative function \(\rho\in C_c(\R)\) with
\(\operatorname{supp}\rho\subseteq[-1/2,1/2]\) and
\(\int_{\R}\rho(u)\,\dd u=1\). For \(m\geq2\), set
\[
\rho_m(s,t)
:=
\frac{\rho(m(s-t))}
{\int_0^1\rho(m(u-t))\,\dd u},
\qquad s,t\in[0,1].
\]
Since \(\rho\) is symmetric and \(m\geq2\), the denominator is at
least \(1/(2m)\). Hence \(\rho_m\) is continuous and bounded on
\([0,1]^2\), and \(\rho_m(\cdot,t)\) is a probability density on
\([0,1]\), supported on \(\{s\in[0,1]:|s-t|\leq 1/(2m)\}\). Define
\[
S_m(\dd s,\dd\xi)
:=
\int_{[0,1]}
\rho_m(s,t)\,B^\perp(\dd t,\dd\xi)\,\dd s.
\]
Then \(S_m(K_\sigma^0)=B^\perp(K_\sigma^0)\), and the time marginal
of \(S_m\) is \(b_m(s)\,\dd s\), where
\[
b_m(s):=
\int_{[0,1]\times K_\sigma}
\rho_m(s,t)\,B^\perp(\dd t,\dd\xi)
\]
is bounded and continuous. Since \(\rho_m(\cdot,t)\) is supported in
\(\{|s-t|\leq 1/(2m)\}\), for every \(f\in C(K_\sigma^0)\),
\[
\left|
\int f\,\dd S_m-\int f\,\dd B^\perp
\right|
\leq
B^\perp(K_\sigma^0)
\sup_{\substack{|s-t|\leq 1/(2m),\ \xi\in K_\sigma}}
|f(s,\xi)-f(t,\xi)|
{\rightarrow}0\quad \text{as }m\to\infty,
\]
so \(S_m\Rightarrow B^\perp\) as \(m\to\infty\).

We now choose the weight of \(S_m\) in \(A_m\). For each
\(m\geq2\), the map \((c,w)\mapsto(1-c)\ell(\frac{w}{1-c})\)
is uniformly continuous on \([0,1/2]\times[m^{-1},m]\).
Hence there exists \(\delta_m\in(0,m^{-1}]\) such that
\begin{equation}\label{eq:svbound}
\sup_{\substack{0\leq c\leq\delta_m,\ m^{-1}\leq w\leq m}}
\left|
(1-c)\ell\left(\frac{w}{1-c}\right)-\ell(w)
\right|
\leq\frac1m.
\end{equation}
Set \(R_m:=\max\{m,\|b_m\|_\infty/\delta_m\}\) and
\(c_m(t):=b_m(t)/R_m\), so that \(R_m\to\infty\) and
\(0\leq c_m\leq\delta_m\leq m^{-1}\). Define
\begin{equation}
	A_m
	:=
	(1-c_m(t))A+\frac1{R_m}S_m,
	\qquad
	B_m
	:=
	v_m^0A+S_m.
	\label{eq:tag-5-28}
\end{equation}
The time marginal of \(A_m\) is
\(
(1-c_m(t))\,\dd t
+\frac{b_m(t)}{R_m}\,\dd t
=
\dd t,
\)
so \((A_m,B_m)\in\mathcal Q_\sigma\). Since \(c_m\to0\) uniformly,
\(S_m(K_\sigma^0)/R_m\to0\), \(S_m\Rightarrow B^\perp\), and
\(v_m^0A\Rightarrow vA\), we obtain \((A_m,B_m)\Rightarrow(A,B)\).

We now verify that \(B_m\) has a density with respect to \(A_m\) that is
bounded above and away from zero. Let \(\lambda_m:=A+S_m\),
and let \(a_m:=\dd A/\dd\lambda_m\in[0,1]\), so that
\(\dd S_m/\dd\lambda_m=1-a_m\). By \eqref{eq:tag-5-28},
\begin{equation}\label{eq:amdef}
\frac{\dd A_m}{\dd\lambda_m}
=
(1-c_m)a_m+\frac{1-a_m}{R_m},
\qquad
\frac{\dd B_m}{\dd\lambda_m}
=
v_m^0a_m+1-a_m.
\end{equation}
The first density is at least \(\min\{1-c_m,1/R_m\}>0\), since it is a
convex combination of \(1-c_m\) and \(1/R_m\) with weights \(a_m\) and
\(1-a_m\). Hence \(B_m\ll A_m\), and
\begin{equation}\label{eq:vmdef}
v_m:=\frac{\dd B_m}{\dd A_m}
=
\frac{(1-c_m)a_m\cdot\frac{v_m^0}{1-c_m}+\frac{1-a_m}{R_m}\cdot R_m}
{(1-c_m)a_m+(1-a_m)/R_m}
\end{equation}
is a weighted average of \(v_m^0/(1-c_m)\) and \(R_m\).
Since \(m^{-1}\leq v_m^0/(1-c_m)\leq m/(1-m^{-1})\) and \(R_m\geq m\),
we obtain
\[
m^{-1}
\leq
v_m
\leq
\max\left\{
\frac{m}{1-m^{-1}},R_m
\right\},
\qquad A_m\text{-a.e.}
\]
Hence \((A_m,B_m)\) is regular.

It remains to show that \(I_{\rm sv}(A_m,B_m)\to I_{\rm sv}(A,B)\). Let
\(g(x,y):=x\ell(y/x)\) for \(x>0\) and \(y\geq0\), with \(g(0,0):=0\). Since
\(\ell\) is convex, \(g\) is convex and positively one-homogeneous, and
hence subadditive. Since \(B_m=v_mA_m\),
\[
I_{\rm sv}(A_m,B_m)
=\int_{K_\sigma^0}\ell(v_m)\,\dd A_m
=\int_{K_\sigma^0}
g\Bigl(\frac{\dd A_m}{\dd\lambda_m},\frac{\dd B_m}{\dd\lambda_m}\Bigr)
\dd\lambda_m.
\]
By \eqref{eq:amdef}, the pair of densities is the sum of
\(((1-c_m)a_m,v_m^0a_m)\) and \(((1-a_m)/R_m,1-a_m)\). Subadditivity and
homogeneity of \(g\) therefore give
\begin{equation}\label{eq:AmBm}
\begin{aligned}
	I_{\rm sv}(A_m,B_m)
	\leq{}&
	\int_{K_\sigma^0}
	(1-c_m)a_m
	\ell\!\left(\frac{v_m^0}{1-c_m}\right)\dd\lambda_m+
	\int_{K_\sigma^0}
	\frac{1-a_m}{R_m}\ell(R_m)\,\dd\lambda_m\\
	={}&\int_{K_\sigma^0}(1-c_m)\ell\!\left(\frac{v_m^0}{1-c_m}\right)\dd A
	+\frac{\ell(R_m)}{R_m}
	B^\perp(K_\sigma^0).
\end{aligned}
\end{equation}
By \eqref{eq:svbound}, \(m^{-1}\leq v_m^0\leq m\), \(c_m\leq\delta_m\)
and \(A(K_\sigma^0)=1\), the first integral differs from
\(\int\ell(v_m^0)\,\dd A\) by at most \(1/m\). Since
\(\ell(R_m)/R_m\to\frac12\) as \(R_m\to\infty\),
\eqref{eq:tag-5-27} and \eqref{eq:AmBm} yield
\[
\limsup_{m\to\infty}I_{\rm sv}(A_m,B_m)
\leq
\int\ell(v)\,\dd A
+
\frac12B^\perp(K_\sigma^0)
=
I_{\rm sv}(A,B).
\]
Conversely, the lower semicontinuity of \(I_{\rm sv}\)
(Proposition~\ref{prop:sv}) and
\((A_m,B_m)\Rightarrow(A,B)\) imply
\(\liminf_{m\to\infty}I_{\rm sv}(A_m,B_m)\geq I_{\rm sv}(A,B)\).
\end{proof}
By Lemma~\ref{lem:regularize}, it suffices to recover controls whose
pair \((A,B)\) is regular. The next lemma approximates such a control by
the control \(\zeta^\eta\) in \eqref{eq:zetaeta}, which is generated by a
piecewise-constant state path \(q^\eta\), a piecewise-constant variance
density \(u^\eta\) and marks at strictly ordered times, and
Lemma~\ref{lem:restoring-tracking} realizes such controls by a diffusion.
Throughout this section, the number \(N\) of marks may be zero, in which
case lists indexed by \(j=1,\ldots,N\) are empty, sums over \(j\) vanish
and conditions on \(j\) are vacuous.

\begin{lemma}[Piecewise-constant approximation of regular marked controls]
\label{lem:marked-approximation}
Let \(\zeta=(A,B,\mathfrak m,\xi_f)\in\mathcal C_\varepsilon\) with
\(\mathfrak m=((t_j,\xi_j,\delta_j))_{j=1}^N\), and suppose that \((A,B)\)
is regular with \(B=vA\) and \(r\leq v\leq R\) \(A\)-a.e. Then, for every
open neighbourhood \(O\) of \(\zeta\) and every \(\eta>0\), there exist a
partition \(0=\tau_0<\cdots<\tau_P=1\) with \(P\geq2\), times
\(0<s_1<\cdots<s_N<1\) among the \(\tau_i\), and functions
\(q^\eta:[0,1]\to\R\) and \(u^\eta:[0,1]\to[r,R]\), constant on
\([\tau_i,\tau_{i+1})\) for \(0\leq i<P-1\) and on \([\tau_{P-1},1]\),
such that the control
\begin{equation}\label{eq:zetaeta}
\zeta^\eta
:=
\bigl(
\dd t\,\delta_{\iota(q^\eta(t))},
\;
u^\eta(t)\,\dd t\,\delta_{\iota(q^\eta(t))},
\;
((s_j,\iota(q^\eta(s_j-)),\delta_j))_{j=1}^N,
\;
\iota(q^\eta(1))
\bigr)
\end{equation}
belongs to \(O\) and satisfies
\(\int_0^1\ell(u^\eta(t))\,\dd t\leq I_{\rm sv}(A,B)+\eta\).
\end{lemma}
\begin{proof}
Fix an open neighbourhood \(O\) of \(\zeta\) and \(\eta>0\).

\medskip
\noindent\textit{Step 1. We reduce \(O\) to a product neighbourhood and fix the target values of the marks and of the terminal tag.}
Since \(\mathcal C_{\varepsilon,N}\) has the topology
inherited from the product space, there exist an open neighbourhood
\(U_{AB}\) of \((A,B)\) in \(\mathcal Q_\sigma\), open neighbourhoods
\(U_j\) of \(\xi_j\) and \(U_f\) of \(\xi_f\) in \(K_\sigma\), and open
intervals \(T_j\ni t_j\) and \(D_j\ni\delta_j\) such that
\begin{equation}
\Bigl(
U_{AB}
\times
\prod_{j=1}^N(T_j\times U_j\times D_j)
\times U_f
\Bigr)
\cap\mathcal C_{\varepsilon,N}
\subset O.
\label{eq:marked-approx-neighbourhood}
\end{equation}
By the density of \(\iota(\R)\) in \(K_\sigma\), choose
\(q_1,\ldots,q_N,q_f\in\R\) with \(\iota(q_j)\in U_j\) and
\(\iota(q_f)\in U_f\). For sufficiently small \(\theta\in(0,1)\), the
times \(s_j:=(1-\theta)t_j+\theta j/(N+1)\) satisfy
\(s_j\in T_j\cap(0,1)\) and \(s_1<\cdots<s_N\), even if some \(t_j\)
coincide or belong to \(\{0,1\}\). Therefore it suffices to
construct piecewise-constant \(q^\eta,u^\eta\), with \(r\leq u^\eta\leq R\),
such that the pair generated by them satisfies

\begin{equation}
\bigl(\dd t\,\delta_{\iota(q^\eta(t))},\,u^\eta(t)\,\dd t\,\delta_{\iota(q^\eta(t))}\bigr)\in U_{AB},
\qquad
q^\eta(s_j-)=q_j,
\qquad
q^\eta(1)=q_f,
\label{eq:marked-approx-targets}
\end{equation}
and \(\int_0^1\ell(u^\eta)\,\dd t\leq I_{\rm sv}(A,B)+\eta\). Indeed, the
control \(\zeta^\eta\) in \eqref{eq:zetaeta} then lies in
\(\mathcal C_\varepsilon\), since the times \(s_j\) are strictly ordered
and the increments \(\delta_j\) are unchanged, and its coordinates lie in
the sets in \eqref{eq:marked-approx-neighbourhood}, so that
\(\zeta^\eta\in O\). Steps~2 and~3 construct \(q^\eta\) and \(u^\eta\).

\medskip
\noindent\textit{Step 2. We approximate \((A,B)\) by a pair generated by piecewise-constant functions.}
We claim that there exist functions \(q^0:[0,1]\to\R\) and
\(u^0:[0,1]\to[r,R]\), piecewise constant on a finite partition and
right-continuous on \([0,1)\), such that
\begin{equation}
(A^0,B^0):=
\bigl(\dd t\,\delta_{\iota(q^0(t))},\,u^0(t)\,\dd t\,\delta_{\iota(q^0(t))}\bigr)\in U_{AB},
\qquad
\int_0^1\ell(u^0(t))\,\dd t
\leq
I_{\rm sv}(A,B)+\eta.
\label{eq:marked-approx-step2}
\end{equation}

This is an instance of the approximation of relaxed controls by
ordinary controls \cite{Warga}. We give the short construction, since it
has to be carried out in the label space \(K_\sigma\). Choosing a version of \(v\) with values in \([r,R]\), we lift \((A,B)\) to the probability measure
\(\Pi(\dd t,\dd\xi,\dd u):=A(\dd t,\dd\xi)\,\delta_{v(t,\xi)}(\dd u)\) on
\([0,1]\times K_\sigma\times[r,R]\). It has Lebesgue time marginal, and
\[
\mathcal T(\Pi)=(A,B),
\qquad
\int\ell(u)\,\Pi(\dd t,\dd\xi,\dd u)=I_{\rm sv}(A,B),
\qquad
\mathcal T(\Pi'):=\bigl((\pi_{t,\xi})_\#\Pi',\,(\pi_{t,\xi})_\#(u\Pi')\bigr).
\]
Since \(\pi_{t,\xi}\) is continuous and \(u\) is bounded and continuous on
\([0,1]\times K_\sigma\times[r,R]\), the map \(\mathcal T\) is continuous
for the weak topologies. Therefore it suffices to find piecewise-constant
\(q^0,u^0\) such that \(\dd t\,\delta_{(\iota(q^0(t)),u^0(t))}\) is close
to \(\Pi\) weakly and \(\int\ell(u^0)\,\dd t\) is close to
\(\int\ell(u)\,\dd\Pi\). Since \(\ell\) is continuous on \([r,R]\), both
follow from
\begin{equation}
\dd t\,\delta_{(\iota(q^{(k)}(t)),u^{(k)}(t))}\Rightarrow\Pi
\qquad\text{as }k\to\infty
\label{eq:marked-approx-weak}
\end{equation}
for suitable piecewise-constant \(q^{(k)},u^{(k)}\), which we now construct.

Fix \(k\geq1\). Divide \([0,1]\) into the intervals
\(\mathcal J_a:=[a/k,(a+1)/k)\), \(a=0,\ldots,k-2\), and \(\mathcal J_{k-1}:=[(k-1)/k,1]\).
Partition the compact space \(K_\sigma\times[r,R]\) into nonempty Borel
sets \(\mathcal V_1,\ldots,\mathcal V_{m_k}\) of diameter at most \(1/k\) for
the metric \(d_{K_\sigma}(\xi,\xi')+|u-u'|\). By the density of
\(\iota(\R)\), each \(\mathcal V_p\) contains a point \((\xi_p',\widehat u_p)\) such that
\(d_{K_\sigma}(\iota(\widehat q_p),\xi_p')\leq1/k\) for some
\(\widehat q_p\in\R\). Within each \(\mathcal J_a\), we
let \((q^{(k)},u^{(k)})\) take the value \((\widehat q_p,\widehat u_p)\) on
consecutive intervals of length \(\Pi(\mathcal J_a\times\mathcal V_p)\),
\(p=1,\ldots,m_k\). These lengths add up to \(|\mathcal J_a|=1/k\), because \(\Pi\) has
Lebesgue time marginal. Hence \(q^{(k)},u^{(k)}\) are well defined,
right-continuous, and piecewise constant, and we set their values at \(1\)
equal to their left limits. On the interval assigned to
\(\mathcal J_a\times\mathcal V_p\), the point
\((t,\iota(q^{(k)}(t)),u^{(k)}(t))\) lies within distance \(3/k\) of every
point of \(\mathcal J_a\times\mathcal V_p\). Thus, for every
\(f\in C([0,1]\times K_\sigma\times[r,R])\) with modulus of continuity
\(\omega_f\),
\[
\Bigl|
\int_0^1
f\bigl(t,\iota(q^{(k)}(t)),u^{(k)}(t)\bigr)\,\dd t
-
\int f\,\dd\Pi
\Bigr|
\leq
\sum_{a=0}^{k-1}\sum_{p=1}^{m_k}
\Pi(\mathcal J_a\times\mathcal V_p)\,\omega_f(3/k)
=\omega_f(3/k)\to0.
\]
This proves \eqref{eq:marked-approx-weak}, and hence the claim, with
\((q^0,u^0):=(q^{(k)},u^{(k)})\) for \(k\) large.

\medskip
\noindent\textit{Step 3. We modify \(q^0\) on short intervals before the mark times and near \(1\), so that the marks and the terminal tag take the target values.}
Since \(U_{AB}\) is open, there exists \(d_0>0\) such that every
\((A',B')\in\mathcal Q_\sigma\) with
\(d_{\mathcal Q_\sigma}((A',B'),(A^0,B^0))<d_0\) lies in \(U_{AB}\). Since
\(0<s_1<\cdots<s_N<1\), we can choose \(\kappa>0\) so small that the
intervals
\[
I_j:=[s_j-\kappa,s_j)\subset(0,1),\quad j=1,\ldots,N,
\qquad
I_f:=[1-\kappa,1]
\]
are pairwise disjoint and their union \(I\) satisfies
\(2(1+R)|I|=2(1+R)(N+1)\kappa<d_0\). Set \(u^\eta:=u^0\), and let \(q^\eta\) equal \(q_j\) on
\(I_j\), \(q_f\) on \(I_f\), and \(q^0\) elsewhere. Then
\(q^\eta(s_j-)=q_j\) and \(q^\eta(1)=q_f\). Since \(q^\eta\) and \(q^0\)
differ only on \(I\), the pair \((A^\eta,B^\eta)\) generated by
\((q^\eta,u^\eta)\) and the pair \((A^0,B^0)\) differ only on times
in \(I\). Hence, for
every \(f\in C(K_\sigma^0)\) with \(\|f\|_\infty\leq1\),
\[
\Bigl|\int f\,\dd A^\eta-\int f\,\dd A^0\Bigr|\leq2|I|,
\qquad
\Bigl|\int f\,\dd B^\eta-\int f\,\dd B^0\Bigr|\leq2R|I|,
\]
and \eqref{eq:dw} gives
\(d_{\mathcal Q_\sigma}((A^\eta,B^\eta),(A^0,B^0))\leq2(1+R)|I|<d_0\).
Thus \((A^\eta,B^\eta)\in U_{AB}\), and \eqref{eq:marked-approx-targets}
holds. Moreover,
\(\int_0^1\ell(u^\eta)\,\dd t=\int_0^1\ell(u^0)\,\dd t\leq
I_{\rm sv}(A,B)+\eta\) by \eqref{eq:marked-approx-step2}.

Finally, let \(0=\tau_0<\cdots<\tau_P=1\) be a partition containing all
breakpoints of \(q^0,u^0\) and all endpoints of \(I_1,\ldots,I_N,I_f\).
It contains every \(s_j\), and \(P\geq2\) since \(1-\kappa\in(0,1)\).
Since \(q^0,u^0\) are right-continuous and the intervals \(I_j\) are
closed on the left and open on the right, \(q^\eta\) and \(u^\eta\) are
constant on each \([\tau_i,\tau_{i+1})\) and on \([\tau_{P-1},1]\), as
stated in the lemma.
\end{proof}

\begin{lemma}[Tracking of the piecewise-constant approximation]
\label{lem:restoring-tracking}
Let \(\zeta^\eta\) be the control \eqref{eq:zetaeta} constructed in
Lemma~\ref{lem:marked-approximation} from \(q^\eta\), \(u^\eta\) and
\((s_j,\delta_j)_{j=1}^N\), and let \(\widetilde W\) be a standard Brownian
motion under a probability measure \(\widetilde{\mathbb P}\).
For each \(K\geq1\), set \(\psi_K(x):=Kx/(1+|x|)\), and let
\(\Gamma^{\eta,K}\) be the c\`adl\`ag solution of
\begin{equation}
\Gamma_t^{\eta,K}
=
y_0
+
\int_0^t
\psi_K(q^\eta(r)-\Gamma_r^{\eta,K})\,\dd r
+
\int_0^t
\sqrt{u^\eta(r)}\,\dd\widetilde W_r
+
\sum_{\substack{1\leq j\leq N,\ s_j\leq t}}
\bigl[
\mathcal S_{\delta_j}(\Gamma_{s_j-}^{\eta,K})
-
\Gamma_{s_j-}^{\eta,K}
\bigr],
\qquad t\in[0,1],
\label{eq:restoring-sde}
\end{equation}
Then the controls
\begin{equation}
\zeta^{\eta,K}
:=
\left(
\dd t\,\delta_{\iota(\Gamma_t^{\eta,K})},
\;
u^\eta(t)\,\dd t\,\delta_{\iota(\Gamma_t^{\eta,K})},
\;
((s_j,\iota(\Gamma_{s_j-}^{\eta,K}),\delta_j))_{j=1}^N,
\;
\iota(\Gamma_1^{\eta,K})
\right)
\rightarrow\zeta^\eta
\label{eq:restoring-control-convergence}
\end{equation}
converge in \(\widetilde{\mathbb P}\)-probability in \(\mathcal C_\varepsilon\)
as \(K\to\infty\).
\end{lemma}

\begin{proof}
For fixed \(K\), the drift in \eqref{eq:restoring-sde} is
bounded and Lipschitz in the state, so \eqref{eq:restoring-sde} has a
unique strong solution, obtained by solving between the times \(s_j\) and
applying the updates at these times. The restoring drift \(\psi_K\) pulls
\(\Gamma^{\eta,K}\) towards \(q^\eta\) at rate of order \(K\), while the
noise has bounded variance. Hence the tracking error
\(\Gamma^{\eta,K}-q^\eta\) tends to zero in \(L^2\) at every non-partition
time and immediately before every partition point, which leads to
\eqref{eq:restoring-control-convergence}. Constants \(C\) may depend on
the fixed data, but not on \(K\geq1\), and \(\widetilde{\mathbb E}\)
denotes expectation under \(\widetilde{\mathbb P}\).

By construction, the partition contains every \(s_j\), and
\(q^\eta,u^\eta\) are constant on each \([\tau_i,\tau_{i+1})\). For
\(i=0,\ldots,P-1\) and \(t\in[\tau_i,\tau_{i+1})\), set
\(e_{i,t}:=\Gamma_t^{\eta,K}-q^\eta(\tau_i)\), and let
\(e_{i,\tau_{i+1}-}\) be its left limit at \(\tau_{i+1}\). Since
\(\psi_K(-x)=-\psi_K(x)\), \eqref{eq:restoring-sde} gives
\begin{equation}
\dd e_{i,t}
=-\frac{Ke_{i,t}}{1+|e_{i,t}|}\,\dd t
+\sqrt{u^\eta(\tau_i)}\,\dd\widetilde W_t,
\qquad t\in[\tau_i,\tau_{i+1}).
\label{eq:restoring-error-sde}
\end{equation}
At an interior partition point \(\tau_{i+1}\), the process
\(\Gamma^{\eta,K}\) jumps by \(\mathcal S_{\delta_j}(\Gamma_{s_j-}^{\eta,K})-\Gamma_{s_j-}^{\eta,K}\)
if \(\tau_{i+1}=s_j\), and is continuous otherwise. Since
\(|\mathcal S_\delta(z)-z|\leq(\sigma_{\mathrm{max}}/\sigma_{\mathrm{min}})|\delta|\)
by \eqref{eq:beta3}, the triangle inequality gives
\begin{equation}
|e_{i+1,\tau_{i+1}}|
\leq
|e_{i,\tau_{i+1}-}|
+\bigl|q^\eta(\tau_{i+1})-q^\eta(\tau_{i})\bigr|
+\frac{\sigma_{\mathrm{max}}}{\sigma_{\mathrm{min}}}
\sum_{j:\,s_j=\tau_{i+1}}|\delta_j|,
\qquad i=0,\ldots,P-2.
\label{eq:restoring-partition-error}
\end{equation}

The drift in \eqref{eq:restoring-error-sde} is linear in \(e_{i,t}\) only
near the origin. To control the region where \(|e_{i,t}|\) is large, we
first show that
\begin{equation}
\sup_{K\geq1}\max_{0\leq i<P}
\sup_{\tau_i\leq t<\tau_{i+1}}
\widetilde{\mathbb E}|e_{i,t}|^4
<\infty.
\label{eq:restoring-uniform-fourth}
\end{equation}
By It\^o's formula and \eqref{eq:restoring-error-sde},
\[
\dd|e_{i,t}|^4
=-\frac{4K|e_{i,t}|^4}{1+|e_{i,t}|}\,\dd t
+6u^\eta(\tau_i)|e_{i,t}|^2\,\dd t
+4\sqrt{u^\eta(\tau_i)}\,e_{i,t}^3\,\dd\widetilde W_t.
\]
There exists \(C_0=C_0(R)\) such that
\(6Rx^2\leq4Kx^4/(1+|x|)+C_0\) for all \(K\geq1\) and \(x\in\R\). Since
\(u^\eta\leq R\), it follows that the \(\dd t\)-terms are at most
\(C_0\,\dd t\). Taking expectations, after localization, and
applying Fatou's lemma give
\[
\sup_{\tau_i\leq t<\tau_{i+1}}
\widetilde{\mathbb E}|e_{i,t}|^4
+\widetilde{\mathbb E}|e_{i,\tau_{i+1}-}|^4
\leq
2\widetilde{\mathbb E}|e_{i,\tau_i}|^4
+2C_0.
\]
By \eqref{eq:restoring-partition-error},
\(\widetilde{\mathbb E}|e_{i+1,\tau_{i+1}}|^4\leq
C\widetilde{\mathbb E}|e_{i,\tau_{i+1}-}|^4+C\), since the partition values
and the increments \(\delta_j\) are fixed. Starting from the deterministic
value \(e_{0,0}=y_0-q^\eta(0)\) and iterating these two bounds over the
\(P\) intervals proves \eqref{eq:restoring-uniform-fourth}.

Using \eqref{eq:restoring-uniform-fourth}, we next show that the tracking
error vanishes in \(L^2\) away from the left endpoints of the partition
intervals, namely, for \(i=0,\ldots,P-1\),
\begin{equation}
\lim_{K\to\infty}\widetilde{\mathbb E}|e_{i,t}|^2=0,
\quad \tau_i<t<\tau_{i+1},
\qquad
\lim_{K\to\infty}\widetilde{\mathbb E}|e_{i,\tau_{i+1}-}|^2=0.
\label{eq:restoring-interior-L2}
\end{equation}
By It\^o's formula,
\[
\frac{\dd}{\dd t}\widetilde{\mathbb E}|e_{i,t}|^2
=-2K\,\widetilde{\mathbb E}\frac{|e_{i,t}|^2}{1+|e_{i,t}|}
+u^\eta(\tau_i)
\qquad\text{for a.e. }t\in(\tau_i,\tau_{i+1}).
\]
Let \(L\geq1\). Since
\(x^2/(1+|x|)\geq x^2\mathbf1_{\{|x|\leq L\}}/(1+L)\), and since
\(\widetilde{\mathbb E}[|e_{i,t}|^2\mathbf1_{\{|e_{i,t}|>L\}}]\leq C/L^2\)
by \eqref{eq:restoring-uniform-fourth}, the dissipation term is at most
\(-\frac{2K}{1+L}(\widetilde{\mathbb E}|e_{i,t}|^2-C/L^2)\). Together with
\(u^\eta\leq R\) and \(\widetilde{\mathbb E}|e_{i,\tau_i}|^2\leq C\),
Gronwall's inequality gives, for \(\tau_i\leq t<\tau_{i+1}\),
\begin{equation}
\widetilde{\mathbb E}|e_{i,t}|^2
\leq
C\exp\left(-\frac{2K}{1+L}(t-\tau_i)\right)
+\frac{C}{L^2}+\frac{R(1+L)}{2K}.
\label{eq:restoring-second-moment-bound}
\end{equation}
By Fatou's lemma, the same bound holds for
\(\widetilde{\mathbb E}|e_{i,\tau_{i+1}-}|^2\) with \(t=\tau_{i+1}\), that
is,
\begin{equation}
\widetilde{\mathbb E}|e_{i,\tau_{i+1}-}|^2
\leq
C\exp\left(-\frac{2K}{1+L}(\tau_{i+1}-\tau_i)\right)
+\frac{C}{L^2}+\frac{R(1+L)}{2K}.
\label{eq:restoring-prejump-endpoint-bound}
\end{equation}
Choosing \(L=K^{1/4}\) in
\eqref{eq:restoring-second-moment-bound} and
\eqref{eq:restoring-prejump-endpoint-bound}, all three terms tend to zero
as \(K\to\infty\), for every \(t>\tau_i\). This proves \eqref{eq:restoring-interior-L2}.

It remains to deduce \eqref{eq:restoring-control-convergence} from
\eqref{eq:restoring-interior-L2}. The mark times and increments of
\(\zeta^{\eta,K}\) and \(\zeta^\eta\) coincide. If \(s_j=\tau_\ell\), then
\(\Gamma_{s_j-}^{\eta,K}-q^\eta(s_j-)=e_{\ell-1,\tau_\ell-}\). Since no
update occurs at \(1\), also
\(\Gamma_1^{\eta,K}-q^\eta(1)=e_{P-1,1-}\). By
\eqref{eq:restoring-interior-L2} and the continuity of \(\iota\), the mark
labels and the terminal label of \(\zeta^{\eta,K}\) converge in
probability to those of \(\zeta^\eta\).

For the measure coordinates, let \(f\in C(K_\sigma^0)\).
By \eqref{eq:restoring-interior-L2} and the continuity of \(\iota\) and
\(f\), for every \(t\) outside the partition,
\(f(t,\iota(\Gamma_t^{\eta,K}))-f(t,\iota(q^\eta(t)))\to0\) in
probability. Since this difference is bounded by \(2\|f\|_\infty\) and
\(u^\eta\leq R\), dominated convergence gives
\[
\widetilde{\mathbb E}\int_0^1
\bigl(1+u^\eta(t)\bigr)
\bigl|f(t,\iota(\Gamma_t^{\eta,K}))-f(t,\iota(q^\eta(t)))\bigr|\,\dd t
\longrightarrow0.
\]
Applying this to the functions \(\varphi_j\) in
\eqref{eq:dw} and controlling the tail of the series, we conclude that the
two measure coordinates of \(\zeta^{\eta,K}\) converge in probability to
those of \(\zeta^\eta\) in \(\mathcal Q_\sigma\). Together with the
convergence of the labels, this proves
\eqref{eq:restoring-control-convergence}.
\end{proof}

By Lemma~\ref{lem:restoring-tracking}, it remains to realize the law of
\(\zeta^{\eta,K}\) as the limit of the laws of the empirical controls
\(\zeta_n^\varepsilon\) under changes of measure. Under the measures
\(Q_n^{\eta,K}\) of the next lemma, the increment of step
\(\lceil ns_j\rceil\) is close to \(\delta_j\) and produces the \(j\)-th
mark. Every other increment \(\Delta_kW\) is conditionally Gaussian with
variance \(u^\eta(k/n)/n\), and its conditional mean compensates the drift
\(\gamma\) and the quadratic correction \(-\frac12\alpha u^\eta\) of the
transformed scheme, so that the scheme follows the restoring drift in
\eqref{eq:restoring-sde}. Recall from \eqref{eq:empirical-control} that
each increment contributes either to \(B_n^\varepsilon\) or to the marked
list, but never to both, so that the two kinds of increments can be
treated separately.

\begin{lemma}[{Tilted approximation and entropy}]
\label{lem:recovery-tilted-law}
Let \(K\geq1\), and let \(\zeta^{\eta,K}\) be the control in
\eqref{eq:restoring-control-convergence}, constructed from the data
\(q^\eta\), \(u^\eta\) and \((s_j,\delta_j)_{j=1}^N\) of
Lemma~\ref{lem:marked-approximation}. Then there exist probability
measures \(Q_n^{\eta,K}\) on \(\mathcal W\) with
\(H(Q_n^{\eta,K}\mid\mathbb P)<\infty\), \(n\geq1\), such that
\begin{equation}
\mathcal L_{Q_n^{\eta,K}}
\bigl(\zeta_n^\varepsilon\bigr)
\Longrightarrow
\mathcal L_{\widetilde{\mathbb P}}
\bigl(\zeta^{\eta,K}\bigr)
\quad\text{on }\mathcal C_\varepsilon,
\qquad n\to\infty,
\label{eq:recovery-joint-convergence}
\end{equation}
and
\begin{equation}
\limsup_{n\to\infty}
\frac1nH(Q_n^{\eta,K}\mid\mathbb P)
\leq
\int_0^1\ell(u^\eta(t))\,\dd t
+\frac12\sum_{j=1}^N\delta_j^2.
\label{eq:recovery-entropy-bound}
\end{equation}
\end{lemma}
\begin{proof}
Throughout the proof, \(\eta\), \(K\), the functions
\(q^\eta,u^\eta\), their common partition
\(0=\tau_0<\cdots<\tau_P=1\), the times \(s_1<\cdots<s_N\) and the
increments \(\delta_1,\ldots,\delta_N\) are fixed, and the partition
contains every \(s_j\). Constants may depend on these data and on the
coefficients, but not on \(n\), and \(\widetilde{\mathbb E}\) denotes
expectation under \(\widetilde{\mathbb P}\).

\medskip\noindent\textit{Step 1. We define \(Q_n^{\eta,K}\) by an explicit density, identify the conditional laws of the increments, and prove \eqref{eq:recovery-entropy-bound}.}
Since \(0<s_1<\cdots<s_N<1\), there exists \(n_0\) such
that, for every \(n\geq n_0\), the integers \(\lceil ns_j\rceil\),
\(j=1,\ldots,N\), are distinct and belong to \(\{1,\ldots,n-1\}\). For
\(n<n_0\), set \(Q_n^{\eta,K}=\mathbb P\). From now on, \(n\geq n_0\),
and
\begin{equation}
\mathcal K_n^1:=\{\lceil ns_j\rceil:j=1,\ldots,N\},\qquad
\mathcal K_n^0:=\{0,\ldots,n-1\}\setminus\mathcal K_n^1.
\label{eq:recovery-index-sets}
\end{equation}
For \(k\in\mathcal K_n^0\), define
\begin{equation}
m_{k,n}^{\eta,K}
:=
\psi_K(q^\eta(k/n)-Z_k^n)
-\gamma(Z_k^n)
+\frac12\alpha(Z_k^n)u^\eta(k/n).
\label{eq:recovery-mean}
\end{equation}
By \eqref{eq:transformed-euler-update}, \(Z_k^n\) is a
continuous function of \(\Delta_{<k}W\). Since \(\psi_K,\gamma,\alpha\)
are continuous, \(m_{k,n}^{\eta,K}\) is
\(\mathcal F_{k/n}\)-measurable. Moreover, \(|\psi_K|\leq K\) and
\(u^\eta\leq R\) give
\begin{equation}
\max_{k\in\mathcal K_n^0}|m_{k,n}^{\eta,K}|
\leq
B_K:=K+\|\gamma\|_\infty+\frac R2\|\alpha\|_\infty.
\label{eq:recovery-mean-bound}
\end{equation}
For \(k=0,\ldots,n-1\), let \(\mu_{k,n}\) be the Gaussian
law
\begin{equation}
\mu_{k,n}
:=
\begin{cases}
N\bigg(
\dfrac{m_{k,n}^{\eta,K}}{n},
\dfrac{u^\eta(k/n)}{n}
\bigg),
&k\in\mathcal K_n^0,\\[2mm]
N\left(
\delta_j+\operatorname{sgn}(\delta_j)n^{-1/4},
n^{-2}
\right),
&k=\lceil ns_j\rceil,\quad j=1,\ldots,N,
\end{cases}
\label{eq:recovery-kernels}
\end{equation}
which depends on \(\omega\) only through
\(\mathcal F_{k/n}\)-measurable quantities. Let \(p_n:=N(0,n^{-1})\), and
define
\begin{equation}
\mathcal D_{k,n}:=\frac{\dd\mu_{k,n}}{\dd p_n}(\Delta_kW),
\qquad
D_k:=\prod_{\ell<k}\mathcal D_{\ell,n},
\qquad
Q_n^{\eta,K}(\dd\omega):=D_n(\omega)\,\mathbb P(\dd\omega).
\label{eq:recovery-density}
\end{equation}
We claim that \(Q_n^{\eta,K}\) is a probability measure
and that
\begin{equation}
\mathcal L_{Q_n^{\eta,K}}
\bigl(\Delta_kW\mid\mathcal F_{k/n}\bigr)
=\mu_{k,n},
\qquad k=0,\ldots,n-1.
\label{eq:recovery-conditional-law}
\end{equation}
Under \(\mathbb P\), the increment \(\Delta_kW\) has law
\(p_n\) and is independent of \(\mathcal F_{k/n}\). Since \(\mu_{k,n}\) is
determined by \(\mathcal F_{k/n}\), for every bounded Borel function \(f\),
\begin{equation}
\mathbb E\bigl[\mathcal D_{k,n}f(\Delta_kW)\mid\mathcal F_{k/n}\bigr]
=\int_\R f\,\dd\mu_{k,n}.
\label{eq:recovery-one-step}
\end{equation}
With \(f\equiv1\), this shows that the positive process
\((D_k)_{k=0}^n\) is a \(\mathbb P\)-martingale with respect to
\((\mathcal F_{k/n})_{k=0}^n\), with \(D_0=1\). Hence
\(\mathbb E D_n=1\), so \(Q_n^{\eta,K}\) is a probability measure, and
\(\mathbb E[D_n\mid\mathcal F_{(k+1)/n}]=D_{k+1}=D_k\mathcal D_{k,n}\). Bayes'
formula and \eqref{eq:recovery-one-step} therefore give
\[
\mathbb E_{Q_n^{\eta,K}}\bigl[f(\Delta_kW)\mid\mathcal F_{k/n}\bigr]
=\frac{\mathbb E[D_n f(\Delta_kW)\mid\mathcal F_{k/n}]}{D_k}
=\mathbb E\bigl[\mathcal D_{k,n}f(\Delta_kW)\mid\mathcal F_{k/n}\bigr]
=\int_\R f\,\dd\mu_{k,n},
\]
which proves \eqref{eq:recovery-conditional-law}.

We now compute the entropy. By \eqref{eq:recovery-density},
\(\log D_n=\sum_{k<n}\log\mathcal D_{k,n}\), and each \(\log\mathcal D_{k,n}\) is a
quadratic polynomial in \(\Delta_kW\) whose coefficients are
\(\mathcal F_{k/n}\)-measurable and, by \eqref{eq:recovery-mean-bound} and
\(r\leq u^\eta\leq R\), bounded. By \eqref{eq:recovery-conditional-law},
it is therefore \(Q_n^{\eta,K}\)-integrable, and conditioning on
\(\mathcal F_{k/n}\) gives the chain rule for relative
entropy \cite{DupuisEllis},
\begin{equation}
H(Q_n^{\eta,K}\mid\mathbb P)
=\mathbb E_{Q_n^{\eta,K}}\log D_n
=
\sum_{k=0}^{n-1}\mathbb E_{Q_n^{\eta,K}}
H\bigl(\mu_{k,n}\mid p_n\bigr)
<\infty.
\label{eq:recovery-entropy-chain}
\end{equation}
For \(a\in\R\) and \(v>0\),
\(H(N(a,v)\mid p_n)=\frac12(nv-1-\log(nv)+na^2)\). By
\eqref{eq:ellv} and \eqref{eq:recovery-kernels}, this gives
\[
H(\mu_{k,n}\mid p_n)
=\ell(u^\eta(k/n))+\frac{(m_{k,n}^{\eta,K})^2}{2n},
\qquad k\in\mathcal K_n^0,
\]
and, for \(k=\lceil ns_j\rceil\),
\[
H(\mu_{k,n}\mid p_n)
=\frac12\Bigl[
n\bigl(\delta_j+\operatorname{sgn}(\delta_j)n^{-1/4}\bigr)^2
+n^{-1}-1+\log n
\Bigr]
=\frac n2\delta_j^2+o(n).
\]
Insert these into \eqref{eq:recovery-entropy-chain}. By
\eqref{eq:recovery-mean-bound}, the terms
\((m_{k,n}^{\eta,K})^2/(2n)\) contribute at most \(B_K^2/(2n)\) to
\(\frac1nH(Q_n^{\eta,K}\mid\mathbb P)\). The function
\(\ell(u^\eta(\cdot))\) is piecewise constant and bounded by
\(\max_{[r,R]}\ell\), and \(\mathcal K_n^0\) misses only \(N\) indices.
Hence \(\frac1n\sum_{k\in\mathcal K_n^0}\ell(u^\eta(k/n))\to
\int_0^1\ell(u^\eta(t))\,\dd t\), and
\[
\lim_{n\to\infty}\frac1nH(Q_n^{\eta,K}\mid\mathbb P)
=\int_0^1\ell(u^\eta(t))\,\dd t+\frac12\sum_{j=1}^N\delta_j^2.
\]
This proves \eqref{eq:recovery-entropy-bound}.

\medskip\noindent\textit{Step 2. We couple the scheme with the restoring diffusion and reduce \eqref{eq:recovery-joint-convergence} to two comparisons.} For \(j=1,\ldots,N\), the
prescribed update of the scheme takes place on the grid cell
\([p_j^n,r_j^n)\), where
\begin{equation}
p_j^n:=\frac{\lceil ns_j\rceil}{n},\qquad
r_j^n:=\frac{\lceil ns_j\rceil+1}{n},\qquad
0\leq p_j^n-s_j\leq\frac1n,
\qquad
0<r_j^n-s_j\leq\frac2n.
\label{eq:joint-jump-grid-times}
\end{equation}
Set \(r_0^n:=0\) and \(p_{N+1}^n:=s_{N+1}:=1\), and
increase \(n_0\), if necessary, so that \(r_j^n<s_{j+1}\) for
\(j=1,\ldots,N\). Let
\begin{equation}
B_n^\circ:=\sum_{k\in\mathcal K_n^0}(\Delta_kW)^2
\delta_{(k/n,\iota(Z_k^n))},
\qquad
\zeta_n^\circ:=\left(
A_n,B_n^\circ,
\bigl((p_j^n,\iota(Z_{n p_j^n}^n),\delta_j)\bigr)_{j=1}^N,
\iota(Z_n^n)\right).
\label{eq:joint-auxiliary-control}
\end{equation}
Thus \(\zeta_n^\circ\) is the empirical control in which the
marks are placed exactly at the prescribed indices and carry the target
increments.

By \eqref{eq:recovery-conditional-law}, each of the
variables
\begin{equation}
G_k:=
\begin{cases}
\displaystyle
\sqrt{\frac{n}{u^\eta(k/n)}}
\left(\Delta_kW-\frac{m_{k,n}^{\eta,K}}{n}\right),
&k\in\mathcal K_n^0,\\[2mm]
\displaystyle
n\left(\Delta_{k}W-\delta_j
-\operatorname{sgn}(\delta_j)n^{-1/4}\right),
&k=n p_j^n,\quad j=1,\ldots,N,
\end{cases}
\label{eq:joint-normal-variables}
\end{equation}
is \(\mathcal F_{(k+1)/n}\)-measurable and has conditional
law \(N(0,1)\) given \(\mathcal F_{k/n}\). By successive conditioning,
\(G_0,\ldots,G_{n-1}\) are independent standard normal variables under
\(Q_n^{\eta,K}\).

We drive
\eqref{eq:restoring-sde} with a Brownian motion whose grid increments are
the variables \(G_k/\sqrt n\). On an extension of
\((\mathcal W,\mathcal B(\mathcal W),Q_n^{\eta,K})\), let
\(\mathbb B_0^n,\ldots,\mathbb B_{n-1}^n\) be independent
standard Brownian bridges on \([0,1]\), independent of
\((G_0,\ldots,G_{n-1})\).
We continue to denote the extended measure by \(Q_n^{\eta,K}\).
For \(k=0,\ldots,n-1\), define
\begin{equation}
\widetilde W_t^n:=\frac1{\sqrt n}\left[
\sum_{\ell<k}G_\ell+n(t-k/n)G_k
+\mathbb B_k^n\bigl(n(t-k/n)\bigr)\right],
\quad t\in[k/n,(k+1)/n].
\label{eq:joint-brownian-coupling}
\end{equation}
Then \(\widetilde W^n\) is a standard Brownian motion under
\(Q_n^{\eta,K}\). Since \(m_{k,n}^{\eta,K}\) is a function
of \(\Delta_{<k}W\), induction on \(k\) in
\eqref{eq:joint-normal-variables} and
\eqref{eq:transformed-euler-update} shows that
\begin{equation}
\sigma(Z_k^n,\Delta_{<k}W)
\subseteq\sigma(G_0,\ldots,G_{k-1})
\subseteq\sigma(\widetilde W_s^n:0\leq s\leq k/n),
\qquad k=0,\ldots,n.
\label{eq:joint-grid-adaptedness}
\end{equation}
Let \(\Gamma^{\eta,K,n}\) solve \eqref{eq:restoring-sde}
with \(\widetilde W\) replaced by \(\widetilde W^n\), and let
\(\zeta^{\eta,K,n}\) be the control in
\eqref{eq:restoring-control-convergence} with \(\Gamma^{\eta,K}\) replaced
by \(\Gamma^{\eta,K,n}\). Since \eqref{eq:restoring-sde} has a unique
strong solution,
\begin{equation}
\mathcal L_{Q_n^{\eta,K}}(\Gamma^{\eta,K,n})
=\mathcal L_{\widetilde{\mathbb P}}(\Gamma^{\eta,K}),
\label{eq:joint-diffusion-law}
\end{equation}
and hence
\begin{equation}
\mathcal L_{Q_n^{\eta,K}}(\zeta^{\eta,K,n})
=\mathcal L_{\widetilde{\mathbb P}}(\zeta^{\eta,K}).
\label{eq:joint-control-law}
\end{equation}

We shall prove that
\begin{equation}
\lim_{n\to\infty}\mathbb E_{Q_n^{\eta,K}}
d_{\mathcal C_\varepsilon}(\zeta_n^\varepsilon,\zeta_n^\circ)=0,
\label{eq:joint-cutoff-control-comparison}
\end{equation}
\begin{equation}
d_{\mathcal C_\varepsilon}(\zeta_n^\circ,\zeta^{\eta,K,n})
\longrightarrow0
\quad\text{in }Q_n^{\eta,K}\text{-probability}.
\label{eq:joint-auxiliary-control-comparison}
\end{equation}
These imply \eqref{eq:recovery-joint-convergence}. Since
\(d_{\mathcal C_\varepsilon}\leq1\),
\eqref{eq:joint-auxiliary-control-comparison} also holds in mean. Hence,
for every bounded Lipschitz function
\(\varphi:\mathcal C_\varepsilon\to\R\), \eqref{eq:joint-control-law} gives
\begin{equation}
\begin{aligned}
&\left|\mathbb E_{Q_n^{\eta,K}}\varphi(\zeta_n^\varepsilon)
-\widetilde{\mathbb E}\varphi(\zeta^{\eta,K})\right|\\
&\quad\leq\operatorname{Lip}(\varphi)\Bigl[
\mathbb E_{Q_n^{\eta,K}}
 d_{\mathcal C_\varepsilon}(\zeta_n^\varepsilon,\zeta_n^\circ)
+\mathbb E_{Q_n^{\eta,K}}
 d_{\mathcal C_\varepsilon}(\zeta_n^\circ,\zeta^{\eta,K,n})
\Bigr]\longrightarrow0,
\end{aligned}
\label{eq:joint-final-test-convergence}
\end{equation}
and bounded Lipschitz functions determine weak convergence
on the metric space \(\mathcal C_\varepsilon\). Step~3 proves \eqref{eq:joint-cutoff-control-comparison}, and Step~4
reduces \eqref{eq:joint-auxiliary-control-comparison} to a comparison of
states and a comparison of measures, which are proved in
Appendix~\ref{app:euler-comparison}.

\medskip\noindent\textit{Step 3. We prove \eqref{eq:joint-cutoff-control-comparison}.} By \eqref{eq:joint-normal-variables},
\begin{equation}
\mathbb E_{Q_n^{\eta,K}}
|\Delta_{n p_j^n}W-\delta_j|^2
=n^{-1/2}+n^{-2},
\qquad j=1,\ldots,N.
\label{eq:joint-prescribed-increment-L2}
\end{equation}
Consider the event that the increments reaching the cutoff
are exactly the prescribed ones,
\begin{equation}
F_n:=\bigl\{
|\Delta_{n p_j^n}W|\geq\varepsilon\text{ for }j=1,\ldots,N,
\text{ and }|\Delta_kW|<\varepsilon\text{ for }k\in\mathcal K_n^0
\bigr\}.
\label{eq:joint-index-identification-event}
\end{equation}
On \(F_n\), \eqref{eq:empirical-variance} and
\eqref{eq:empirical-marks} give \(B_n^\varepsilon=B_n^\circ\) and
\(\mathfrak m_n^\varepsilon
=((p_j^n,\iota(Z_{n p_j^n}^n),\Delta_{n p_j^n}W))_{j=1}^N\). Thus
\(\zeta_n^\varepsilon\) and \(\zeta_n^\circ\) differ only in their
increment coordinates, and \eqref{eq:joint-control-metric} gives
\(d_{\mathcal C_\varepsilon}(\zeta_n^\varepsilon,\zeta_n^\circ)
\leq\sum_{j}|\Delta_{n p_j^n}W-\delta_j|\). Since
\(d_{\mathcal C_\varepsilon}\leq1\), the Cauchy--Schwarz inequality and
\eqref{eq:joint-prescribed-increment-L2} give
\[
\mathbb E_{Q_n^{\eta,K}}
d_{\mathcal C_\varepsilon}(\zeta_n^\varepsilon,\zeta_n^\circ)
\leq Q_n^{\eta,K}(F_n^c)+N\bigl(n^{-1/2}+n^{-2}\bigr)^{1/2}.
\]
It remains to show that \(Q_n^{\eta,K}(F_n^c)\to0\). For
all sufficiently large \(n\), \(B_K/n\leq\varepsilon/2\). Hence, by
\eqref{eq:recovery-mean-bound}, \eqref{eq:recovery-conditional-law}, the
Gaussian tail estimate, and a union bound,
\begin{equation}
Q_n^{\eta,K}\left\{
\max_{k\in\mathcal K_n^0}|\Delta_kW|\geq\varepsilon
\right\}\leq Q_n^{\eta,K}\left\{
\max_{k\in\mathcal K_n^0}
\bigg|\Delta_kW-\frac{m_{k,n}^{\eta,K}}{n}\bigg|
\geq\frac\varepsilon2
\right\}
\leq2n\exp\left(-\frac{n\varepsilon^2}{8R}\right).
\label{eq:joint-no-additional-indices}
\end{equation}
For \(j=1,\ldots,N\), \eqref{eq:joint-normal-variables}
gives
\(\operatorname{sgn}(\delta_j)\Delta_{n p_j^n}W
=|\delta_j|+n^{-1/4}+n^{-1}\operatorname{sgn}(\delta_j)G_{n p_j^n}\).
Since \(|\delta_j|\geq\varepsilon\) by \eqref{eq:ordered-mark-list},
\[
Q_n^{\eta,K}\{|\Delta_{n p_j^n}W|<\varepsilon\}
\leq Q_n^{\eta,K}\left\{
\operatorname{sgn}(\delta_j)G_{n p_j^n}<-n^{3/4}
\right\}\leq\exp(-n^{3/2}/2).
\]
Together with \eqref{eq:joint-no-additional-indices}, this
gives \(Q_n^{\eta,K}(F_n^c)\leq2n\exp(-n\varepsilon^2/(8R))
+N\exp(-n^{3/2}/2)\to0\), and \eqref{eq:joint-cutoff-control-comparison}
follows.

\medskip\noindent\textit{Step 4. We reduce \eqref{eq:joint-auxiliary-control-comparison} to a comparison of states and a comparison of measures.}
We now reduce \eqref{eq:joint-auxiliary-control-comparison}.
Let \(Z^n(t):=Z_k^n\) for \(t\in[k/n,(k+1)/n)\), \(k=0,\ldots,n-1\), and
\(Z^n(1):=Z_n^n\), so that
\(A_n=\dd t\,\delta_{\iota(Z^n(t))}\) by \eqref{eq:empirical-occupation}.
Set
\begin{equation}
A_{\Gamma,n}:=\dd t\,\delta_{\iota(\Gamma_t^{\eta,K,n})},\qquad
B_{\Gamma,n}:=u^\eta(t)\,\dd t\,
\delta_{\iota(\Gamma_t^{\eta,K,n})}.
\label{eq:joint-auxiliary-measures}
\end{equation}
The controls \(\zeta_n^\circ\) and \(\zeta^{\eta,K,n}\) both
belong to \(\mathcal C_{\varepsilon,N}\) and have identical increment coordinates. Thus
\eqref{eq:joint-control-metric} and \(|p_j^n-s_j|\leq1/n\) give
\begin{equation}
\begin{aligned}
d_{\mathcal C_\varepsilon}(\zeta_n^\circ,\zeta^{\eta,K,n})
\leq{}&d_{\mathrm w}(A_n,A_{\Gamma,n})
+d_{\mathrm w}(B_n^\circ,B_{\Gamma,n})
+\frac Nn\\
&+\sum_{j=1}^{N+1}
d_{K_\sigma}\bigl(
\iota(Z_{n p_j^n}^n),\iota(\Gamma_{s_j-}^{\eta,K,n})\bigr),
\end{aligned}
\label{eq:joint-auxiliary-control-comparison0}
\end{equation}
where the term \(j=N+1\) is the terminal label, since
\(p_{N+1}^n=s_{N+1}=1\) and \(\Gamma^{\eta,K,n}\) has no update at \(1\).
Let \(E_n:=\bigcup_{j=1}^N[s_j,r_j^n)\), so that \(|E_n|\leq2N/n\) by
\eqref{eq:joint-jump-grid-times}. In Appendix~\ref{app:euler-comparison} we prove the state comparison
\begin{equation}
\mathbb E_{Q_n^{\eta,K}}\bigg[
\sup_{t\in[0,1]\setminus E_n}
|Z^n(t)-\Gamma_t^{\eta,K,n}|^2
+\sum_{j=1}^{N+1}
|Z_{n p_j^n}^n-\Gamma_{s_j-}^{\eta,K,n}|^2
\bigg]\leq Cn^{-1/2},
\label{eq:joint-coupled-state-comparison}
\end{equation}
and the measure comparison
\begin{equation}
d_{\mathrm w}(A_n,A_{\Gamma,n})
+d_{\mathrm w}(B_n^\circ,B_{\Gamma,n})
\longrightarrow0
\quad\text{in }Q_n^{\eta,K}\text{-probability}.
\label{eq:joint-occupation-comparison}
\end{equation}

To pass from states to labels, we use the following
elementary fact. Let \(T_n\subset[0,1]\), let \(\Lambda\) be a metric
space, and let \(g:[0,1]\times\R\to\Lambda\) be continuous. If
\(\sup_{t\in T_n}|x_n(t)-y_n(t)|\to0\) in probability and the laws of
\(\sup_{t\in T_n}|y_n(t)|\) are tight, then
\begin{equation}
\sup_{t\in T_n}
d_\Lambda\bigl(g(t,x_n(t)),g(t,y_n(t))\bigr)\longrightarrow0
\quad\text{in probability}.
\label{eq:joint-localization}
\end{equation}
For a proof, fix \(a>0\) and \(L\geq1\). By uniform
continuity of \(g\) on \([0,1]\times[-L-1,L+1]\), there is
\(\theta\in(0,1]\) such that the supremum in \eqref{eq:joint-localization}
is at most \(a\) unless \(\sup_{T_n}|y_n|>L\) or
\(\sup_{T_n}|x_n-y_n|>\theta\). Letting first \(n\to\infty\) and then
\(L\to\infty\) proves \eqref{eq:joint-localization}.

By \eqref{eq:joint-diffusion-law}, the law of
\(\sup_{0\leq t\leq1}|\Gamma_t^{\eta,K,n}|\) does not depend on \(n\), and
it is finite because \(\Gamma^{\eta,K}\) has c\`adl\`ag paths. Hence
\eqref{eq:joint-coupled-state-comparison},
\eqref{eq:joint-localization} with \(g(t,x)=\iota(x)\), and Markov's
inequality show that the label terms in
\eqref{eq:joint-auxiliary-control-comparison0} tend to zero in
probability. Together with \eqref{eq:joint-occupation-comparison}, this
proves \eqref{eq:joint-auxiliary-control-comparison}, and hence
\eqref{eq:recovery-joint-convergence} by Step~2.

\end{proof}

\subsection{Proof of Theorem~\ref{thm:cutoffldp}}
\label{subsec:proof-cutoffldp}

The upper bound combines the dual representation
\eqref{eq:Isv-dual} with an exponential Chebyshev inequality.
For a test function \(\vartheta\) in \eqref{eq:Isv-dual}, the
weights \(\theta_k:=\vartheta(t_k,\iota(Z_k^n))\) are
\(\mathcal F_{t_k}\)-measurable, and the sum
\(\sum_k[\theta_kn(\Delta_kW)^2+\frac12\log(1-2\theta_k)]\) equals \(n\)
times the dual functional in \eqref{eq:Isv-dual} at
\((A_n,B_n)\), up to a discretization error. The next lemma shows that
the exponential of such a sum has expectation one.

\begin{lemma}[Exponential identity for squared Gaussian increments]
\label{lem:gaussianexp}
Let \(n\geq1\), and let \(\theta_k\), \(0\leq k<n\), be
\(\mathcal F_{t_k}\)-measurable random variables with \(\theta_k<\frac12\)
almost surely. Then, for every \(0\leq k<n\),
\begin{equation}
\mathbb E\left[
\exp\left\{
\theta_k n(\Delta_kW)^2
+\frac12\log(1-2\theta_k)
\right\}
\,\middle|\,\mathcal F_{t_k}
\right]
=1,
\label{eq:conditional-gaussian-identity}
\end{equation}
and hence
\begin{equation}
\mathbb E\exp\left\{
\sum_{k=0}^{n-1}
\left[
\theta_k n(\Delta_kW)^2
+\frac12\log(1-2\theta_k)
\right]
\right\}
=1.
\label{eq:predictable-gaussian-identity}
\end{equation}
\end{lemma}
\begin{proof}
Under \(\mathbb P\), the variable \(\sqrt n\,\Delta_kW\) is standard normal
and independent of \(\mathcal F_{t_k}\), and
\(\mathbb Ee^{\theta G^2}=(1-2\theta)^{-1/2}\) for a standard normal \(G\)
and \(\theta<\frac12\). Since \(\theta_k\) is
\(\mathcal F_{t_k}\)-measurable, it follows that
\[
\mathbb E\bigl[\exp\{\theta_kn(\Delta_kW)^2\}\,\big|\,\mathcal F_{t_k}\bigr]
=(1-2\theta_k)^{-1/2},
\]
which is \eqref{eq:conditional-gaussian-identity}. For
\(m=0,\ldots,n\), let
\[
M_m:=\exp\Bigl\{\sum_{k<m}\Bigl[\theta_kn(\Delta_kW)^2
+\frac12\log(1-2\theta_k)\Bigr]\Bigr\},
\qquad M_0=1.
\]
Each \(M_m\) is nonnegative and \(\mathcal F_{t_m}\)-measurable, and
\(M_{m+1}/M_m\) is the exponential in
\eqref{eq:conditional-gaussian-identity} with \(k=m\). Hence
\eqref{eq:conditional-gaussian-identity} yields
\(\mathbb E[M_{m+1}\mid\mathcal F_{t_m}]=M_m\), and the tower property
gives \(\mathbb EM_n=M_0=1\), which is
\eqref{eq:predictable-gaussian-identity}.
\end{proof}

\begin{proof}[Proof of Theorem~\ref{thm:cutoffldp}]
\medskip\noindent\textit{Step 1. We prove the lower bound for open sets.}
Let \(O\subset\mathcal C_\varepsilon\) be open. We show
that, for every \(\zeta=(A,B,\mathfrak m,\xi_f)\in O\) with
\(I_\varepsilon(\zeta)<\infty\),
\begin{equation}
\liminf_{n\to\infty}
\frac1n\log
\mathbb P\{\zeta_n^\varepsilon\in O\}
\geq
-I_\varepsilon(\zeta).
\label{eq:local-lower}
\end{equation}
Taking the supremum over \(\zeta\in O\) then gives the
lower bound.

Fix such a \(\zeta\) and \(0<\eta<1\). By
Lemma~\ref{lem:regularize}, there are regular pairs \((A_m,B_m)\) with
\((A_m,B_m)\Rightarrow(A,B)\) and \(I_{\rm sv}(A_m,B_m)\to
I_{\rm sv}(A,B)\). Let \(N\) be the length of \(\mathfrak m\). With the marks
and terminal tag fixed, \((A_m,B_m,\mathfrak m,\xi_f)\to\zeta\) in
\(\mathcal C_{\varepsilon,N}\). Since \(O\) is open and
\(I_\varepsilon(A_m,B_m,\mathfrak m,\xi_f)=I_{\rm sv}(A_m,B_m)
+\frac12\sum_j\delta_j^2\to I_\varepsilon(\zeta)\) by \eqref{eq:IAB}, there is \(m\) such that
\(\zeta_m:=(A_m,B_m,\mathfrak m,\xi_f)\in O\) and
\(I_\varepsilon(\zeta_m)\leq I_\varepsilon(\zeta)+\eta\). Lemma~\ref{lem:marked-approximation}, applied to \(\zeta_m\), \(O\)
and \(\eta\), yields functions \(q^\eta,u^\eta\), times \(s_1<\cdots<s_N\), and the control \(\zeta^\eta\)
in \eqref{eq:zetaeta}, with the same increments
\(\delta_1,\ldots,\delta_N\) as \(\zeta\), such that
\begin{equation}
\zeta^\eta\in O,
\qquad
\int_0^1\ell(u^\eta(t))\,\dd t
\leq I_{\rm sv}(A_m,B_m)+\eta.
\label{eq:finite-recovery-target}
\end{equation}
By Lemma~\ref{lem:restoring-tracking},
\(\zeta^{\eta,K}\to\zeta^\eta\) in \(\widetilde{\mathbb P}\)-probability as
\(K\to\infty\). Since \(O\) is open and contains \(\zeta^\eta\), we can fix
\(K\geq1\) such that
\(\widetilde{\mathbb P}\{\zeta^{\eta,K}\in O\}\geq1-\eta\). Let
\(Q_n^\eta:=Q_n^{\eta,K}\) be the measures of
Lemma~\ref{lem:recovery-tilted-law} for these data. By
\eqref{eq:recovery-joint-convergence} and the portmanteau theorem,
\(\liminf_{n}Q_n^\eta\{\zeta_n^\varepsilon\in O\}
\geq\widetilde{\mathbb P}\{\zeta^{\eta,K}\in O\}\geq1-\eta\). By
\eqref{eq:recovery-entropy-bound}, \eqref{eq:finite-recovery-target} and
\eqref{eq:IAB},
\[
\begin{aligned}
\limsup_{n\to\infty}\frac1nH(Q_n^\eta\mid\mathbb P)
&\leq\int_0^1\ell(u^\eta(t))\,\dd t+\frac12\sum_{j=1}^N\delta_j^2
\leq I_{\rm sv}(A_m,B_m)+\frac12\sum_{j=1}^N\delta_j^2+\eta\\
&=I_\varepsilon(\zeta_m)+\eta\leq I_\varepsilon(\zeta)+2\eta.
\end{aligned}
\]
We have thus shown that
\begin{equation}
\liminf_{n\to\infty}Q_n^\eta\{\zeta_n^\varepsilon\in O\}\geq1-\eta,
\qquad
\limsup_{n\to\infty}
\frac1nH(Q_n^\eta\mid\mathbb P)
\leq I_\varepsilon(\zeta)+2\eta.
\label{eq:tag-5-35e}
\end{equation}
We transfer this to \(\mathbb P\) by the entropy
inequality. Let \(Q\ll\mathbb P\) and let \(E\) be an event with
\(Q(E)>0\). By~\cite[Lemma~1.4.3(g)]{DupuisEllis}, \(H(Q\mid\mathbb P)\)
is the supremum of the discrete relative entropies over all finite
measurable partitions. Taking \(\{E,E^c\}\) gives
\[
H(Q\mid\mathbb P)
\geq
Q(E)\log\frac{Q(E)}{\mathbb P(E)}
+
(1-Q(E))
\log\frac{1-Q(E)}{1-\mathbb P(E)}.
\]
Since \(q\log q+(1-q)\log(1-q)\geq-\log2\) for
\(q\in[0,1]\) and \(-\log(1-\mathbb P(E))\geq0\), this yields
\begin{equation}
\log\mathbb P(E)
\geq
-\frac{H(Q\mid\mathbb P)+\log2}{Q(E)}.
\label{eq:tag-5-35f}
\end{equation}
Applying \eqref{eq:tag-5-35f} with \(Q=Q_n^\eta\) and
\(E=\{\zeta_n^\varepsilon\in O\}\), and using \eqref{eq:tag-5-35e}, we
obtain
\[
\liminf_{n\to\infty}
\frac1n\log
\mathbb P\{\zeta_n^\varepsilon\in O\}
\geq
-\frac{I_\varepsilon(\zeta)+2\eta}{1-\eta}.
\]
Letting \(\eta\downarrow0\) proves \eqref{eq:local-lower}.

\medskip\noindent\textit{Step 2. We prove a local upper bound around every control.}
We show that for every \(\zeta\in\mathcal C_\varepsilon\)
and every \(r<I_\varepsilon(\zeta)\), there is an open neighbourhood \(U_\zeta\) of
\(\zeta\) such that
\begin{equation}
\limsup_{n\to\infty}
\frac1n\log
\mathbb P\{\zeta_n^\varepsilon\in U_\zeta\}
\leq-r.
\label{eq:local-upper}
\end{equation}
Let
\(\zeta=(A,B,\mathfrak m,\xi_f)\) with
\(\mathfrak m=((t_j,\xi_j,\delta_j))_{j=1}^N\). For
\(\vartheta\in C(K_\sigma^0)\) with \(\sup\vartheta<1/2\), define $g_\vartheta\in C(K_\sigma^0)$ by
\(g_\vartheta(t,\xi):=\frac12\log(1-2\vartheta(t,\xi))\) and set
\[
L_\vartheta(A',B')
:=\int_{K_\sigma^0}\vartheta\,\dd B'
+\int_{K_\sigma^0}g_\vartheta\,\dd A',
\qquad (A',B')\in\mathcal Q_\sigma.
\]
Since \(I_\varepsilon(\zeta)=I_{\rm sv}(A,B)+\frac12\sum_j\delta_j^2\)
by \eqref{eq:IAB}, the dual representation
\eqref{eq:Isv-dual} gives such a \(\vartheta\) with
\(L_\vartheta(A,B)+\frac12\sum_j\delta_j^2>r\). Choose bounded open
intervals \(D_j\ni\delta_j\) so small that \(L_\vartheta(A,B)+c>r\), where
\(c:=\frac12\sum_{j}\inf_{x\in D_j}x^2\), and fix
\(r'\in(r,L_\vartheta(A,B)+c)\). Since \(L_\vartheta\) is continuous on
\(\mathcal Q_\sigma\), the set \(U_\zeta\) of controls
\((A',B',((t_j',\xi_j',\delta_j'))_{j=1}^N,\xi_f')\) with exactly \(N\)
marks, \(\delta_j'\in D_j\) for all \(j\), and
\begin{equation}
L_\vartheta(A',B')+c>r'
\label{eq:local-upper-margin}
\end{equation}
is an open neighbourhood of \(\zeta\).

On \(\{\zeta_n^\varepsilon\in U_\zeta\}\), exactly \(N\)
increments reach the cutoff. For indices \(0\leq k_1<\cdots<k_N<n\), let
\begin{equation}
E_n(k_1,\ldots,k_N)
:=
\left\{
\zeta_n^\varepsilon\in U_\zeta,\;
L_n^\varepsilon=\{k_1,\ldots,k_N\}
\right\}.
\label{eq:fixed-mark-event}
\end{equation}
There are at most \(n^N\) such events, and their union is
\(\{\zeta_n^\varepsilon\in U_\zeta\}\). Therefore it suffices to prove,
uniformly in \(k_1,\ldots,k_N\),
\begin{equation}
\mathbb P\bigl(E_n(k_1,\ldots,k_N)\bigr)
\leq
C\,n^{N/2}\exp\{-nr'+o(n)\},
\label{eq:local-upper-exp}
\end{equation}
since then
\(\mathbb P\{\zeta_n^\varepsilon\in U_\zeta\}\leq
Cn^{3N/2}\exp\{-nr'+o(n)\}\), and \eqref{eq:local-upper} follows from
\(r'>r\).

Fix \(k_1<\cdots<k_N\), let \(S:=\{k_1,\ldots,k_N\}\), and
set \(\theta_k:=\vartheta(t_k,\iota(Z_k^n))\), which is
\(\mathcal F_{t_k}\)-measurable. Consider
\[
Y_k:=
\exp\left\{
\theta_k\,n(\Delta_kW)^2
+g_\vartheta(t_k,\iota(Z_k^n))
\right\}
\mathbf1_{\{|\Delta_kW|<\varepsilon\}},
\quad k\notin S,
\qquad
H_j:=\mathbf1_{\{\Delta_{k_j}W\in D_j\}}.
\]
By Lemma~\ref{lem:gaussianexp},
\(\mathbb E[Y_k\mid\mathcal F_{t_k}]\leq1\). Since \(\Delta_{k_j}W\) is
\(N(0,n^{-1})\) and independent of \(\mathcal F_{t_{k_j}}\), its density
is at most \(\sqrt{n/(2\pi)}\exp\{-\frac n2\inf_{D_j}x^2\}\) on \(D_j\),
so that \(\mathbb E[H_j\mid\mathcal F_{t_{k_j}}]\leq
C_j\sqrt n\exp\{-\frac n2\inf_{D_j}x^2\}\) with
\(C_j:=|D_j|/\sqrt{2\pi}\). Taking successive conditional expectations given
\(\mathcal F_{t_{n-1}},\mathcal F_{t_{n-2}},\ldots,\mathcal F_{t_0}=\mathcal F_0\) yields
\begin{equation}
\mathbb E\left[
\prod_{k\notin S}Y_k
\prod_{j=1}^NH_j
\right]
\leq
C\,n^{N/2}e^{-nc}.
\label{eq:local-upper-tower}
\end{equation}
We now bound the product from below on
\(E_n(k_1,\ldots,k_N)\). There, all indicators equal one, and
\(B_n^\varepsilon=\sum_{k\notin S}(\Delta_kW)^2\delta_{(t_k,\iota(Z_k^n))}\)
by \eqref{eq:empirical-variance}. Hence
\(\sum_{k\notin S}\theta_k\,n(\Delta_kW)^2=n\int\vartheta\,\dd
B_n^\varepsilon\). Moreover, let \(\omega_\vartheta\) be the modulus of
continuity of \(g_\vartheta\) in the time variable, uniformly in
\(\xi\in K_\sigma\). Since \(A_n\) has the state tag \(\iota(Z_k^n)\) on
\([t_k,t_{k+1})\) by \eqref{eq:empirical-occupation},
\[
\bigg|\sum_{k\notin S}g_\vartheta(t_k,\iota(Z_k^n))
-n\int g_\vartheta\,\dd A_n\bigg|
\leq n\,\omega_\vartheta(1/n)+N\|g_\vartheta\|_\infty=o(n),
\]
uniformly in the values of \(Z_k^n\) and in
\(k_1,\ldots,k_N\).

On \(E_n(k_1,\ldots,k_N)\), the definitions of
\(Y_k\) and \(H_j\) give
\[
\log\!\left(\prod_{k\notin S}Y_k\prod_{j=1}^{N}H_j\right)
=
n\sum_{k\notin S}\theta_k(\Delta_k W)^2
+\sum_{k\notin S}g_\vartheta(t_k,\iota(Z_k^n)).
\]
Using the preceding estimate for the second sum and
\eqref{eq:local-upper-margin}, we obtain after exponentiation
\[
\prod_{k\notin S}Y_k\prod_{j=1}^NH_j
\geq\exp\bigl\{nL_\vartheta(A_n,B_n^\varepsilon)-o(n)\bigr\}
\geq\exp\bigl\{n(r'-c)-o(n)\bigr\}.
\]
Hence Markov's inequality and \eqref{eq:local-upper-tower} yield
\[
\mathbb P\bigl(E_n(k_1,\ldots,k_N)\bigr)
\leq e^{-n(r'-c)+o(n)}\,
\mathbb E\Bigl[\prod_{k\notin S}Y_k\prod_{j=1}^NH_j\Bigr]
\leq C\,n^{N/2}e^{-nr'+o(n)},
\]
which is \eqref{eq:local-upper-exp}.

\medskip\noindent\textit{Step 3. We show that the empirical controls are exponentially tight and that \(I_\varepsilon\) is a good rate function.}
For \(M>1\), let \(\mathcal K_M\) be the set
of controls
\((A,B,((t_j,\xi_j,\delta_j))_{j=1}^N,\xi_f)\in\mathcal C_\varepsilon\)
with \(B(K_\sigma^0)\leq M\), \(N\leq M/\varepsilon^2\), and
\(|\delta_j|\leq\sqrt M\) for all \(j\). The set \(\mathcal K_M\) is
compact. Only finitely many intersections
\(\mathcal K_M\cap\mathcal C_{\varepsilon,N}\) are nonempty, since
\(N\leq M/\varepsilon^2\) on \(\mathcal K_M\). Each is closed in
\(\mathcal C_{\varepsilon,N}\). In each such intersection, the \(A\)-coordinates lie in the compact set of probability measures
on \(K_\sigma^0\) with time marginal \(\lambda\), the \(B\)-coordinates lie
in a compact set by Prokhorov's theorem, and the remaining coordinates lie
in compact sets.

On \(\{\mu_n([0,1])\leq M\}\), we have
\(\zeta_n^\varepsilon\in\mathcal K_M\), since each \(k\in
L_n^\varepsilon\) contributes at least \(\varepsilon^2\) to
\(\sum_k(\Delta_kW)^2=\mu_n([0,1])\), every increment satisfies
\(|\Delta_kW|\leq\mu_n([0,1])^{1/2}\), and
\(B_n^\varepsilon(K_\sigma^0)\leq\mu_n([0,1])\). Hence, by
\eqref{eq:tag-13},
\begin{equation}
\limsup_{n\to\infty}\frac1n\log
\mathbb P\{\zeta_n^\varepsilon\notin\mathcal K_M\}
\leq-\ell(M)\longrightarrow-\infty,
\qquad M\to\infty.
\label{eq:cutoff-exp-tight}
\end{equation}

We next show that \(I_\varepsilon\) is good on \(\mathcal C_\varepsilon\). The function \(I_{\rm sv}\)
is lower semicontinuous by Proposition~\ref{prop:sv}, and the mark cost
\(\frac12\sum_j\delta_j^2\) is continuous on \(\mathcal C_{\varepsilon,N}\) for every \(N\).
Their sum \(I_\varepsilon\) is therefore lower semicontinuous
on each \(\mathcal C_{\varepsilon,N}\). If \(\zeta_m\to\zeta\) in
\(\mathcal C_\varepsilon\) and \(\zeta\in\mathcal C_{\varepsilon,N}\),
then \(\zeta_m\in\mathcal C_{\varepsilon,N}\) for all sufficiently large \(m\),
since controls with different numbers of marks are at distance \(1\).
Lower semicontinuity on \(\mathcal C_{\varepsilon,N}\) gives
\(I_\varepsilon(\zeta)\leq\liminf_{m\to\infty}I_\varepsilon(\zeta_m)\),
which proves lower semicontinuity on \(\mathcal C_\varepsilon\). Let \(a<\infty\) and
\(\zeta=(A,B,((t_j,\xi_j,\delta_j))_{j=1}^N,\xi_f)\) with
\(I_\varepsilon(\zeta)\leq a\). Then \(B(K_\sigma^0)\leq 4a+2\log2\) by \eqref{eq:IsvAB},
and \(\varepsilon^2\leq\delta_j^2\leq\sum_i\delta_i^2\leq2a\) gives
\(N\leq2a/\varepsilon^2\) and \(|\delta_j|\leq\sqrt{2a}\). Thus the sublevel set
\(\{\zeta\in\mathcal C_\varepsilon:I_\varepsilon(\zeta)\leq a\}\) is a closed subset of \(\mathcal K_M\) with
\(M:=\max\{4a+2\log2,2a,2\}\), and hence compact.

\medskip\noindent\textit{Step 4. We combine Steps~2 and~3 into the upper bound for closed sets.}
Let \(F\subset\mathcal C_\varepsilon\) be closed, let
\(r<\inf_F I_\varepsilon\), and let \(M>1\). The set
\(F\cap\mathcal K_M\) is compact, and \(I_\varepsilon>r\) on it. By Step~2, it is
covered by finitely many neighbourhoods \(U_{\zeta_1},\ldots,U_{\zeta_L}\)
satisfying \eqref{eq:local-upper}. Since
\[
\mathbb P\{\zeta_n^\varepsilon\in F\}
\leq\sum_{\ell=1}^L\mathbb P\{\zeta_n^\varepsilon\in U_{\zeta_\ell}\}
+\mathbb P\{\zeta_n^\varepsilon\notin\mathcal K_M\},
\]
and the number \(L\) of neighbourhoods is finite, \eqref{eq:local-upper} and \eqref{eq:cutoff-exp-tight}
yield
\[
\limsup_{n\to\infty}\frac1n\log
\mathbb P\{\zeta_n^\varepsilon\in F\}
\leq\max\{-r,-\ell(M)\}.
\]
Letting \(M\to\infty\) and then \(r\uparrow\inf_F I_\varepsilon\)
proves the upper bound for closed sets.
\end{proof}

\subsection*{Acknowledgment}
This work was supported by the National Natural Science Foundation of China (Nos.\ 12671472 and 12201228) and the Fundamental Research Funds for the Central Universities (No.\ 2026BRSXB006). The author acknowledges that ChatGPT was used as an AI tool for exploratory discussions and for language polishing during the preparation of this work. The author independently checked, revised, and finalized all mathematical statements, proofs, and references, and bears full responsibility for the content of the paper.

\appendix

\section{Well-posedness and continuity of the controlled equation}
\label{app:wellposedness}
\begin{lemma}[Helly--Bray for cumulative signed measures]
\label{lem:hellybray}
Let \(E\) be a compact metric space, and let
\(\nu_m,\nu\) be finite signed measures on \([0,1]\times E\) such that
\[
\nu_m\Rightarrow\nu,
\qquad
\sup_m\|\nu_m\|_{\rm TV}<\infty.
\]
Define
\[
H_m(t):=\nu_m([0,t]\times E),
\qquad
H(t):=\nu([0,t]\times E).
\]
Then
\[
H_m{\rightarrow} H
\quad\text{in }L^1([0,1]),
\qquad
H_m(1){\rightarrow} H(1).
\]
\end{lemma}

\begin{proof}
Testing \(\nu_m\Rightarrow\nu\) with the constant function
\(1\) gives \(H_m(1)\to H(1)\). For \(g\in C^1([0,1])\), Fubini's theorem
gives
\[
\int_{[0,1]\times E}g(s)\,\nu_m(\dd s,\dd\xi)
=g(1)H_m(1)-\int_0^1g'(t)H_m(t)\,\dd t,
\]
and similarly for \(\nu\). Hence
\(\int_0^1\varphi H_m\,\dd t\to\int_0^1\varphi H\,\dd t\) for every
\(\varphi\in C([0,1])\), taking \(\varphi=g'\).

Let \(V:=\sup_m\|\nu_m\|_{\rm TV}\). Each \(H_m\) is bounded
by \(V\) and has total variation at most \(V\) on \([0,1]\). By Helly's
selection theorem, every subsequence of \((H_m)\) has a further
subsequence that converges pointwise on \([0,1]\) to some function \(G\).
By dominated convergence, this subsequence converges to \(G\) in
\(L^1([0,1])\), and \(\int_0^1\varphi G\,\dd t=\int_0^1\varphi H\,\dd t\)
for every \(\varphi\in C([0,1])\). Hence \(G=H\) almost everywhere. Since
the limit does not depend on the subsequence, \(H_m\to H\) in
\(L^1([0,1])\).
\end{proof}

Uniform convergence may fail, even when the time marginal of
\(\nu\) is atomless. For instance, with \(\xi_0\in E\), the measures
\(\nu_m:=\delta_{(1/2+1/m,\xi_0)}-\delta_{(1/2,\xi_0)}\), \(m\geq2\), converge weakly to
\(\nu=0\), while \(H_m=-1\) on \([1/2,1/2+1/m)\).
\begin{proposition}[Well-posedness and continuity of the controlled equation]
\label{prop:wellposed}
For every \((A,B)\in\mathcal Q_\sigma\) and every list
\(\mathfrak m=((t_j,\xi_j,\delta_j))_{j=1}^N\) with
\(0\leq t_1\leq\cdots\leq t_N\leq1\), \(\xi_j\in K_\sigma\) and
\(\delta_j\in\R\),
equation \eqref{eq:controlled-error} has a unique c\`adl\`ag solution,
denoted by \(d_{A,B,\mathfrak m}\).
Let \((A_k,B_k)\in\mathcal Q_\sigma\) and
\(\mathfrak m_k=((t_{k,j},\xi_{k,j},\delta_{k,j}))_{j=1}^N\), with the same
\(N\), satisfy \(A_k\Rightarrow A\), \(B_k\Rightarrow B\), and
\((t_{k,j},\xi_{k,j},\delta_{k,j})\to(t_j,\xi_j,\delta_j)\) for each \(j\).
Then
\[
d_{A_k,B_k,\mathfrak m_k}
{\rightarrow}
d_{A,B,\mathfrak m}
\quad\text{in }L^1([0,1]),
\qquad
d_{A_k,B_k,\mathfrak m_k}(1)
{\rightarrow}
d_{A,B,\mathfrak m}(1).
\]
In particular, the map
\((A,B,\mathfrak m,\xi_f)\mapsto
\bar{\mathfrak g}(\xi_f,d_{A,B,\mathfrak m}(1))\) is continuous on
\(\mathcal C_\varepsilon\) for every \(0<\varepsilon\leq1\).
\end{proposition}

\begin{proof}
The equation is a Volterra equation whose drift is
Lipschitz in the state and is integrated against \(A\), which has
Lebesgue time marginal, while all jumps come from a forcing term of finite
variation. Step~1 proves existence and uniqueness by Picard iteration.
Step~2 proves continuity. The forcing terms converge by
Lemma~\ref{lem:hellybray}, the drift terms converge after approximating
the limit path by continuous functions, and Gronwall's inequality
combines the two. Let \(L_\beta\) be the Lipschitz constant of
\(\bar\beta\) in its second argument from the proof of
Lemma~\ref{lem:compactification}.

\medskip\noindent\textit{Step 1.}
Let
\[
H_{B,\mathfrak m}(t)
:=
\frac12
\int_{[0,t]\times K_\sigma}
\bar\alpha(\xi)\,B(\dd r,\dd\xi)
+
\sum_{j:t_j\leq t}
\bigl[
\delta_j-\bar{\mathfrak j}(\xi_j,\delta_j)
\bigr].
\]
This function is c\`adl\`ag and bounded, since it is the
distribution function of a finite signed measure on \([0,1]\). Then
\eqref{eq:controlled-error} reads \(d=\Phi(d)\), where, for
bounded Borel functions \(d\) on \([0,1]\),
\[
\Phi(d)(t)
:=
H_{B,\mathfrak m}(t)
+
\int_{[0,t]\times K_\sigma}
\bigl[
\bar\beta(\xi,d(r))-\bar\gamma(\xi)
\bigr]
A(\dd r,\dd\xi).
\]
The integrand is bounded and Borel, since \(\bar\beta\) is
jointly continuous. Since the time marginal of \(A\) is Lebesgue measure,
the integral term is Lipschitz in \(t\), with constant
\(\|\bar\beta\|_\infty+\|\bar\gamma\|_\infty\), and
\[
|\Phi(d_1)(t)-\Phi(d_2)(t)|
\leq
L_\beta\int_0^t|d_1(r)-d_2(r)|\,\dd r,
\qquad t\in[0,1].
\]
By induction, the \(\ell\)-th iterate \(\Phi^\ell\) is
Lipschitz in the uniform norm with constant \(L_\beta^\ell/\ell!\). Hence
the Picard iterates \(\Phi^\ell(H_{B,\mathfrak m})\) converge uniformly to a
fixed point \(d\), which is c\`adl\`ag because \(H_{B,\mathfrak m}\) is
c\`adl\`ag and the integral term is continuous. If \(d_1,d_2\) are two
solutions, the displayed bound and Gronwall's inequality give
\(d_1=d_2\). Since the integral term is continuous, the jumps of \(d\) are
those of \(H_{B,\mathfrak m}\), coming from atoms of the time marginal of
\(B\) and from the marks.

\medskip\noindent\textit{Step 2.}
Let
\(d_k:=d_{A_k,B_k,\mathfrak m_k}\), \(d:=d_{A,B,\mathfrak m}\),
\(H_k:=H_{B_k,\mathfrak m_k}\) and
\(H:=H_{B,\mathfrak m}\). We first show that
\begin{equation}
H_k{\rightarrow} H
\quad\text{in }L^1([0,1]),
\qquad
H_k(1){\rightarrow} H(1).
\label{eq:controlled-forcing-convergence}
\end{equation}
The function \(H_k\) is the distribution function of the
finite signed measure \(\nu_k\) on \([0,1]\) with
\[
\int_{[0,1]}f\,\dd\nu_k
=\frac12\int_{K_\sigma^0}f(t)\bar\alpha(\xi)\,B_k(\dd t,\dd\xi)
+\sum_{j=1}^Nf(t_{k,j})
\bigl[\delta_{k,j}-\bar{\mathfrak j}(\xi_{k,j},\delta_{k,j})\bigr],
\]
and \(H\) is that of the analogous measure \(\nu\). For
\(f\in C([0,1])\), the function \(f(t)\bar\alpha(\xi)\) is continuous on
\(K_\sigma^0\), and \(\bar{\mathfrak j}\) is jointly continuous by
Lemma~\ref{lem:compactification}. Hence \(\nu_k\Rightarrow\nu\). Moreover,
\(\|\nu_k\|_{\rm TV}\leq\frac12\|\bar\alpha\|_\infty B_k(K_\sigma^0)
+\sum_j|\delta_{k,j}-\bar{\mathfrak j}(\xi_{k,j},\delta_{k,j})|\) is
bounded in \(k\), because \(B_k(K_\sigma^0)\to B(K_\sigma^0)\) and the
marks converge. Lemma~\ref{lem:hellybray}, with \(E\) a single point,
gives \eqref{eq:controlled-forcing-convergence}.

We next compare the drift terms along the limit path \(d\),
\[
G_k(t)
:=
\int_{[0,t]\times K_\sigma}
\bigl[
\bar\beta(\xi,d(r))-\bar\gamma(\xi)
\bigr]
A_k(\dd r,\dd\xi),
\]
and let \(G\) be defined with \(A\) in
place of \(A_k\). We claim that
\begin{equation}
G_k{\rightarrow} G
\quad\text{in }L^1([0,1]),
\qquad
G_k(1){\rightarrow} G(1).
\label{eq:controlled-drift-convergence}
\end{equation}
Since \(d\) is only c\`adl\`ag, the weak convergence
\(A_k\Rightarrow A\) cannot be applied directly. Let \(g\in C([0,1])\),
and define \(G_k^g\) and \(G^g\) as \(G_k\) and \(G\) with \(d\) replaced by
\(g\). The integrand \(\bar\beta(\xi,g(r))-\bar\gamma(\xi)\) is continuous
and bounded on \(K_\sigma^0\). Hence the signed measures
\([\bar\beta(\xi,g(r))-\bar\gamma(\xi)]A_k(\dd r,\dd\xi)\) converge weakly,
with uniformly bounded total variation, and Lemma~\ref{lem:hellybray}
gives
\begin{equation}
G_k^g{\rightarrow} G^g
\quad\text{in }L^1([0,1]),
\qquad
G_k^g(1){\rightarrow} G^g(1).
\label{eq:controlled-drift-continuous-approx}
\end{equation}
Since the time marginals of \(A_k\) and \(A\) are Lebesgue
measure,
\[
\sup_{0\leq t\leq1}|G_k(t)-G_k^g(t)|
+
\sup_{0\leq t\leq1}|G(t)-G^g(t)|
\leq
2L_\beta\|d-g\|_{L^1([0,1])}.
\]
Since \(d\) is bounded and Borel, \(\|d-g\|_{L^1}\) can be
made arbitrarily small. Letting first \(k\to\infty\) and then \(g\to d\)
in \(L^1\) proves \eqref{eq:controlled-drift-convergence}.

Finally, subtracting the equations for \(d_k\) and \(d\),
and using the Lipschitz bound for \(\bar\beta\), we obtain
\[
|d_k(t)-d(t)|
\leq
a_k(t)
+
L_\beta\int_0^t|d_k(s)-d(s)|\,\dd s,
\qquad
a_k:=|H_k-H|+|G_k-G|.
\]
Gronwall's inequality gives
\(|d_k(t)-d(t)|\leq a_k(t)+L_\beta e^{L_\beta}\int_0^1a_k(s)\,\dd s\). By
\eqref{eq:controlled-forcing-convergence} and
\eqref{eq:controlled-drift-convergence}, \(\|a_k\|_{L^1}\to0\) and
\(a_k(1)\to0\). Hence \(d_k\to d\) in \(L^1([0,1])\) and
\(d_k(1)\to d(1)\).

For the last assertion, suppose that
\(\zeta_k=(A_k,B_k,\mathfrak m_k,\xi_{f,k})\) converges in \(\mathcal C_\varepsilon\)
to \(\zeta=(A,B,\mathfrak m,\xi_f)\in\mathcal C_{\varepsilon,N}\).
Then \(\zeta_k\in\mathcal C_{\varepsilon,N}\) for all sufficiently large \(k\), with \(A_k\Rightarrow A\), \(B_k\Rightarrow B\), convergence of
the marks, and \(\xi_{f,k}\to\xi_f\). By what we have shown,
\(d_k(1)\to d(1)\), and the joint continuity of \(\bar{\mathfrak g}\) in
Lemma~\ref{lem:compactification} gives
\(\bar{\mathfrak g}(\xi_{f,k},d_k(1))\to\bar{\mathfrak g}(\xi_f,d(1))\).
\end{proof}

\section[Proof of the comparisons in Lemma 4.6]{Proof of the comparisons in Lemma~\ref{lem:recovery-tilted-law}}
\label{app:euler-comparison}
We work in the setting of Step~4 of the proof of
Lemma~\ref{lem:recovery-tilted-law}, with \(\eta\), \(K\) and the data of
Lemma~\ref{lem:marked-approximation} fixed, and with constants \(C\)
independent of \(n\). Section~\ref{app:interval-comparison} compares the
tilted scheme with the restoring diffusion between prescribed indices.
Section~\ref{app:state-comparison} propagates this comparison through the
prescribed updates, which proves the state comparison
\eqref{eq:joint-coupled-state-comparison}, and
Section~\ref{app:measure-comparison} proves the measure comparison
\eqref{eq:joint-occupation-comparison}.

\subsection{An interval comparison}
\label{app:interval-comparison}
Between prescribed indices, the tilted scheme is an Euler scheme for
\eqref{eq:restoring-sde} up to a small remainder, so it stays close to the
diffusion on each such interval.

Let \(\mathcal E_n^{\mathrm{part}}\) be the union of the
grid cells \([k/n,(k+1)/n)\) whose interior contains one of
\(\tau_1,\ldots,\tau_{P-1}\). Then
\(|\mathcal E_n^{\mathrm{part}}|\leq(P-1)/n\), and outside
\(\mathcal E_n^{\mathrm{part}}\) the functions \(q^\eta,u^\eta\) agree
with their values at the left grid point. We claim the following. Let
\(0\leq a<b\leq n\) be integers with
\(\{a,\ldots,b-1\}\subset\mathcal K_n^0\), and let \(\widehat X^n\) solve
\begin{equation*}
\widehat X_t^n=\widehat X_{a/n}^n
+\int_{a/n}^t\psi_K(q^\eta(s)-\widehat X_s^n)\,\dd s
+\int_{a/n}^t\sqrt{u^\eta(s)}\,\dd\widetilde W_s^n,
\qquad t\in[a/n,b/n],
\end{equation*}
where \(\widehat X_{a/n}^n\) is square integrable and
\(\sigma(\widetilde W_s^n:s\leq a/n)\)-measurable. Then, with \(C\)
independent of \(a\), \(b\) and \(\widehat X_{a/n}^n\),
\begin{equation}
\mathbb E_{Q_n^{\eta,K}}
\sup_{a/n\leq t\leq b/n}|Z^n(t)-\widehat X_t^n|^2
\leq C\mathbb E_{Q_n^{\eta,K}}|Z_a^n-\widehat X_{a/n}^n|^2+Cn^{-1/2}.
\label{eq:joint-block-comparison}
\end{equation}

We first bound the remainder by which the tilted scheme differs from the
Euler scheme for \eqref{eq:restoring-sde}, then compare the tilted scheme
with the Euler scheme driven by the same Gaussian variables, and finally
compare that Euler scheme with the diffusion.

For \(k\in\mathcal K_n^0\), substituting
\eqref{eq:joint-normal-variables} into
\eqref{eq:lamperti-one-step-remainder} and using
\eqref{eq:recovery-mean} gives
\begin{equation}
Z_{k+1}^n-Z_k^n
=\frac1n\psi_K(q^\eta(k/n)-Z_k^n)
+\sqrt{\frac{u^\eta(k/n)}n}\,G_k+\epsilon_{k,n},
\label{eq:joint-discrete-expansion}
\end{equation}
where
\begin{equation}
\begin{aligned}
\epsilon_{k,n}
={}&-\frac{\alpha(Z_k^n)u^\eta(k/n)}{2n}(G_k^2-1)
-\frac{\alpha(Z_k^n)m_{k,n}^{\eta,K}\sqrt{u^\eta(k/n)}}{n^{3/2}}G_k\\
&-\frac{\alpha(Z_k^n)(m_{k,n}^{\eta,K})^2}{2n^2}
+r_n(Z_k^n,\Delta_kW).
\end{aligned}
\label{eq:joint-local-remainder}
\end{equation}
By \eqref{eq:recovery-mean-bound} and
\eqref{eq:recovery-conditional-law},
\(\mathbb E_{Q_n^{\eta,K}}[|\Delta_kW|^p\mid\mathcal F_{k/n}]\leq
C n^{-p/2}\) for \(k\in\mathcal K_n^0\) and \(1\leq p\leq6\). By
\eqref{eq:lamperti-remainder-global}, this gives
\(\mathbb E_{Q_n^{\eta,K}}[|r_n(Z_k^n,\Delta_kW)|\mid\mathcal F_{k/n}]\leq
Cn^{-3/2}\) and
\(\mathbb E_{Q_n^{\eta,K}}[|r_n(Z_k^n,\Delta_kW)|^2\mid\mathcal F_{k/n}]
\leq Cn^{-3}\). Since \(G_k\) and \(G_k^2-1\) are conditionally centered,
\eqref{eq:joint-local-remainder} and \eqref{eq:recovery-mean-bound} then
give
\begin{equation}
\left|\mathbb E_{Q_n^{\eta,K}}
[\epsilon_{k,n}\mid\mathcal F_{{k/n}}]\right|
\leq Cn^{-3/2},
\qquad
\mathbb E_{Q_n^{\eta,K}}
[|\epsilon_{k,n}|^2\mid\mathcal F_{{k/n}}]
\leq Cn^{-2},
\qquad k\in\mathcal K_n^0.
\label{eq:joint-remainder-moments}
\end{equation}
Set \(\epsilon_{k,n}=0\) for \(k\in\mathcal K_n^1\).
Decomposing each remainder into its conditional mean and a
martingale difference, Doob's inequality and
\eqref{eq:joint-remainder-moments} imply
\begin{equation}
\mathbb E_{Q_n^{\eta,K}}
\max_{0\leq\ell\leq n}
\bigg|\sum_{k<\ell}\epsilon_{k,n}\bigg|^2
\leq
2\bigg(\sum_{k=0}^{n-1}Cn^{-3/2}\bigg)^2
+8\sum_{k=0}^{n-1}\mathbb E_{Q_n^{\eta,K}}|\epsilon_{k,n}|^2
\leq\frac Cn.
\label{eq:joint-cumulative-remainder}
\end{equation}
All grid states are square integrable, since
\(|\mathcal S_{n,x}(z)-z|\leq C(|x|+n^{-1})\) by
\eqref{eq:transformed-euler-update} and the Gaussian laws in
\eqref{eq:recovery-conditional-law} have finite moments.

To prove \eqref{eq:joint-block-comparison}, we pass through
the Euler scheme driven by the same Gaussian variables,
\begin{equation}
S_a:=Z_a^n,\qquad
S_{k+1}:=S_k+\frac1n\psi_K(q^\eta({k/n})-S_k)
+\sqrt{\frac{{u^\eta(k/n)}}{n}}G_k,
\qquad a\leq k<b,
\label{eq:joint-auxiliary-euler}
\end{equation}
and its continuous interpolation
\begin{equation}
\bar S_t:=S_k+(t-{k/n})\psi_K(q^\eta({k/n})-S_k)
+\sqrt{{u^\eta(k/n)}}(\widetilde W_t^n-\widetilde W_{{k/n}}^n),
\quad t\in[{k/n},{(k+1)/n}],\quad a\leq k<b.
\label{eq:joint-euler-interpolation}
\end{equation}
First, subtracting \eqref{eq:joint-auxiliary-euler} from
\eqref{eq:joint-discrete-expansion} and using the \(K\)-Lipschitz
continuity of \(\psi_K\), the discrete Gronwall inequality and
\eqref{eq:joint-cumulative-remainder} give
\begin{equation}
\mathbb E_{Q_n^{\eta,K}}\max_{a\leq\ell\leq b}|Z_\ell^n-S_\ell|^2
\leq e^{2K}\mathbb E_{Q_n^{\eta,K}}
\max_{a\leq\ell\leq b}
\bigg|\sum_{k=a}^{\ell-1}\epsilon_{k,n}\bigg|^2
\leq\frac Cn.
\label{eq:joint-euler-grid-comparison}
\end{equation}
Second, \(\bar S_{k/n}=S_k\), and \(|\psi_K|\leq K\),
\(u^\eta\leq R\) and Doob's \(L^4\) inequality give
\begin{equation}
\sup_{{k/n}\leq t\leq {(k+1)/n}}
\mathbb E_{Q_n^{\eta,K}}|\bar S_t-S_k|^2\leq Cn^{-1},
\qquad
\mathbb E_{Q_n^{\eta,K}}\max_{a\leq k<b}
\sup_{{k/n}\leq t\leq {(k+1)/n}}|\bar S_t-S_k|^2
\leq Cn^{-1/2},
\label{eq:joint-local-interpolation-moments}
\end{equation}
where the second bound follows from the fourth moments
\(\mathbb E\sup_t|\bar S_t-S_k|^4\leq Cn^{-2}\) and
\(\mathbb E\max_k(\cdot)^2\leq(\sum_k\mathbb E(\cdot)^4)^{1/2}\). Third,
\(\bar S-\widehat X^n\) has drift
\(\psi_K(q^\eta(k/n)-S_k)-\psi_K(q^\eta(s)-\widehat X_s^n)\) and diffusion
coefficient \(\sqrt{u^\eta(k/n)}-\sqrt{u^\eta(s)}\) on
\([k/n,(k+1)/n)\). The drift is bounded by
\(K|S_k-\bar S_s|+K|\bar S_s-\widehat X_s^n|\) outside
\(\mathcal E_n^{\mathrm{part}}\) and by \(2K\) on it, and the diffusion
coefficient vanishes outside \(\mathcal E_n^{\mathrm{part}}\) and is
bounded by \(\sqrt R\). Hence the first bound in
\eqref{eq:joint-local-interpolation-moments}, It\^o's isometry and Doob's
inequality give, for \(t\in[a/n,b/n]\),
\begin{equation}
\begin{aligned}
\mathbb E_{Q_n^{\eta,K}}
\sup_{a/n\leq r\leq t}|\bar S_r-\widehat X_r^n|^2
\leq{}&C\mathbb E_{Q_n^{\eta,K}}|Z_a^n-\widehat X_{a/n}^n|^2+\frac Cn\\
&+C\int_{a/n}^t\mathbb E_{Q_n^{\eta,K}}
\sup_{a/n\leq r\leq s}|\bar S_r-\widehat X_r^n|^2\,\dd s.
\end{aligned}
\label{eq:joint-euler-diffusion-inequality}
\end{equation}
For \(t\in[k/n,(k+1)/n)\),
\(Z^n(t)-\widehat X_t^n=(Z_k^n-S_k)+(S_k-\bar S_t)
+(\bar S_t-\widehat X_t^n)\). Gronwall's inequality in
\eqref{eq:joint-euler-diffusion-inequality},
\eqref{eq:joint-euler-grid-comparison} and the second bound in
\eqref{eq:joint-local-interpolation-moments} therefore prove
\eqref{eq:joint-block-comparison}.

\subsection{Propagation through the prescribed updates}
\label{app:state-comparison}
At a prescribed index, both processes perform the update
\(\mathcal S_{\delta_j}\) up to small errors. We now propagate \eqref{eq:joint-block-comparison} through
the prescribed updates. The diffusion jumps at \(s_j\), whereas the scheme
completes the corresponding update at \(r_j^n\). Set
\begin{equation}
D_{0,n}:=0,\qquad
D_{j,n}:=
\mathbb E_{Q_n^{\eta,K}}
|Z_{n r_j^n}^n-\Gamma_{r_j^n}^{\eta,K,n}|^2,
\qquad j=1,\ldots,N.
\label{eq:joint-restarting-error-definition}
\end{equation}
For \(j=1,\ldots,N+1\), with \(r_{N+1}^n:=1\), let \(X_j^n\)
be the solution of \eqref{eq:restoring-sde} without updates, driven by
\(\widetilde W^n\) and started at time \(r_{j-1}^n\),
\begin{equation}
X_j^n(t)
=\Gamma_{r_{j-1}^n}^{\eta,K,n}
+\int_{r_{j-1}^n}^t
\psi_K(q^\eta(r)-X_j^n(r))\,\dd r+\int_{r_{j-1}^n}^t
\sqrt{u^\eta(r)}\,\dd\widetilde W_r^n,
\qquad t\in[r_{j-1}^n,r_j^n].
\label{eq:joint-continuous-extension}
\end{equation}
Since \(\Gamma^{\eta,K,n}\) has no update in
\((r_{j-1}^n,s_j)\), pathwise uniqueness gives \(X_j^n=\Gamma^{\eta,K,n}\)
on \([r_{j-1}^n,s_j)\), and continuity gives
\(X_j^n(s_j)=\Gamma_{s_j-}^{\eta,K,n}\). Since \(r_j^n-s_j\leq2/n\),
\(|\psi_K|\leq K\), \(u^\eta\leq R\) and Doob's inequality give
\begin{equation}
\mathbb E_{Q_n^{\eta,K}}\bigg[
\sup_{s_j\leq t\leq r_j^n}
|X_j^n(t)-X_j^n(s_j)|^2+\sup_{s_j\leq t\leq r_j^n}
|\Gamma_t^{\eta,K,n}-\Gamma_{s_j}^{\eta,K,n}|^2
\bigg]\leq\frac Cn,
\qquad j=1,\ldots,N.
\label{eq:joint-short-interval-comparison}
\end{equation}
Fix \(j\in\{1,\ldots,N+1\}\). The indices
\(n r_{j-1}^n,\ldots,n p_j^n-1\) belong to \(\mathcal K_n^0\), and
\(\Gamma_{r_{j-1}^n}^{\eta,K,n}\) is measurable with respect to
\(\sigma(\widetilde W_s^n:s\leq r_{j-1}^n)\). Hence
\eqref{eq:joint-block-comparison} with \(\widehat X^n=X_j^n\) gives
\begin{equation}
\mathbb E_{Q_n^{\eta,K}}
\sup_{r_{j-1}^n\leq t\leq p_j^n}
|Z^n(t)-X_j^n(t)|^2
\leq C(D_{j-1,n}+n^{-1/2}).
\label{eq:joint-before-update-comparison0}
\end{equation}
Since \(Z^n(p_j^n)=Z_{np_j^n}^n\),
\(X_j^n(s_j)=\Gamma_{s_j-}^{\eta,K,n}\) and \(s_j\leq p_j^n\leq r_j^n\),
\eqref{eq:joint-before-update-comparison0} and
\eqref{eq:joint-short-interval-comparison} give
\begin{equation}
\mathbb E_{Q_n^{\eta,K}}
|Z_{n p_j^n}^n-\Gamma_{s_j-}^{\eta,K,n}|^2
\leq C(D_{j-1,n}+n^{-1/2}),
\qquad j=1,\ldots,N+1.
\label{eq:joint-before-update-comparison}
\end{equation}
For \(j=N+1\), we have \(p_{N+1}^n=s_{N+1}=1\), so only
\eqref{eq:joint-before-update-comparison0} is needed, and the left-hand
side is \(\mathbb E_{Q_n^{\eta,K}}|Z_n^n-\Gamma_1^{\eta,K,n}|^2\).

Now let \(j\leq N\). The two updates are
\(\Gamma_{s_j}^{\eta,K,n}
=\mathcal S_{\delta_j}(\Gamma_{s_j-}^{\eta,K,n})\)
and
\(Z_{n r_j^n}^n
=\mathcal S_{n,\Delta_{n p_j^n}W}(Z_{n p_j^n}^n)\).
By \eqref{eq:lamperti-drift-perturbation},
\eqref{eq:beta3}, and since the derivative of \(\mathcal S_{\delta_j}\) is bounded
by \(L_j:=(\sigma_{\mathrm{max}}/\sigma_{\mathrm{min}})
(1+\|\sigma'\|_\infty|\delta_j|)\),
\[
|\mathcal S_{n,x}(z)-\mathcal S_{\delta_j}(y)|
\leq\frac Cn+\frac{\sigma_{\mathrm{max}}}{\sigma_{\mathrm{min}}}|x-\delta_j|
+L_j|z-y|,
\qquad x,y,z\in\R.
\]
Hence \eqref{eq:joint-before-update-comparison} and
\eqref{eq:joint-prescribed-increment-L2} give
\(\mathbb E_{Q_n^{\eta,K}}|Z_{n r_j^n}^n-\Gamma_{s_j}^{\eta,K,n}|^2
\leq C(D_{j-1,n}+n^{-1/2})\). Combined with
\eqref{eq:joint-short-interval-comparison} for
\(\Gamma_{r_j^n}^{\eta,K,n}-\Gamma_{s_j}^{\eta,K,n}\), this yields
\begin{equation}
D_{j,n}\leq C(D_{j-1,n}+n^{-1/2}),
\qquad j=1,\ldots,N.
\label{eq:joint-after-update-recursion}
\end{equation}
Since \(D_{0,n}=0\) and \(N\) is fixed, iteration gives
\(\max_{0\leq j\leq N}D_{j,n}\leq Cn^{-1/2}\). Finally,
\([0,1]\setminus E_n\) is the union of the intervals
\([r_{j-1}^n,s_j)\), \(j=1,\ldots,N+1\), and the point \(1\). On
\([r_{j-1}^n,s_j)\), we have \(X_j^n=\Gamma^{\eta,K,n}\) and
\(s_j\leq p_j^n\). Therefore
\eqref{eq:joint-before-update-comparison0} and
\eqref{eq:joint-before-update-comparison}, together with the bound on
\(D_{j,n}\), prove \eqref{eq:joint-coupled-state-comparison}.

\subsection{The measure comparison}
\label{app:measure-comparison}
The measure \(B_n^\circ\) records the squared increments on the ordinary
grid steps \(k\in\mathcal K_n^0\). To compare it with the continuous-time
variance measure \(B_{\Gamma,n}\) in \eqref{eq:joint-auxiliary-measures},
we introduce the intermediate measure
\[
\overline B_n:=\mathbf1_{\{t\notin\bigcup_{j}[p_j^n,r_j^n)\}}
\,u^\eta(t)\,\dd t\,\delta_{\iota(Z^n(t))}.
\]

Fix \(f\in C(K_\sigma^0)\). By
\eqref{eq:joint-coupled-state-comparison}, Markov's inequality, the
tightness noted in Step~4 of the proof of Lemma~\ref{lem:recovery-tilted-law}, and \eqref{eq:joint-localization} with
\(g(t,x)=f(t,\iota(x))\) and \(T_n=[0,1]\setminus E_n\),
\begin{equation}
\sup_{t\in[0,1]\setminus E_n}
\bigl|f(t,\iota(Z^n(t)))-f(t,\iota(\Gamma_t^{\eta,K,n}))\bigr|
\longrightarrow0
\quad\text{in }Q_n^{\eta,K}\text{-probability}.
\label{eq:joint-uniform-test-comparison}
\end{equation}
Since \(u^\eta\leq R\), \(|E_n|\leq2N/n\), and each
\([p_j^n,r_j^n)\) is contained in \(E_n\), the differences
\(|\int f\,\dd A_n-\int f\,\dd A_{\Gamma,n}|\) and
\(|\int f\,\dd\overline B_n-\int f\,\dd B_{\Gamma,n}|\) are bounded by
\(R+1\) times the supremum in \eqref{eq:joint-uniform-test-comparison}
plus \(4(R+1)N\|f\|_\infty/n\). Applying this to the functions in
\eqref{eq:dw} and controlling the tail of the series gives
\[
d_{\mathrm w}(A_n,A_{\Gamma,n})
+d_{\mathrm w}(\overline B_n,B_{\Gamma,n})\longrightarrow0
\quad\text{in }Q_n^{\eta,K}\text{-probability}.
\]

It remains to compare \(B_n^\circ\) with \(\overline B_n\).
For \(f\in C(K_\sigma^0)\) and \(k\in\mathcal K_n^0\), set
\begin{equation}
V_{k,n}(f):=\bigg[
(\Delta_kW)^2-\frac{{u^\eta(k/n)}}{n}
-\frac{(m_{k,n}^{\eta,K})^2}{n^2}
\bigg]f({k/n},\iota(Z_k^n)).
\label{eq:joint-variance-difference}
\end{equation}
By \eqref{eq:joint-auxiliary-control} and the definition of
\(\overline B_n\),
\begin{equation}
\int f\,\dd B_n^\circ-\int f\,\dd\overline B_n
=\sum_{k\in\mathcal K_n^0}V_{k,n}(f)
+\frac1{n^2}\sum_{k\in\mathcal K_n^0}
(m_{k,n}^{\eta,K})^2 f({k/n},\iota(Z_k^n))+R_n(f),
\label{eq:joint-variance-decomposition}
\end{equation}
where
\begin{equation}
R_n(f):=\sum_{k\in\mathcal K_n^0}\int_{{k/n}}^{{(k+1)/n}}
\left[{u^\eta(k/n)} f({k/n},\iota(Z_k^n))
-u^\eta(t)f(t,\iota(Z_k^n))\right]\,\dd t.
\label{eq:joint-variance-riemann-remainder}
\end{equation}
The second sum in \eqref{eq:joint-variance-decomposition} is
at most \(B_K^2\|f\|_\infty/n\) by \eqref{eq:recovery-mean-bound}. Let
\(\omega_f(\theta):=\sup_{|t-s|\leq \theta,\,\xi\in K_\sigma}
|f(t,\xi)-f(s,\xi)|\), which tends to zero as \(\theta\downarrow0\) by uniform
continuity of \(f\). For \(t\in[t_k,t_{k+1})\setminus\mathcal E_n^{\mathrm{part}}\),
we have \(u^\eta(t)=u^\eta(t_k)\). Since \(|\mathcal E_n^{\mathrm{part}}|\leq(P-1)/n\),
\eqref{eq:joint-variance-riemann-remainder} gives
\[
|R_n(f)|\leq R\omega_f(1/n)
+\frac{2R(P-1)\|f\|_\infty}{n}.
\]
By \eqref{eq:recovery-conditional-law}, the variables
\(V_{k,n}(f)\) are martingale differences with
\[
\mathbb E_{Q_n^{\eta,K}}[|V_{k,n}(f)|^2\mid\mathcal F_{{k/n}}]
=f({k/n},\iota(Z_k^n))^2
\bigg[\frac{2({u^\eta(k/n)})^2}{n^2}
+\frac{4{u^\eta(k/n)}(m_{k,n}^{\eta,K})^2}{n^3}\bigg]
\leq\frac{C\|f\|_\infty^2}{n^2},
\]
so that
\(\mathbb E_{Q_n^{\eta,K}}|\sum_{k\in\mathcal K_n^0}V_{k,n}(f)|^2\leq
C\|f\|_\infty^2/n\). Hence all three terms in
\eqref{eq:joint-variance-decomposition} tend to zero in
\(L^2(Q_n^{\eta,K})\). Applying this to the functions in \eqref{eq:dw} and
controlling the tail of the series gives
\(d_{\mathrm w}(B_n^\circ,\overline B_n)\to0\) in
\(Q_n^{\eta,K}\)-probability. The triangle inequality for
\(d_{\mathrm w}\) completes the proof of
\eqref{eq:joint-occupation-comparison}.

\end{document}